\documentclass[reqno,a4paper,11pt]{amsart}

\usepackage[all,poly]{xy}
\usepackage{amsfonts}
\usepackage[mathcal]{eucal}
\usepackage{eufrak}
\usepackage{amssymb}
\usepackage{amsmath}
\usepackage{mathrsfs}
\usepackage{color}
\usepackage[pagebackref,colorlinks]{hyperref}
\usepackage{enumerate}

\theoremstyle{plain}
\newtheorem {lemma}{Lemma}[section] % add section to get consequence numbering
\newtheorem {theorem}[lemma]{Theorem}

\newtheorem {corollary}[lemma]{Corollary}

\newtheorem {proposition}[lemma]{Proposition}

\theoremstyle{definition}
\newtheorem {remark}[lemma]{Remark}

\newtheorem {example}[lemma]{Example}

\theoremstyle{definition}

\newtheorem{deff}[lemma]{Definition}{}

\newcommand{\op}{\operatorname{op}}
\newcommand{\M}{\operatorname{\mathbb M}}
\newcommand{\LL}{\operatorname{\mathcal L}}
\newcommand{\I}{\operatorname{\mathcal I}}
\newcommand{\gr}{\operatorname{gr}}
\newcommand{\V}{\operatorname{\mathcal V}}

\newcommand{\Ga}{\Gamma}
\newcommand{\ga}{\gamma}

\newcommand{\de}{\delta}
\newcommand{\la}{\lambda}

\newcommand{\Aut}{\operatorname{Aut}}

\newcommand{\End}{\operatorname{End}}

\newcommand{\id}{\operatorname{id}}

\newcommand{\Esc}{\operatorname{Esc}}
\newcommand{\st}{\operatorname{st}}
\newcommand{\outdeg}{\operatorname{outdeg}}
\newcommand{\totdeg}{\operatorname{totdeg}}
\newcommand{\card}{\operatorname{card}}
\newcommand{\sink}{\operatorname{Sink}}
\newcommand{\coker}{\operatorname{coker}}

\newcommand{\A}{\underline{\alpha}}

\newcommand{\Pp}{\underline{p}}

\title{The graded structure of Leavitt Path algebras}

\begin{thanks}
{I would like to thank Gene Abrams for his illuminating remarks, Andrew Steele for his formulation of using exit edges rather than exit vertices in Definition~\ref{exiedg} and Gonzalo Aranda Pino for introducing me to this subject and for his invitation to Malaga in Summer 2009. Part of this work has been done in Summers of 2009 and 2011 at ICTP,
Trieste, Italy. The author acknowledges the support of EPSRC first grant
scheme EP/D03695X/1 and the standard scheme EP/I007784/1.}
\end{thanks}

\author{Roozbeh Hazrat}\address{
Department of Pure Mathematics\\
Queen's University\\
Belfast BT7 1NN\\
United Kingdom} \email{r.hazrat@qub.ac.uk}


\begin{document}

\begin{abstract}
A Leavitt path algebra associates to a directed graph a $\mathbb Z$-graded algebra and in its simplest form it recovers the Leavitt algebra $L(1,k)$. In this note, we first study this $\mathbb Z$-grading and characterize the ($\mathbb Z$-graded) structure of Leavitt path algebras, associated to finite acyclic graphs, $C_n$-comet, multi-headed graphs and a mixture of these graphs (i.e., polycephaly graphs). The last two types  are examples of graphs whose Leavitt path algebras are strongly graded.  
We give a criterion when a Leavitt path algebra is strongly graded and in particular characterize unital Leavitt path algebras which are strongly graded completely, along the way obtaining classes of algebras which are group rings or crossed-products.  In an attempt to generalize the grading, we introduce weighted Leavitt path algebras associated to directed weighted graphs which have natural $\textstyle{\bigoplus} \mathbb Z$-grading and in their simplest form recover  the Leavitt algebras $L(n,k)$. We then show that the basic properties of Leavitt path algebras can be naturally carried over to weighted Leavitt path algebras. 
\end{abstract}

\maketitle

%\tableofcontents

\section{Introduction} \label{introf}

A Leavitt path algebra (LPA for short), introduced by Abrams and Aranda Pino~\cite{aap05}, and Ara, Moreno and Pardo~\cite{amp}  
associates to a directed graph $E$ a $\mathbb Z$-graded algebra $\LL(E)$ which is equipped with an (anti-graded) involution. In its simplest form, when the graph $E$ has only one vertex and $k+1$ loops, $\LL(E)$ recovers  the algebra constructed by Leavitt~\cite[p.118]{vitt62} which is of type $(1,k)$.  
The characterization of  Leavitt path algebras (such as simplicity, finite dimensionality, locally finiteness, exchange, etc.) in terms of intrinsic properties of the underlying graph has been the subject of recent studies \cite{aap05,aap06,aap07,aalp,ap08,amp}. However, many of these characterizations have been carried out without taking into account the natural $\mathbb Z$-graded structure of these algebras, i.e., not considering $\LL$ as a functor from the category of graphs to graded algebras. 

In this note we are concerned with the graded structure of  Leavitt path algebras. The note is divided into two parts. In the first part, we study the natural $\mathbb Z$-grading of the Leavitt path algebras, characterizing when these algebras are strongly graded. The first main theorem of this part (Theorem~\ref{sthfin}) states that for a finite graph $E$, the Leavitt path algebra $\LL_R(E)$, with coefficients in any unital ring $R$,  is strongly graded if and only if any vertex is connected to a cycle. Two distinguished types of strongly graded Leavitt path algebras  arise from  $C_n$-comet graphs (see Definition~\ref{cometi}) and multi-headed rose graphs (see Fig.~\ref{monster2}).  Considering the graded structure of Leavitt path algebras reveals much more about the structure of these algebras. 
In~\cite{ajis}, Abrams, Aranda Pino and Siles Molina characterize  locally finite just infinite algebras (which turned out to be LPAs associated to $C_n$-comets).  In~\cite[Theorem~3.3]{ajis}, it is shown that if $E$ is a $C_n$-comet, with the unique cycle $C$, then $\LL_K(E)\cong \M_d(K[x,x^{-1}])$, where $K$ is a field, $d$ is the number of paths in $E$ which do not contain the cycle $C$ and which end in $v$, for a fixed vertex  $v$  in $C$.
However, in this characterization the natural $\mathbb Z$-grading of Leavitt path algebras is not taken into account. 
 For this reason for the following graphs
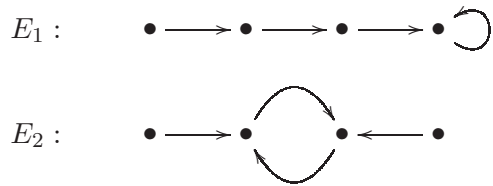
\begin{equation*}
\xymatrix{
E_1:  &  \bullet \ar[r] &   \bullet \ar[r] &   \bullet \ar[r] & \bullet \ar@(rd,ru) &  \\
E_2:  & \bullet \ar[r] & \bullet \ar@/^1.5pc/[r] & \bullet \ar@/^1.5pc/[l] & \bullet \ar[l]&
}
\end{equation*}
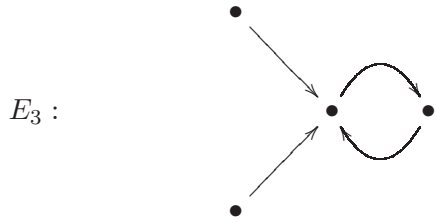
\smallskip
\begin{equation*}
\xymatrix{
& \bullet \ar[dr] &\\
E_3: \qquad \quad  &  & \bullet \ar@/^1.5pc/[r] & \bullet \ar@/^1.5pc/[l] &\\
& \bullet \ar[ur] &&\\
}
\end{equation*}
and 
\begin{equation*}
E_4: \qquad  
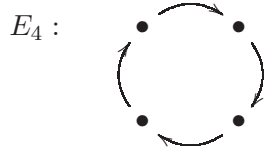
\xymatrix{
\bullet  \ar@/^0.75pc/[r] & \bullet \ar@/^0.75pc/[d] \\
\bullet \ar@/^0.75pc/[u] & \bullet \ar@/^0.75pc/[l]}
\end{equation*}
we obtain  
\[
\LL_K(E_1)\cong \LL_K(E_2) \cong \LL_K(E_3) \cong \LL_K(E_4) \cong \M_4(K[x,x^{-1}]).\]

In this paper, we shall see that, building on their approach, by taking into account the $\mathbb Z$-graded structure of the 
 LPAs, it turns out  that $\LL_K(E_1)$ is a group ring (Theorem~\ref{polyheadnew}), $\LL_K(E_2)$ and $\LL_K(E_4)$ are skew-group rings whereas $\LL_K(E_3)$ is not even a crossed-product (Examples~\ref{noncori} and ~\ref{cycskew}).

The second main theorem of this part, Theorem~\ref{polyhead},  characterizes the Leavitt path algebras associated to multi-headed graphs (see Figure~\ref{monster} and 
Definition~\ref{popyt}) and gives a criterion when  this type of graphs produces group rings 
(Theorem~\ref{polyheadnew}). As a corollary we obtain necessary and sufficient conditions for Leavitt path algebras of finite acyclic graphs and $C_n$-comets to be graded isomorphic (Theorems~\ref{acyclicc}, ~\ref{grcomet}, see also Example~\ref{niroi}). 

Whereas a Leavitt path algebra associates to a directed graph a $\mathbb Z$-graded algebra and it covers Leavitt's algebra of type $(1,k)$ when the underlying graph has one vertex and $k+1$ loops, a natural question  is whether one can associate to a directed graph an algebra with an arbitrary grading rather than $\mathbb Z$ grading, in  such a way that it recovers Leavitt path algebras when we restrict the grading to $\mathbb Z$, and also recovers Leavitt's algebra of type $(n,k)$ for a suitable graph. In the second part of this note, we do just this.  

In an attempt to define Leavitt path algebras with more general natural grading,  we define  the weighted Leavitt path algebras (wLPA for short) which are equipped with $\textstyle{\bigoplus} \mathbb Z$ grading (in fact we can define weighted Leavitt path algebras with $G$-grading where $G$ is any group, see Remark~\ref{ggrad}).  In the special case of a graph with  weights $1$ (or unweighted), this construction gives the usual Leavitt path algebras and in its simplest form, when the graph $E$ has only one vertex and $n+k$ loops of weights $n$,  the weighted Leavitt path algebra recovers the algebra constructed by Leavitt which is of type $(n,k)$ (see Example~\ref{lat}). The weighted Leavitt path algebras provide new classes of algebras which could not be obtained using unweighted graphs (i.e., using Leavitt path algebras). For example, note that except 
$\xymatrix{\LL_R(\!\!\!\!\!&\!\bullet\ar@(lu,ld)})=R[x,x^{-1}]$ 
(which happens to be the only commutative LPA along with $\LL_R(\bullet)=R$, where $R$ is an integral domain), all LPA have plenty of zero divisors. (In this note we only consider connected graphs, as a disjoint graph would simply produce direct sum of the corresponding Leavitt path algebras of its connected subgraphs.) However all wLPA with one vertex and with weights greater than one and less than the number of loops are non-commutative domains (see Example~\ref{lat}). We then establish basic properties of weighted Leavitt path algebras  in the remaining of this section. We shall see that, one can naturally adapt and re-write several of the theorems on the setting of Leavitt path algebras in the setting of weighted Leavitt path algebras (see in particular Theorems~\ref{wpolyhead},~\ref{kspcts},~\ref{wke}),  an indication that the weighted approach  could be the right generalization for these graph algebras.

\section{Preliminaries}

\subsection{Grading on rings and matrices} \label{gradingsec} 
A Leavitt path algebra has a natural $\mathbb Z$-graded structure which is the focus of this paper. In Section~\ref{secitwo} we define weighted Leavitt path algebras which have $\textstyle{\bigoplus} \mathbb Z$-grading. In fact we can define weighted Leavitt path algebras with $G$-grading, where $G$ is any arbitrary group (see Remark~\ref{ggrad}). Then setting $G=\mathbb Z$ and the weight map the constant map assigning $1$ to each edge we recover the usual Leavitt path algebra.

In this section we recall some basic definitions and (recent) results in the graded setting. For a graded ring $A$, the grading on a matrix ring of $A$ is also recalled and a theorem on classifying these gradings is proved (Theorem~\ref{clasigr}). This is needed in studying the grading on Leavitt path algebras of $C_n$-comets (see Definition~\ref{cometi}). Since we are ultimately dealing with $\mathbb Z$-gradings (and  $\textstyle{\bigoplus} \mathbb Z$-gradings),  we assume our grade groups are abelian, although all the concepts below can be arranged and defined for arbitrary (non-abelian) grade groups as well.  

A ring $A = \textstyle{\bigoplus_{ \ga \in \Ga}} A_{\ga}$ is called a
\emph{$\Ga$-graded ring}, or simply a \emph{graded ring},
\index{graded ring} if $\Ga$ is an (abelian) group, each $A_{\ga}$ is
an additive subgroup of $A$ and $A_{\ga}  A_{\delta} \subseteq
A_{\ga + \delta}$ for all $\ga, \delta \in \Ga$. %Note that, at this stage, we
%do not assume $\Ga$ to be abelian or totally ordered.
%Each $x \in R$ can be uniquely expressed as a finite sum $x=
%\sum_{\ga \in \Ga} x_{\ga}$ with each $x_{\ga} \in R_{\ga}$. 
The
elements of $A_\ga$ are called \emph{homogeneous of degree $\ga$}
and we write deg$(a) = \ga$ if $a \in A_{\ga}$. We let $A^{h} =
\bigcup_{\ga \in \Ga} A_{\ga}$ be the set of homogeneous elements of
$A$. 

%The set
%$$
%\Ga_{R} = \big \{ \ga \in \Ga : R_{\ga} \neq \{0 \} \big \},
%$$
%which is also denoted by $\Supp (R)$, is called the \emph{support}
%(or grade set) of $R$. We note that the support of $R$ is not
%necessarily a group.

A $\Ga$-graded ring $A=\textstyle{\bigoplus_{ \ga \in \Ga}} A_{\ga}$
is called a \emph{strongly graded ring} if $A_{\ga} A_{\de} = A_{\ga +\de}$
for all $\ga, \de \in \Ga$. A graded ring $A$ is called a
\emph{crossed-product} if there is an invertible element in every
homogeneous component $A_\ga$ of $A$; that is, $A^* \cap A_\ga \neq
\emptyset$ for all $\ga \in \Ga$ where $A^*$ is the group of all invertible elements of $A$. It is immediate that a crossed-product is a strongly graded ring.  The structure of crossed-product algebras are known (see~\cite[\S 1.4]{grrings}) and can be described as follows: for any $\ga \in \Gamma$ choose $u_\ga \in A^* \cap A_\ga$ and define $\sigma:\Gamma\rightarrow \Aut(A_0)$ by  $\sigma(\ga)(a)=u_\ga a u_\ga^{-1}$ for $\ga \in \Gamma$ and $a\in A_0$. Furthermore define the cocycle map $\alpha:\Gamma\times \Gamma \rightarrow A_0^*$ by $\alpha(\zeta,\eta)=u_\zeta u_\eta u_{\zeta \eta}^{-1}$. Then $A={A_0}_\alpha^\sigma[\Gamma]=\textstyle{\bigoplus_{\ga \in \Gamma}} A_0 \ga$ with multiplication $(a_1 \zeta)(a_2 \eta)=a_{1} {}^{\zeta} a_2 \alpha(\zeta,\eta) \zeta\eta$ where ${}^\zeta a$ is defined as $\sigma(\zeta)(a)$. Note that when $\Gamma$ is cyclic, which is the case in our paper,  one can choose $u_i=u_1^i$ for $u_1 \in A^* \cap A_1$ and thus the cocycle map $\alpha$ is trivial, $\sigma$ is a homomorphism and the crossed-product is a skew group ring denoted by $A_0 \star_\sigma \mathbb Z$. %In this case $R=\textstyle{\bigoplus_{\i \in \mathbb Z}} R_0 x^i$. 
Moreover if $u_1$ is in the center of $A$, then $\sigma$ is the identity  map and  the crossed-product reduces to the group ring $A_0[\Gamma]$. 

For a group $\Gamma$ and an arbitrary ring $A$ with identity, the group ring $A[\Gamma]$ has a natural $\Gamma$-grading 
$A[\Gamma]=\textstyle{\bigoplus_{\ga\in \Gamma}} A\ga$. If $\Gamma$ is abelian and  $A$ itself is a $\Gamma$-graded ring
$A=\textstyle{\bigoplus_{\ga\in \Gamma}} A_{\ga}$, then $A[\Gamma]$ has a $\Gamma$-grading 
\begin{equation}\label{hgogt}
A[\Gamma]=\textstyle{\bigoplus_{\ga\in \Gamma}} A^\ga, \text{  where  } A^\ga=\textstyle{\bigoplus_{\ga=\zeta+\zeta'}}  A_\zeta \zeta'.
\end{equation}
This grading will be used to describe the grading of Leavitt path algebras in Theorems~\ref{polyhead} and ~\ref{polyheadnew}. 

Let $A$ be a $\Gamma$-graded ring. A \emph{graded left $A$-module} $M$ is defined to be a left $A$-module
with a direct sum decomposition $M=\bigoplus_{\ga \in \Omega}
M_{\ga}$, where each $M_{\ga}$ is an additive subgroup of $M$ and
$\Ga \subseteq \Omega$, such that $A_{\ga} \cdot
M_{\la} \subseteq M_{\ga + \la}$ for all $\ga \in \Ga, \la \in
\Omega$. For some $\de \in \Omega$, we define the
$\de$-shifted $A$-module $M(\de)$ \label{deshiftedmodule} as $M(\de)
=\bigoplus_{\ga \in \Omega} M(\de)_\ga$ where $M(\de)_\ga =
M_{\ga+\de}$. 
For two graded $A$-modules $M$ and $N$, a {\it graded $A$-module homomorphism of degree $\delta$}, is an $A$-module homomorphism $f:M\rightarrow N$, such that $f(M_\ga)\subseteq N_{\ga+\delta}$ for any $\ga \in \Omega$.

To study the graded structure of Leavitt path algebras of finite acyclic graphs (see Definition~\ref{mulidef} and Theorem~\ref{acyclicc}), we need the concept of the group grading on matrix algebras. 
Given a group $\Omega$ and a ring $A$ (which is not graded, i.e., has a trivial grading), it is
possible to define a grading on $\M_n (A)$. Such a grading is
called a \emph{good grading} of $\M_n (A)$ if the matrices
$e_{ij}$ are homogeneous, where $e_{ij}$ is the matrix with $1$ in
the $ij$-position and zero elsewhere. These group gradings on matrix rings have been
studied by D\u{a}sc\u{a}lescu et al.~\cite{dascalescu}. Let $\{\delta_1,\dots,\delta_n\}$ be a subset of $\Omega$, and define a grading on $\M_n(A)$ by assigning 
\begin{equation}\label{oiuytr}
\deg(e_{ij})=\delta_i - \delta_j
\end{equation} and extend it linearly.  We denote the graded matrix ring obtained this way by $$\M_n(A)(\delta_1,\dots,\delta_n).$$ Note that by the definition above this is a good grading on $\M_n(A)$ and furthermore if $A$ is a field, any good grading on $\M_n(A)$ is obtained in this way (see~\cite{dascalescu}, Prop.~1.2). Also note that if 
$\pi \in S_n$ is a permutation and $\sigma \in \Omega$, then the map sending $e_{ij}$ to 
$e_{\pi(i)\pi(j)}$ induces a graded isomorphism 
\begin{equation}\label{kjhgf}
\M_n(A)(\delta_1,\dots,\delta_n) \cong_{\gr}  \M_n(A)(\delta_{\pi(1)}+\sigma,\dots,\delta_{\pi(n)}+\sigma).
\end{equation}
 In fact if $A$ is a division ring, then any two isomorphic graded matrix algebras are of this form (see~\cite[Theorem 2.1]{can2}). This will be used in Theorem~\ref{acyclicc} to classify Leavitt path algebras of acyclic graphs. Note that the aforementioned results from the  papers \cite{can2,dascalescu} were established when $A$ is a field. But one can easily observe that these theorems are also valid for division rings as well. In fact, in order to classify Leavitt path algebras of $C_n$-comets, we  need to establish  a similar statement  for a graded division algebra $A$ (Theorem~\ref{clasigr}) which covers Theorem~2.1 in ~\cite{can2} when $A$ has a trivial grading. 

Next, we define a grading on a matrix ring which comes from the grading of the base ring. Let $B=\textstyle{\bigoplus_{w\in\Omega}}B_w$ be a $\Omega$-graded ring and $\Ga$ a subgroup of $\Omega$ and let 
$A=\textstyle{\bigoplus_{ \omega \in \Ga}} B_{\omega}$ be a graded subring of $B$, i.e., 
$A$ is a $\Omega$-graded ring with $A_\gamma=0$ if $\gamma\not \in \Ga$ and $A_\gamma=B_\gamma$ otherwise. In this paper, on several occasions, $B=K[x,x^{-1}]=\textstyle{\bigoplus_{i\in \mathbb Z}} Kx^i$ and for a fixed $s$, $A=K[x^s,x^{-s}]= \textstyle{\bigoplus_{i\in s\mathbb Z}} Kx^i$, where $K$ is a field or a division ring.
Let $(\de_1 , \ldots , \de_n) \in \Omega^n$ and let $x$ be a homogeneous element of $A$. Define a grading on $\M_n(A)$ by assigning 
 \begin{equation}\label{hogr}
\deg(e_{ij}(x))=\deg(x)+\delta_i-\delta_j,
\end{equation} and extend it linearly. 
 One can see that  for $\la \in \Omega$, ${\M_n (A)}_{\la}$ is the $n \times n$-matrices over
$A$ with the degree shifted as follows:
\begin{equation}\label{mmkkhh}
{\M_n(A)}_{\la} =
\begin{pmatrix}
A_{ \la+\de_1 - \de_1} & A_{\la+\de_2  - \de_1} & \cdots &
A_{\la +\de_n - \de_1} \\
A_{\la + \de_1 - \de_2} & A_{\la + \de_2 - \de_2} & \cdots &
A_{\la+\de_n  - \de_2} \\
\vdots  & \vdots  & \ddots & \vdots  \\
A_{\la + \de_1 - \de_n} & A_{ \la + \de_2 - \de_n} & \cdots &
A_{\la + \de_n - \de_n}
\end{pmatrix}.
\end{equation}
Thus ${\M_n(A)}_{\la}$ consists of matrices with the $ij$-entry in
$A_{\la+\delta_j  - \delta_i}$. This defines a grading on $\M_n(A)$ as follows:
\[
\M_n (A)=\bigoplus_{\la \in \Omega} {\M_n (A)}_{\la}.
\] We denote this matrix ring with this grading (also) by $\M_n (A)(\de_1 , \ldots , \de_n)$. 
Note that if $A$ has a trivial grading,  this construction reduces to the group grading described in~(\ref{oiuytr}).

The following two statements can be proved easily (see~\cite[pp.~60-61]{grrings},  or Theorem~\ref{clasigr}):
\begin{itemize}
\item[$\circ$] If $\alpha \in \Omega$, and $\pi \in S_n$ is a permutation then 
\begin{equation} \label{pqow1}
\M_n (A)(\de_1 , \ldots , \de_n)\cong_{\gr} \M_n (A)(\de_{\pi(1)}+\alpha , \ldots , \de_{\pi(n)}+\alpha).
\end{equation}

\item[$\circ$]   If $\alpha_1,\dots, \alpha_n \in \Ga$ with $\alpha_i=\deg(u_i)$ for some units $u_i \in A^h$, then 
\begin{equation}\label{pqow2}
\M_n (A)(\de_1 , \ldots , \de_n)\cong_{\gr} \M_n (A)(\de_1+\alpha_1 , \ldots , \de_n+\alpha_n).
\end{equation}
\end{itemize}

In Theorem~\ref{grcomet}, we will show that for a $C_n$-comet graph $E$  we have a $\mathbb Z$-graded isomorphism  
\[\LL_K(E) \cong_{\gr} \M_{m}\big(K[x^n,x^{-n}] \big)\big (|p_1|,\dots, |p_m|),\] where $n$ is the length of the cycle $C$, the set $\{p_i \mid 1\leq i \leq m\}$ are all paths which end in an arbitrary but fixed vertex $u$ on $C$ and do not contain $C$, and $|p_i|$ are length of these paths.  Thus for the graph 

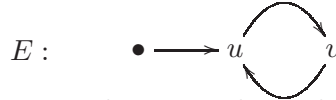
\begin{equation*}
\xymatrix{
E:  & \bullet \ar[r] & u \ar@/^1.5pc/[r] & v \ar@/^1.5pc/[l] &  
}
\end{equation*}

\medskip
\noindent  depending on which fixed vertex, $u$ or $v$, we choose on the cycle, we obtain \[\LL_K(E)\cong_{\gr}\M_3(K[x^2,x^{-2}])(0,1,1)\cong_{\gr}\M_3(K[x^2,x^{-2}])(0,1,2),\] which is justified by the above discussion on shifting of grading. However note that $\LL_K(E)$ is not $\mathbb Z$-graded isomorphic to $A=\M_3(K[x^2,x^{-2}])(0,0,0)$, as one can easily see that $\dim_K\LL_K(E)_0=5$ whereas $\dim_KA_0=9$. 

%Note that if $\Omega=\Gamma$ we get the usual grading on matrix rings which is available in the literature. However the simple modification above seems to introduce many more cases and is needed in the paper. 

Let $A= \textstyle{\bigoplus_{\ga \in \Ga'}} A_{\ga}$ and $B= \textstyle{\bigoplus_{\ga \in
\Ga''}} B_{\ga}$ be graded rings, such that there is a group
$\Gamma$ containing $\Ga'$ and $\Ga''$ as subgroups.  Then $A\times B$ has a natural grading  given by $A\times B  = \textstyle{\bigoplus_{\ga \in \Gamma}} (A \times B)_{\ga}$ where $(A \times B)_{\ga}=A_\ga \times B_\ga$.  
Similarly,  if $A$ and $B$ are $K$-modules for a field $K$ (where here $K$ has a trivial grading), then $A
\otimes_K B$ has a natural grading  given by
$A \otimes_K B = \textstyle{\bigoplus_{\ga \in \Gamma}} (A \otimes_K B)_{\ga}$
where,
\begin{equation}\label{tengr}
(A \otimes_K B)_{\ga} = \Big \{ \sum_i a_i \otimes b_i   \mid   a_i \in
A^h, b_i \in B^h, \deg(a_i)+\deg(b_i) = \ga \Big\}.
\end{equation}
Let $A$ be a $K$-algebra where $K$ is a field,  and moreover, let $A$ be a $\Gamma$-graded ring.  Let $\Omega$ be an abelian group such that $\Gamma \subseteq \Omega$. Then by the definition of grading on tensor products and  by~(\ref{oiuytr}) and~(\ref{hogr}) there is a natural graded isomorphism 
\begin{equation}\label{mozart}
\M_n(K)(\delta_1,\dots,\delta_n) \otimes_K A \cong_{\gr} \M_n(A)(\delta_1,\dots,\delta_n),
\end{equation}
where $\{\delta_1,\dots,\delta_n\} \subseteq \Omega$. This grading on tensor products will be used in Theorem~\ref{polyhead}. 

As we mentioned earlier, some of the rings we are  dealing with in this note are of the form $K[x,x^{-1}]$ where $K$ is a division ring. This is an example of a graded division ring. A nonzero $\Ga$-graded ring $A$  is called a {\it graded division ring} if  every nonzero homogeneous element has an inverse.  It follows that $A_0$ is a division ring. A commutative graded division ring is called a {\it graded field}. Similar to the non-graded setting, one can show that any graded module $M$ over a graded division ring $A$ is graded free, i.e., it is generated by a homogeneous basis and the graded bases have the same cardinality (see~\cite[Proposition~4.6.1]{grrings}). Moreover, if $N$ is a graded submodule of $M$, then 
\begin{equation}\label{dimcouti}
\dim_A(N)+\dim_A(M/N)=\dim_A(M).
\end{equation} In this note, all graded division rings and fields have  torsion free abelian group gradings (in fact, in all our statements $\Ga=n \mathbb Z$ for some $n\in \mathbb N$). However, this assumption is not necessary for the statements below.  

Let $A$ be a $\Ga$-graded division ring, $\Omega$ be a group  such that $\Ga \subseteq \Omega$ and $\M_n(A)(\lambda_1,\dots,\lambda_n)$ be a  {\it graded central simple ring}, where $\lambda_i \in \Omega$, $1\leq i\leq n$. 
Consider the quotient group $\Omega/\Ga$ and let $\Ga+\varepsilon_1,\dots,\Ga+\varepsilon_k$ be the distinct elements in $\Omega/\Ga$ representing the cosets $\Ga+\lambda_i$, $1\leq i\leq n$, and for each $\varepsilon_l$, let $r_l$ be the number of $i$ with $\Ga+\lambda_i=\Ga+\varepsilon_l$.  It was observed in \cite[Proposition~1.4]{hwad} that 
\begin{equation}\label{urnha}
\M_n(A)_0 \cong \M_{r_1}(A_0)\times \dots \times \M_{r_k}(A_0)
\end{equation}
and in particular $\M_n(A)_0$ is a simple ring if and only if $k=1$. Indeed, using (\ref{pqow1}) and (\ref{pqow2}) we get 
\begin{equation}\label{wdeild}
\M_n(A)(\lambda_1,\dots,\lambda_n)\cong_{\gr}\M_n(A)(\varepsilon_1,\dots,\varepsilon_1,\varepsilon_2,\dots,\varepsilon_2,\dots,\varepsilon_k,\dots,\varepsilon_k),
\end{equation}
with each $\varepsilon_l$ occurring $r_l$ times.  Now (\ref{mmkkhh}) for $\lambda=0$ and $(\delta_1,\dots,\delta_n)=(\varepsilon_1,\dots\varepsilon_1,\varepsilon_2,\dots,\varepsilon_2,\dots,\varepsilon_k,\dots,\varepsilon_k)$ immediately gives (\ref{urnha}).

The following statement  is the graded version of
a similar statement on simple rings (see \cite[\S IX.1]{hungerford}). This
is required for the proof of Theorems~\ref{clasigr} and~\ref{grcomet}.

\begin{proposition}\label{prohungi}
Let $A$ and $B$ be $\Ga$-graded division rings and $\Omega$ a group containing $\Ga$. If \[\M_n(A)(\lambda_1,\dots,\lambda_n) \cong_{\gr} \M_{m}(B)(\gamma_1,\dots,\gamma_m)\] as graded
rings, where $\lambda_i,\gamma_j \in \Omega$, $1\leq i \leq n$, $1\leq j \leq m$, then $n = m$ and $A \cong_{\gr} B$.
\end{proposition}
\begin{proof}
The proof follows the ungraded case (see \cite[\S IX.1]{hungerford}) with an extra attention given to the grading. We refer the reader to ~\cite[\S4.3]{millar} for the proof. 
\end{proof}

We can further determine the relations between the graded shifting $(\lambda_1,\dots,\lambda_n)$ and $(\gamma_1,\dots,\gamma_m)$ in the above proposition. For this we need to extend ~\cite[Theorem 2.1]{can2} (see also \cite[Theorem 9.2.18]{grrings}) from fields (with trivial grading) to graded division algebras. The following theorem states that two graded matrix algebras over a graded division ring with two shiftings are isomorphic if and only if one can obtain one shifting from the other by applying (\ref{pqow1}) and (\ref{pqow2}). This will be used in Theorem~\ref{grcomet} to classify Leavitt path algebras of multi-headed comets (in particular $C_n$-comets). 

\begin{theorem}\label{clasigr}
Let $A$ be a $\Ga$-graded division ring and $\Omega$ a group containing $\Ga$.
%were both $\Omega$ and $\Gamma$ are torsion free abelian groups. 
Then for $\lambda_i,\gamma_j \in \Omega$, $1\leq i \leq n$, $1\leq j \leq m$,  
\begin{equation}\label{uytrtyu}
\M_n(A)(\lambda_1,\dots,\lambda_n)\cong_{\gr}\M_m(A)(\gamma_1,\dots,\gamma_m)
\end{equation} if and only if $n=m$ and for a suitable permutation $\pi \in S_n$, we have $\lambda_i=\gamma_{\pi(i)}+\tau_i+\sigma$, $1\leq i \leq n$, where  $\tau_i\in \Ga$ and a fixed $\sigma \in \Omega$, i.e., $(\lambda_1,\dots,\lambda_n)$ is obtained from $(\gamma_1,\dots,\gamma_m)$ by applying {\upshape (\ref{pqow1})} and {\upshape(\ref{pqow2})}.
\end{theorem}
\begin{proof}
In the course of this proof, we work with right $A$-modules. For one direction we need to prove statements (\ref{pqow1}) and (\ref{pqow2}). These are known and hold for any graded ring $A$. We give a short proof. So let $n=m$. Let $V$ be a graded vector space over $A$ with a homogeneous basis $v_1,\dots,v_n$ of degree $\lambda_1,\dots,\lambda_n$, respectively. Defining $E_{ij}(v_t)=\delta_{j,t}v_i$, $1\leq i,j,t \leq n$, we have $E_{ij}\in \End_A(V)_{\lambda_i-\lambda_j}$. It is easy to see that $\End_A(V)\cong_{\gr} \M_n(A)(\lambda_1,\dots,\lambda_n)$ and that $E_{ij}$ corresponds to the matrix $e_{ij}$ (see~\ref{hogr}). Now let  $\pi \in S_n$. Rearranging the homogeneous basis as $v_{\pi(1)},\dots,v_{\pi(n)}$ and 
defining the $A$-graded isomorphism $\phi:V \rightarrow V$ by $\phi(v_i)=v_{\pi(i)}$, we get a graded isomorphism in the level of endomorphism rings \[ \M_n(A)(\lambda_1,\dots,\lambda_n)\cong_{\gr} \End_A(V) \stackrel{\phi}{\longrightarrow} \End_A(V)\cong_{\gr}  \M_n(A)(\lambda_{\pi(1)},\dots,\lambda_{\pi(n)}).\]
Moreover, (\ref{mmkkhh}) shows that it does not make any difference adding a fixed $\alpha \in \Omega$ to each of the entries in the shifting. This gives us (\ref{pqow1}). 

For (\ref{pqow2}), let $\alpha_i\in \Ga$, $1\leq i \leq n$ with $\alpha_i=\deg(u_i)$ for some units $u_i \in A^h$ (here if $A$ is a graded division ring, then all homogeneous elements are invertible and thus any set of $\alpha_i\in \Ga$, $1\leq i \leq n$ can be chosen). Consider the basis $v_iu_i$, $1\leq i \leq n$ for $V$. With this basis, 
$\End_A(V)\cong_{\gr} \M_n (A)(\de_1+\alpha_1 , \ldots , \de_n+\alpha_n).$ Consider the $A$-graded isomorphism $\id:V\rightarrow V$, by $\id(v_i)= (v_i u_i)u_i^{-1}$. A similar argument as above now gives (\ref{pqow2}). 

We now prove the converse of this theorem. That $n=m$ follows from Proposition~\ref{prohungi}.  Applying the same procedure mentioned above, one can find $\epsilon=(\varepsilon_1,\dots,\varepsilon_1,\varepsilon_2,\dots,\varepsilon_2,\dots,\varepsilon_k,\dots,\varepsilon_k)$ in $\Omega$  such that  $\M_n(A)(\lambda_1,\dots, \lambda_n)\cong_{\gr}\M_n(A)(\epsilon)$ as in (\ref{wdeild}). 

Now set $V=A(-\varepsilon_1)\times \dots \times A(-\varepsilon_1) \times  \dots \times A(-\varepsilon_k)\times \dots \times A(-\varepsilon_k)$ and  pick the (standard) homogeneous basis $e_i$, $1\leq i \leq n$  and define $E_{ij} \in \End_A(V)$ by $E_{ij}(e_t)=\delta_{j,t}e_i$, $1\leq i,j,t \leq n$.
One can easily see that $E_{ij}$ is a $A$-graded homomorphism of degree $\varepsilon_{s_i}-\varepsilon_{s_j}$ where $\varepsilon_{s_i}$ and $\varepsilon_{s_j}$ are $i$-th and $j$-th elements in $\epsilon$. Moreover, 
$\End_A(V) \cong_{\gr} \M_n(A)(\epsilon)$ and $E_{ij}$ corresponds to the matrix $e_{ij}$ in $\M_n(A)(\epsilon)$. In a similar manner, one can find $\epsilon'=(\varepsilon_1',\dots,\varepsilon_1',\varepsilon_2',\dots,\varepsilon_2',\dots,\varepsilon_{k'}',\dots,\varepsilon_{k'}')$ and a graded $A$-vector space $W$ such that $\M_n(A)(\gamma_1,\dots,\gamma_n)\cong_{\gr}\M_n(A)(\epsilon')$, and $\End_A(W)\cong_{\gr} \M_n(A)(\epsilon')$. Therefore (\ref{uytrtyu}) provides  a graded ring isomorphism $\theta:\End_A(V)\rightarrow \End_A(W)$. Define $E_{ij}':=\theta(E_{ij})$ and $E_{ii}'(W)=Q_i$, for $1\leq i,j \leq n$. Since $\{E_{ii} \mid 1\leq i \leq n\}$ is a complete system of orthogonal idempotents, so is $\{E_{ii}' \mid 1\leq i \leq n\}$. It follows that $W\cong_{\gr} \textstyle{\bigoplus_{1\leq j \leq n}} Q_j$. Furthermore, $E_{ij}'E_{tr}'=\delta_{j,t}E_{ir}'$ and $E_{ii}'$ acts as identity on $Q_i$. These relations show that restricting $E_{ij}'$ on $Q_j$ induces an $A$-graded isomorphism $E_{ij}':Q_j \rightarrow Q_i$ of degree 
$\varepsilon_{s_i}-\varepsilon_{s_j}$ (same degree as $E_{ij}$). 
So   $Q_j\cong_{\gr} Q_1(\varepsilon_{s_1}-\varepsilon_{s_j})$ for any $1\leq j \leq n$.  
Therefore 
$W\cong_{\gr} \textstyle{\bigoplus_{1\leq j \leq n}} Q_1(\varepsilon_{s_1}-\varepsilon_{s_j})$. By dimension count (which is valid here, see~\ref{dimcouti}), it follows that $\dim_AQ_1=1$. 

A similar argument for the identity map $\id:\End_A(V)\rightarrow \End_A(V)$ produces $V\cong_{\gr} \textstyle{\bigoplus_{1\leq j \leq n}} P_1(\varepsilon_{s_1}-\varepsilon_{s_j})$, where $P_1=E_{11}(V)$, with $\dim_AP_1=1$. 

Since $P_1$ and $Q_1$ are $A$-graded vector spaces of dimension $1$, there is $\sigma \in \Omega$, such that $Q_1\cong_{\gr}P_1(\sigma)$. Using  the fact that for an $A$-graded module $P$ and $\alpha,\beta \in \Omega$,  $P(\alpha)(\beta)=P(\alpha+\beta)=P(\beta)(\alpha)$, we have 
\[W\cong_{\gr} \textstyle{\bigoplus_{1\leq j \leq n}} Q_1(\varepsilon_{s_1}-\varepsilon_{s_j})\cong_{\gr}
\textstyle{\bigoplus_{1\leq j \leq n}} P_1(\sigma)(\varepsilon_{s_1}-\varepsilon_{s_j})\cong_{\gr} \textstyle{\bigoplus_{1\leq j \leq n}} P_1(\varepsilon_{s_1}-\varepsilon_{s_j})(\sigma)\cong_{\gr}V(\sigma).\] We denote this $A$-graded isomorphism with $\phi:W \rightarrow V(\sigma)$.
Let $e_i'$, $1\leq i \leq n$  be a (standard) homogeneous basis of degree $\varepsilon_{s_i}'$ in $W$. Then $\phi(e_i')=\sum_{1\leq j \leq n} e_j a_j$, where $a_j \in A^h$ and  $e_j$ are homogeneous of degree $\varepsilon_{s_j}-\sigma$ in $V(\sigma)$.  Since $\deg(\phi(e_i'))=\varepsilon_{s_i}'$, all $e_j$'s with non-zero $a_j$ in the sum have the same degree. For if $\varepsilon_{s_j}-\sigma=\deg(e_j)\not = \deg(e_l)=\varepsilon_{s_l}-\sigma$, then since $\deg(e_j a_j)=\deg(e_l a_l)=\varepsilon_{s_i}'$ it follows that $\varepsilon_{s_j}-\varepsilon_{s_l} \in \Ga$ which is a contradiction as $\Ga+\varepsilon_{s_j}$ and $\Ga+\varepsilon_{s_l}$ are distinct. Thus $\varepsilon_{s_i}'=\varepsilon_{s_j}+\tau_j-\sigma$ where $\tau_j=\deg(a_j) \in \Ga$. In the same manner one can show that, $\varepsilon_{s_i}'=\varepsilon_{s_{i'}}'$  in $\epsilon'$ if and only if  $\varepsilon_{s_j}$ and $\varepsilon_{s_{j'}}$ assigned to them by the previous argument are also equal. 
This shows that $\epsilon'$ can be obtained from $\epsilon$ by applying (\ref{pqow1}) and (\ref{pqow2}). Since $\epsilon'$ and $\epsilon$ are also obtained from $\gamma_1,\dots,\gamma_n$ and $\lambda_1,\dots,\lambda_n$, respectively, by applying  (\ref{pqow1}) and (\ref{pqow2}), putting these together shows that $\lambda_1,\dots,\lambda_n$ and $\gamma_1,\dots,\gamma_n$ have the similar relations, i.e.,   $\lambda_i=\gamma_{\pi(i)}+\tau_i+\sigma$, $1\leq i \leq n$, where  $\tau_i\in \Ga$ and a fixed $\sigma \in \Omega$.
\end{proof}

\subsection{Graphs and Leavitt path algebras}

In this subsection we gather some  graph-theoretic definitions and recall the basics on Leavitt path algebras. The reader familiar with this topic can skip to Section~\ref{sec2}. 

A {\it directed graph} $E=(E^0,E^1,r,s)$ consists of two countable sets $E^0$, $E^1$ and maps $r,s:E^1\rightarrow E^0$. The elements of $E^0$ are called {\it vertices} and the elements of $E^1$ {\it edges}. If $s^{-1}(v)$ is a finite set for every $v \in E^0$, then the graph is called {\it row-finite}. In this note we will consider only row-finite graphs. In this setting, if the number of vertices, i.e.,  $|E^0|$,  is finite, then the number of edges, i.e.,  $|E^1|$, is finite as well and we call $E$ a {\it finite} graph. 

 For a graph $E=(E^0,E^1,r,s)$, a vertex $v$ for which $s^{-1}(v)$ is empty is called a {\it sink}, while a vertex $w$ for which $r^{-1}(w)$ is empty is called a {\it source}. An edge with the same source and range is called a {\it loop}. A path $\mu$ in a graph $E$ is a sequence of edges $\mu=\mu_1\dots\mu_k$, such that $r(\mu_i)=s(\mu_{i+1}), 1\leq i \leq k-1$. In this case, $s(\mu):=s(\mu_1)$ is the {\it source} of $\mu$, $r(\mu):=r(\mu_k)$ is the {\it range} of $\mu$, and $k$ is the {\it length} of $\mu$ which is  denoted by $|\mu|$. We consider a vertex $v\in E^0$ as a {\it trivial} path of length zero with $s(v)=r(v)=v$. 
If $\mu$ is a nontrivial path in $E$, and if $v=s(\mu)=r(\mu)$, then $\mu$ is called a {\it closed path based at} $v$. If $\mu=\mu_1\dots\mu_k$ is a closed path based at $v=s(\mu)$ and $s(\mu_i) \not = s(\mu_j)$ for every $i \not = j$, then $\mu$ is called a {\it cycle}.  %Throughout, we denote a cycle of length $n$ by $C_n$. 

For two vertices $v$ and $w$, the existence of a path with the source $v$ and the range $w$ is denoted by $v\geq w$. Here we allow paths of length zero. By $v\geq_n w$, we mean there is a path of length $n$ connecting these vertices. Therefore $v\geq_0 v$ represents the vertex $v$. Also, by $v>w$, we mean a path from $v$ to $w$ where $v\not = w$. In this note, by $v\geq w' \geq w$, it is understood that there is a path connecting $v$ to $w$ and going through $w'$ (i.e., $w'$ is on the path connecting $v$ to $w$). For $n\geq 2$, we define $E^n$ to be the set of paths of length $n$ and $E^*=\bigcup_{n\geq 0} E^n$, the set of all paths.

\begin{deff}\label{llkas}{\sc Leavitt path algebras.} \label{LPA} \\For a graph $E$ and a ring $R$ with identity, we define the {\it Leavitt path algebra of $E$}, denoted by $\LL_R(E)$, to be the algebra generated by the sets $\{v \mid v \in E^0\}$, $\{ \alpha \mid \alpha \in E^1 \}$ and $\{ \alpha^* \mid \alpha \in E^1 \}$ with the coefficients in $R$, subject to the relations 

\begin{enumerate}
\item $v_iv_j=\delta_{ij}v_i \textrm{ for every } v_i,v_j \in E^0$.

\item $s(\alpha)\alpha=\alpha r(\alpha)=\alpha \textrm{ and }
r(\alpha)\alpha^*=\alpha^*s(\alpha)=\alpha^*  \textrm{ for all } \alpha \in E^1$.

\item $\alpha^* \alpha'=\delta_{\alpha \alpha'}r(\alpha)$, for all $\alpha, \alpha' \in E^1$.

\item $\sum_{\{\alpha \in E^1, s( \alpha)=v\}} \alpha \alpha^*=v$ for every $v\in E^0$ for which $s^{-1}(v)$ is non-empty.

\end{enumerate}
\end{deff}
Here the ring $R$ commutes with the generators $\{v,\alpha, \alpha^* \mid v \in E^0,\alpha \in E^1\}$. When the coefficient ring $R$ is clear from the context, we simply write $\LL(E)$ instead of $\LL_R(E)$. When $R$ is not commutative, then we consider $\LL_R(E)$ as a left $R$-module. The elements $\alpha^*$ for $\alpha \in E^1$ are called {\it ghost edges}.

Setting $\deg(v)=0$, for $v\in E^0$, $\deg(\alpha)=1$ and $\deg(\alpha^*)=-1$ for $\alpha \in E^1$, we obtain a natural $\mathbb Z$-grading on the free $R$-ring generated by  $\{v,\alpha, \alpha^* \mid v \in E^0,\alpha \in E^1\}$. Since the relations in the above definition are all homogeneous, the ideal generated by these relations is homogeneous and thus we have a natural $\mathbb Z$-grading on $\LL_R(E)$. 

If $\mu=\mu_1\dots\mu_k$, where $\mu_i \in E^1$, is an element of $\LL(E)$, then we denote by $\mu^*$ the element $\mu_k ^*\dots \mu_1^* \in \LL(E)$. Since $\alpha^* \alpha'=\delta_{\alpha \alpha'}r(\alpha)$, for all $\alpha, \alpha' \in E^1$, any word can be written as $\mu \gamma ^*$ where $\mu$ and $\gamma$ are paths in $E$.  The elements of the form $\mu\gamma^*$ are called {\it monomials}. 

Taking the grading into account, one can write $\LL_R(E) =\textstyle{\bigoplus_{k \in \mathbb Z}} \LL_R(E)_k$ where,
\[\LL_R(E)_k=  \Big \{ \sum_i r_i \alpha_i \beta_i^*\mid \alpha_i,\beta_i \textrm{ are paths}, r_i \in R, \textrm{ and } |\alpha_i|-|\beta_i|=k \textrm{ for all } i \Big\}.\] 
For simplicity we denote $\LL_R(E)_k$, the homogeneous elements of degree $k$, by $\LL_k$.

We define an (anti-graded) involution 
on $\LL_R(E)$ by $\overline {\mu\gamma^*}=\gamma\mu^*$ for the monomials and extend it to the whole $\LL_R(E)$ in the obvious manner. Note that if $x\in \LL_R(E)_n$, then $\overline x \in \LL_R(E)_{-n}$.

  When $R$ is a division ring, by constructing a representation of $\LL_R(E)$ in $\End(V)$, for a suitable vector space $V$, one can show that the vertices of a graph $E$ are linearly independent in $\LL_R(E)$ and the edges and ghost edges are not zero (see Lemma~1.5 in \cite{goodearl}). We will carry this over to the generalized setting of weighted Leavitt path algebras (Theorem~\ref{liniin}) and therefore cover the special case of LPA as a corollary (by setting the weight map the constant map $1$).

Throughout the note we need some more definitions which we gather here.

\begin{deff}\hfill  \label{mulidef}
\begin{enumerate}

\item A path  which does not contain a cycle is called a {\it acyclic} path. 

\item A graph without cycles is called a {\it acyclic} graph.

\item Let $v \in E^0$. Then the {\it out-degree} and the {\it total-degree} of $v$ are defined as $\outdeg(v)=\card(s^{-1}(v))$ and 
$\totdeg(v)=\card(s^{-1}(v) \cup r^{-1}(v))$, respectively.

\item A finite graph $E$ is called a {\it line graph} if it is connected, acyclic and $\totdeg(v)\leq 2$ for every 
$v \in E^0$.  If we want to emphasize the 
number of vertices, we say that $E$ is an $n$-line graph whenever $n=\card(E^0)$. An {\it oriented} $n$-line graph $E$ is an $n$-line graph such that  $E^{n-1} \not =\emptyset$.  

\item For  any vertex $v \in E^0$, the cardinality of the set $R(v)=\{\alpha \in E^* \mid r(\alpha)=v\}$ is denoted by $n(v)$. 

\item For any vertex $v \in E^0$, the {\it tree} of $v$, denoted by $T(v)$,  is the set $\{w\in E^0 \mid v\geq w \}$. Furthermore, for $X\subseteq E^0$, we define $T(X)=\bigcup_{x \in X}T(x)$. 

\end{enumerate}
\end{deff}

\begin{deff}\label{petaldef}
A {\it rose with $k$-petals} is a graph which consists of one vertex and $k$ loops. We denote this graph by $L_k$ and its vertex by $s(L_k)$. The Leavitt path algebra of this graph with coefficient in $R$ is denoted by $\LL_R(1,k)$. (Cohn's notation in~\cite{cohn11} for this algebra is $V_{1,k}$ and this algebra is of type $(1,k-1)$.)
We allow $k$ to be zero and in this case $L_0$ is just a vertex with no loops. With this convention, one can easily establish that $\LL_R(1,0)\cong R$ and $\LL_R(1,1)\cong R[x,x^{-1}]$.  In a graph which contains a rose $L_k$,  we say $L_k$ does not have an exit, if there is no edge $e$ with $s(e)=s(L_k)$ and $r(e) \not = s(L_k)$. 
\end{deff} 

 We need to recall the definition of morphisms between two graphs in order to consider the category of directed graphs.

For two directed graphs $E$ and $F$, a {\it complete  graph homomorphism} $f:E\rightarrow F$ consists of a map $f^0:E^0 \rightarrow F^0$ and $f^1:E^1 \rightarrow F^1$ such that $r(f^1(\alpha))=f^0(r(\alpha))$ and $s(f^1(\alpha))=f^0(s(\alpha))$ 
for any $\alpha \in E^1$, additionally, $f^0$ is injective and $f^1$ restricts to a bijection from $s^{-1}(v)$ to $s^{-1}(f^0(v))$ for every $v\in E^0$ which emits edges. One can check that such a map induces a graded homomorphism on the level of LPAs. i.e, there is a graded homomorphism $\LL(E) \rightarrow \LL(F)$.

\section{Strongly graded Leavitt path algebras}\label{sec2}

In this section we characterize those graphs whose Leavitt path algebras are strongly graded (Theorems~\ref{sth} and ~\ref{sthfin}).  In this section, unless it is noted otherwise, $R$ is a ring with identity.

\begin{deff}\label{infipa}
We say a vertex $v$ is on an {\it infinite path} if there are real paths and ghost paths of arbitrary length starting from $v$, i.e., for any $n\in \mathbb N$ there are vertices $w$ and $w'$ such that $v\geq_n w$ and $v\leq_n w'$, respectively.  
A vertex $v$ is {\it connected to an infinite path} if there is a vertex $w$ such that $v\geq w$ and $w$ is on an infinite path. 
\end{deff}

 The following are easy to observe and will be used in the text. The proofs are left to the reader.  
 
 \begin{lemma}\label{jh543p}\hfill
 
 \begin{enumerate}[\upshape(1)]

\item If a vertex  is not connected to an infinite path, then the graph has a vertex which is either a source or a sink. 

\item If $x\in \LL(E)$ is a monomial of the form, $x=\alpha\beta^*$, then $x\overline x x =x$. 

\item \label{psp} If $\LL_n \LL_m \not = 0$, $n >0$, then there is a vertex $v$ and a path of length at least $n$ emitting from $v$. 

\item \label{ss} If there is a sink, say $v$, then $v \in \LL_0$, however $v \not \in \LL_1\LL_{-1}$. Therefore  if $\LL(E)$ is strongly graded and $E^1$ is not empty, then the support of $\LL(E)$ is $\mathbb Z$. 
\end{enumerate}
\end{lemma}

\begin{deff}\label{exiedg}
For an acyclic path $p=\mu_1\mu_2\dots\mu_k$, an {\it exit edge} $e$ (or an {\it escape} edge), is an edge such that $s(e)=s(\mu_i)$ for some $1\leq i \leq k$ but $e \not = \mu _i$ unless $i=k$. By this definition, $\mu_k$ is an exit edge for $p$. A path $q$ is called an {\it escape path} of $p$, if $q$ is an exit edge for $p$ with $s(p)=s(q)$ or $q=\mu_1 \dots \mu_te$, with $t<k$ and $e$ is an exit edge of $p$. The set of all escape paths of $p$ is denoted by $\Esc(p)$. Note that $p \in \Esc(p)$.
 \end{deff}
 
\begin{example} 
Let $p=\mu_1\mu_2\mu_3\mu_4$ be the path in the graph below
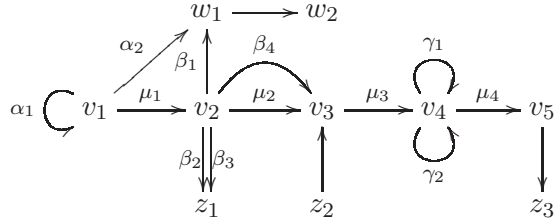
\begin{equation*}
\xymatrix{
 & w_1 \ar[r] & w_2\\
 v_1 \ar@(lu,ld)_{\alpha_1} \ar[r]^{\mu_1} \ar[ur]^{\alpha_2}& v_2 \ar@/^1.5pc/[r]^{\beta_4} \ar[r]^{\mu_2} \ar@<1.5pt>[d]_{\beta_2}  \ar@<-1.5pt>[d]^{\beta_3} \ar[u]^{\beta_1}&  v_3 \ar[r]^{\mu_3} &  v_4 \ar[r]^{\mu_4} \ar@(ul,ur)^{\gamma_1}  \ar@(dl,dr)_{\gamma_2}&  v_5\ar[d]\\
 & z_1  &  z_2 \ar[u] & & z_3}
\end{equation*} 
Then we have $\Esc(p)=\{\alpha_1,\alpha_2,\mu_1\beta_1,\mu_1\beta_2,\mu_1\beta_3,\mu_1\beta_4, \mu_1\mu_2\mu_3\gamma_1,\mu_1\mu_2\mu_3\gamma_2,\mu_1\mu_2\mu_3\mu_4\}$. 
\end{example}

\begin{example}\label{lpa2011}
If $p=\mu_1$ is an edge with $s(\mu_1)=v$ and $r(\mu_1)=u$ then by Definition~\ref{exiedg}, 
$\Esc(p)=\{\alpha \in E^1 \mid s(\alpha)=s(\mu_1)\}$.
 It follows that 
 \[v=\sum_{\{\alpha\in E^1\mid s(\alpha)=v\}} \alpha \alpha^*=\sum_{\alpha \in \Esc(p)}\alpha\alpha^*.\] This will be used as the first step of induction in Lemma~\ref{lelele}. 
 \end{example}
 
 Equation~\ref{edgi} in the lemma below will be used in an essential way in the Lemma~\ref{nicelem}.

\begin{lemma}\hfill \label{lelele}
\begin{enumerate}[\upshape(1)]

\item If $p$ and $q$ are two finite acyclic paths with $r(p)=s(q)$ such that  
$p$ and $q$ are not vertices and $pq$ is also acyclic then 
\begin{equation}\label{verti}
\Esc(pq)=\Esc(p) \backslash \{p\} \cup \{p\alpha \mid \alpha \in \Esc(q)\}. 
\end{equation}

\item For a finite acyclic  path $p$, if $v=s(p)$  then 
\begin{equation}\label{edgi}
v=\sum_{\alpha \in  \Esc(p)}\alpha \alpha^*.
\end{equation}
\end{enumerate}
 \end{lemma}

\begin{proof}
 (1) This follows easily and is left to the reader. 
 The only care needs to be given is, by Definition~\ref{exiedg}, $p$ is an escape path for $p$, which needs to be removed as it is not an escape path for $pq$. 
 
(2) We prove the statement by using an induction on the length of $p$. Let $|p|=1$. Then~(\ref{edgi}) reduces to $v=\sum_{\{\alpha\in E^1\mid s(\alpha)=v\}} \alpha \alpha^*$ which is Relation~(4) in Definition~\ref{llkas} and so is valid (see Example~\ref{lpa2011}). Now let $p=\mu_1\dots\mu_n$ and suppose the statement is valid for paths of length $n-1$. Considering the path $q=\mu_2\dots\mu_n$ with $s(q)=u$ we then have
\begin{equation}
v=s(\mu_1)=\sum_{\{\alpha\in E^1\mid s(\alpha)=v, \alpha \not = \mu_1\}} \alpha \alpha^*+\mu_1\mu_1^*,
\textrm{         and         }\,\,\,\,\,\,\, u=s(q)=\sum_{\beta \in \Esc(q)}\beta\beta^*.
\end{equation}
Writing $\mu_1  \mu_1^*$ as $\mu_1u \mu_1^*$ and replacing $u$  with the above equation we have 
\begin{align}\label{sed}
v & =\sum_{\{\alpha\in E^1\mid s(\alpha)=v, \alpha \not = \mu_1\}} \alpha \alpha^*+
\mu_1\Big(\sum_{\beta \in \Esc(q)}\beta\beta^*\Big )\mu_1^*\notag\\
& =\sum_{\alpha \in \Esc(\mu_1)\backslash \{\mu_1\}} \alpha \alpha^*+
\Big(\sum_{\beta \in \Esc(q)}\mu_1\beta (\mu_1\beta)^*\Big )
\end{align}
Now by ~(\ref{verti}), \[\Esc(p)=\Esc(\mu_1 q)=\Esc(\mu_1) \backslash \{\mu_1\}  \cup 
\{\mu_1\beta \mid \beta \in \Esc(q)\} 
,\] which 
guarantees that Equation~\ref{sed} coincides with~(\ref{edgi}), so we are done. 
\end{proof}

\begin{deff}\label{immedi}
Let $\mathcal P$ be a subset of vertices (with a certain property).  For a vertex $v\in E^0$, we say that $w\in \mathcal P$ is {\it immediate} to $v$, if there is a path $p$ with the source $v$ and range $w$ such that no vertices on $p$ is in $\mathcal P$ except $w$. More formally, $w$ is immediate to $v$ if there is a path $p=\alpha_1\dots\alpha_k$ such that $s(p)=v$, $r(p)=w$ and $s(\alpha_i) \not \in \mathcal P, 1\leq i \leq k$.  Such a path is called an {\it immediate path} from $v$ to $\mathcal P$.  A non-empty subset $\mathcal P \subseteq E^0$ is called {\it dense}, if for any vertex $v\in E^0$, there is a path connecting $v$ to $\mathcal P$, i.e, $T(v) \cap \mathcal P \not = \emptyset$. 
%i.e., there is a path $p$ with $s(p)=v$ and $r(p) \in \mathcal P$. 

For $v\in E^0$, the {\it orbit} of $v$ with respect to $\mathcal P$, denoted by $O_{\mathcal P}(v)$, is the set of all immediate paths from  $v$ to $\mathcal P$. 
We define {\it the bound} of $O_{\mathcal P}(v)$ to be the minimum $n\in \mathbb N$ such that all paths in the set $O_{\mathcal P}(v)$ have length at most $n$. If $O_{\mathcal P}(v)$ has a bound, then we say $O_{\mathcal P}(v)$ is {\it bounded}.
 \end{deff}

From the definition, it is clear that  if $v \in \mathcal P$, then $v$ itself is the only vertex which is immediate to $v$. 
Therefore the bound of $O_{\mathcal P}(v) $ is zero if and only if $v\in \mathcal P$. Also, the set $\mathcal P$ is dense if and only if for any $v\in E^0$, $O_{\mathcal P}(v)$ is nonempty.  

\begin{example}
Consider the following graph with $\mathcal P=\{w,w',w''\}$, 
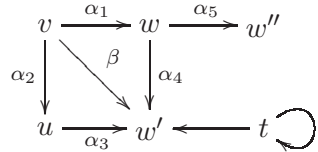
\[
\xymatrix{
v \ar[r]^{\alpha_1} \ar[d]_{\alpha_2} \ar[dr]^\beta &   w \ar[r]^{\alpha_5} \ar[d]^{\alpha_4} &   w''  \\
u \ar[r]_{\alpha_3} & w' & t\ar[l] \ar@(ru,rd)}
\] Then $w$ and $w'$ are immediate to $v$ with the immediate paths $\alpha_1$ and $\alpha_2\alpha_3$ and $\beta$.  However $w''$ is not immediate to $v$. Hence $O_{\mathcal P}(v)=\{\alpha_1,\alpha_2\alpha_3,\beta\}$ and $O_{\mathcal P}(u)=\{\alpha_3\}$. We can write 
\begin{equation*}
 v=\sum_{\alpha \in O_{\mathcal P}(v)}\alpha\alpha^*=
 \alpha_1\alpha_1^*+\beta \beta^* +\alpha_2\alpha_3\alpha_3^*\alpha_2^*.
 \end{equation*}
We prove this equation for dense subsets in Lemma~\ref{nicelem}. Note that $O_{\mathcal P}(t)$ is not bounded. 
\end{example}

 Recall that  $v\geq v' \geq w$, denotes  a path connecting $v$ to $w$ and going through $v'$ (i.e., $v'$ is on the path connecting $v$ to $w$). 

\begin{lemma}\label{counting}
If the set $O_{\mathcal P}(v) $ is bounded then it is finite. Furthermore,  if  $w\in \mathcal P$ is an immediate vertex to $v$ with an immediate path $v> w$ and $v'$ is a vertex on this path, i.e.,  $v>v'\geq w$,  then the bound of $O_{\mathcal P}(v')$ is strictly less than the bound of $O_{\mathcal P}(v)$ and $|O_{\mathcal P}(v')|\leq |O_{\mathcal P}(v)|$.
\end{lemma}
\begin{proof}
Since the graph is row-finite, there are only a finite number of paths of given length emitting from each vertex. This gives the first statement. For the second assertion, let $v>v'\geq w$ where $v>w$ is an immediate path. If $v'=w$, then the bound of $O_{\mathcal P}(v')$ is zero and thus is strictly less than the bound of $O_{\mathcal P}(v)$ (which contains $v>w$). If $v'\not = w$, then none of vertices on  the part $v>v'$ is immediate (otherwise $w$ is not immediate). 
Let $w'$ be an immediate vertex to $v'$, with the immediate path $v'\geq w'$ connecting $v'$ to $w'$. Then $v>v'\geq w'$   makes $w'$ an immediate vertex for $v$ with the length strictly greater the length of $v' \geq w'$. This shows that the bound of $O_{\mathcal P}(v')$ is strictly less than the bound $O_{\mathcal P}(v)$. 
This argument also shows that $|O_{\mathcal P}(v')|\leq |O_{\mathcal P}(v)|$.
\end{proof}

\begin{lemma}\label{nicelem}
Let $\mathcal P$ be a dense subset of $E^0$, $v\in E^0$ and $O_{\mathcal P}(v)$ be the orbit of $v$. If all paths in  $O_{\mathcal P}(v)$  are acyclic, then 
 \begin{equation}\label{gfg}
 v=\sum_{\alpha \in O_{\mathcal P}(v)}\alpha\alpha^*.
 \end{equation}
\end{lemma}
\begin{proof}
 We proceed by induction on the bound of $O_{\mathcal P}(v)$.  If the bound is zero, then $v\in \mathcal P$ and $O_{\mathcal P}(v)$ consists only of the path of length zero $v$ and so  Equation~\ref{gfg} trivially holds. 
 
Let the bound of $O_{\mathcal P}(v)$ be $1$.  (This prevents $v$ having a loop.) We show that 
\[O_{\mathcal P}(v)=\big \{\alpha\in E^1\mid s(\alpha)=v\big \}\] and therefore ~(\ref{gfg}) reduces to $v=\sum_{\{\alpha\in E^1\mid s(\alpha)=v\}} \alpha \alpha^*$ which is Relation~(4) in Definition~\ref{llkas} and so is valid. Let $p \in O_{\mathcal P}(v)$. Since the bound of $O_{\mathcal P}(v)$ is  $1$,   $|p|$ has to be $1$ (if $|p|=0$, then $v \in \mathcal P$ so the bound is $0$). So $p \in \{\alpha\in E^1\mid s(\alpha)=v\}$. Now let $\alpha \in E^1$, with $s(\alpha)=v$. If $r(\alpha) \not \in \mathcal P$, then since $\mathcal P$ is dense, any vertex is connected to $\mathcal P$, including $r(\alpha)$ and so there is an immediate vertex to $v$ with an immediate path of length greater than $1$ which is a contradiction. So $r(\alpha)$ is in $\mathcal P$, so $\alpha \in O_{\mathcal P}(v)$.

Suppose~(\ref{gfg}) is valid for any vertex  with a bound less than the bound of $O_{\mathcal P}(v)$. 
Fix an immediate vertex $w$ to $v$ and consider the immediate path $q$ with $s(q)=v$ and $r(q)=w$, which is acyclic by the assumption of the lemma. Let $\Esc(q)=\{p_1,\dots,p_h,q\} $. Then by Lemma~\ref{lelele}(2), 
\begin{equation}\label{hyp1}
v=p_1p_1^*+\dots+p_hp_h^*+qq^*=p_1v_1p_1^*+\dots+p_hv_hp_h^*+qq^*,
\end{equation}
where $r(p_i)=v_i$. By Lemma~\ref{counting}, the bound of $O_{\mathcal P}(v_i)$ is smaller than the bound of $O_{\mathcal P}(v)$ and clearly all its paths are acyclic, so by induction we have 
\[v_i=\sum_{\alpha \in O_{\mathcal P}(v_i)}\alpha\alpha^*.\] 
Plugging these sums for $v_i$ in Equation~\ref{hyp1} and  observing that for any $v_i$, if $w'$ is an immediate vertex to $v_i$, then it is immediate to $v$ and also, for any immediate vertex $w''\not = w$ of $v$, there is a path connecting $v$ to $w''$, thus $w''$ is an immediate vertex for some $v_i$, we obtain~(\ref{gfg}) for $O_{\mathcal P}(v)$. 
\end{proof}

Specializing $\mathcal P$ to the set of all vertices on an infinite path (see Definition~\ref{infipa}), we can give a criterion for a Leavitt path algebra to be strongly graded. 

\begin{theorem}\label{sth}
Let $E$ be a row-finite graph and $\mathcal P$  be the set of vertices on an infinite path.  
The Leavitt path algebra with coefficients in a ring $R$, associated to a graph $E$ which the orbit of any vertex  is nonempty and bounded, is strongly graded. 
\end{theorem}
\begin{proof}
Suppose any vertex of the graph $E$ is connected to a vertex on an infinite path, i.e., the orbits are not empty (in other words, the set ``vertices on infinite paths'' is dense in $E^0$).
Let $\LL=\LL_R(E)$.  First note that $\LL_n\not = 0$ for any $n \in \mathbb Z$. We need to show that $ \LL_{n+m}=\LL_n \LL_m$ for all $n,m \in \mathbb Z$. 

We will use the followings two facts: 

\begin{enumerate}
\item \label{ext} For any vertex $v$, and $n \in \mathbb N$, one can write $v=\sum p_i q_i^*$ with $|p_i|=|q_i|=n$.  
This can be proved by an easy induction: Clearly $v=v.v$. Since any vertex is connected to an infinite path, the graph does not have a sink. It is also a row-finite graph. So one can write 
\begin{equation}\label{vi}
v= \sum_{\substack{e\in E^1\\ s(e)=v}} e e^*.
\end{equation} So the statement is valid for $n=1$.  
 Now if $v=\sum p_i q_i^*$ with $|p_i|=|q_i|=n-1$, write $v=\sum p_ir(p_i)q_i^*$ and use~(\ref{vi}) for each $r(p_i)$ to conclude by induction. 
 
 \item \label{inf} For a vertex $v$ on an infinite path, and any $n \in \mathbb N$, it is easy to see that one can write $v=p^*p$ where $p$ is a path of length $n$ with $r(p)=v$.  
\end{enumerate}

In order to prove the theorem, since $\LL_n$, where $n\in \mathbb Z$, is an (left) $R$-module, it is enough to show that any monomial $\alpha\beta^* \in \LL_{n+m}$ is in $\LL_n \LL_m$. Writing $\alpha\beta^*=\alpha_1\dots\alpha_k \beta_1^* \dots \beta_l^*$, we have $k-l=n+m$.  
We need to consider two cases.

\noindent {\bf Case $\mathbf{n\geq0}$.} If $k\geq n$ then writing $\alpha\beta^*=(\alpha_1\dots\alpha_n)(\alpha_{n+1}\dots\alpha_k \beta_1^* \dots \beta_l^*)$, it is clear that $\alpha\beta^* \in \LL_n \LL_m$. 

If $k < n$, then by (\ref{ext}), we can write $r(\alpha_k)=\sum p_i q_i^*$ with $|p_i|=|q_i|=n-k$. So 
\[\alpha\beta^*=\alpha_1\dots\alpha_k r(\alpha_k)\beta_1^* \dots \beta_l^*=\sum (\alpha_1\dots \alpha_k p_i)(q_i^* \beta_1^* \dots \beta_l^*) \in \LL_n \LL_m.\] 

\noindent {\bf Case $\mathbf{n<0}$.} If $m<0$, consider $\overline {\alpha\beta^*} \in \LL_{-m-n}$. By the previous case now, 
 $\overline {\alpha\beta^*} \in \LL_{-m}\LL_{-n}$. Applying the involution to this element again we have 
 $\alpha\beta^* \in \LL_n\LL_m$. 
 
 Now for the remaining case (in fact the following argument is valid when $m<0$ as well), let $v=r(\alpha_k)$ and consider $O_{\mathcal P}(v)$ which consists of a finite number of paths by Lemma~\ref{counting}. If $\gamma \in O_{\mathcal P}(v)$, $s(\gamma)=v$ and $r(\gamma)=w$ and so $w$ is on an infinite path and by the definition of immediate path, none of the other vertices on $\gamma$ has this property. This forces $\gamma$ to be an acyclic path. Furthermore, since $\mathcal P$ is dense, by Lemma~\ref{nicelem}, we have  
 \begin{equation}\label{gfgii}
 v=\sum_{\gamma \in O_{\mathcal P}(v)}\gamma\gamma^*=
 \sum_{\gamma \in O_{\mathcal P}(v)}\gamma r(\gamma)\gamma^*.
 \end{equation}

Now in Equation~\ref{gfgii}, since each $r(\gamma)$ is on an infinite path, by ~(\ref{inf}), we can write $r(\gamma)=p_\gamma^*p_\gamma$, where $|p_\gamma|=k+|\gamma|+|n|$. Thus
\[\alpha\beta^*=\alpha_1\dots\alpha_k v \beta_1^* \dots \beta_l^*=
\sum_{\gamma \in O_{\mathcal P}(v)}\Big(\alpha_1\dots\alpha_k\gamma p_\gamma^*\Big)\Big(p_\gamma \gamma^* \beta_1^* \dots \beta_l^*\Big).\] A quick inspection now shows that each term in the sum is in $\LL_n\LL_m$. 
\qedhere
\end{proof}

\begin{example}
By Theorem~\ref{sth}, the Leavitt algebra $\LL(1,n)$ (and so by \cite[Corollary~2.10.8]{grrings} the matrix algebra over $\LL(1,n)$) is strongly graded as the only vertex in its graph is on cycles. 
\end{example}

\begin{example}
By Theorem~\ref{sth}, the Leavitt path algebra of the following graph is strongly graded. 
$$\xymatrix@=13pt{
            &    & \bullet  \ar[d] &  &  \bullet  \ar[d] & \\ 
\ar@{.}[r] & {\bullet} \ar[r] & {\bullet} \ar[r] & {\bullet}\ar[r] & {\bullet}\ar[r] & {\bullet} \ar@{.}[r] & \\ 
           &   \bullet \ar[u]   &  & \bullet \ar[u] &  &  \bullet \ar[u] & }$$
\end{example}

\begin{remark}
For the converse of Theorem~\ref{sth}, we have the following statement: If $\LL(E)$ is strongly graded, then for any vertex $v$ and any natural number $n \in \mathbb N$, there exists vertices $w_1,w_2,w_3$, and paths $v\geq w_1$ of length $n$, and $v\geq w_2$ of length $s$, and $w_3\geq w_2$ of length $s+n$. Although this formulation does not seem to look elegant, when the graph is finite, we have a simple characterization (see Theorem~\ref{sthfin}). 
\end{remark}
 
When the graphs have a finite number of vertices, i.e., their associated Leavitt path algebras are unital, we can give a complete characterization of strongly graded Leavitt path algebras. The following theorem shows that a Leavitt path algebra associated to a finite graph is strongly graded if and only if the graph does not have a sink. 

\begin{theorem}\label{sthfin}
Let $E$ be a finite graph. The Leavitt path algebra $\LL_R(E)$ with coefficients in a ring $R$ is strongly graded if and only if any vertex is connected to a cycle.
\end{theorem}
\begin{proof}
If any vertex is connected to a cycle then the orbit of any vertex is nonempty and bounded (as the graph is finite), thus by Theorem~\ref{sth}, $\LL(E)$ is strongly graded. For the converse, let $\LL(E)$ be strongly graded.
Then the graph $E$ does not have any sink (see Lemma~\ref{jh543p}(4)). Let $|E^0|=n$. For any vertex $v$, consider the path 
$\mu_1\mu_2\dots\mu_{n}$ of length  $n$ emitting from $v$ (this is possible as there is no sink). Since the number of vertices are $n$, this forces $s(\mu_i)=s(\mu_j)$ for some $i,j$. That is $v$ is connected to $\mu_i\mu_{i+1}\dots\mu_j$ which is a closed path. Now the following easy argument, based on an induction on the length, shows that any vertex on a closed path is connected to a cycle. If the vertex $w$ is on a closed path of length $1$, then this a loop and there is nothing to prove. Assume the statement is correct for any vertex on a closed path of length less than $k$. Let 
$\gamma_1\dots\gamma_k$ be a closed path of length $k$, with $s(\gamma_1)=r(\gamma_k)$.  If $s(\gamma_i)\not = s(\gamma_j)$ for all $1\leq i\not = j\leq k$, then the path is a cycle and there is nothing to prove. Otherwise, suppose $s(\gamma_i)=s(\gamma_j)$ and consider the path $\gamma_1\dots\gamma_{i-1}\gamma_j \dots \gamma_k$. This is clearly a closed path of smaller length and we are done by induction. 
\end{proof}

\begin{example}\label{opex}
According to Theorem~\ref{sthfin}, the graph $E$ below produces a strongly graded Leavitt path algebra, but the opposite graph $E^{\op}$ does not.  
\begin{equation*}
\xymatrix{
E: &   \bullet\ar@(lu,ld) & \bullet \ar[l] && E^{\op}:& \bullet\ar@(ld,lu)\ar[r] & \bullet 
}
\end{equation*}
 Also combining Theorem~\ref{sthfin} with~\cite[Theorem~11]{aap06}, it follows that unital purely infinite simple Leavitt path algebras are strongly graded. 
\end{example}

\section{Polycephaly graphs}
Two distinguished types of strongly graded Leavitt path algebras are $C_n$-comet graphs and multi-headed rose graphs (see Figure~\ref{monster2}). We consider  a polycephaly graph (Definition~\ref{popyt}) which is a mixture of these graphs (and so include all these types of graphs) and then determine the graded structure of Leavitt path algebras of polycephaly graphs in Theorem~\ref{polyhead}. 
Building on this, we also study graphs whose Leavitt path algebras are group rings or crossed-products.

\begin{deff}\label{cometi}
A finite graph $E$ is called a {\it $C_n$-comet}, if $E$ has exactly one cycle $C$ (of length $n$), and $C^0$ is dense, i.e., 
$T(v) \cap (C)^0 \not = \emptyset$ for any vertex $v \in E^0$. 
\end{deff} 

Note that the uniqueness of the cycle $C$ in the definition of $C_n$-comet together with its density implies that the
cycle has no exits.

\begin{deff}\label{wert}
A finite graph $E$ is called a {\it multi-headed comets} if $E$ consists of $C_{l_s}$-comets, $1\leq s\leq t$, of length $l_s$, such that cycles are mutually disjoint and any vertex is connected to at least a cycle. (Recall that the graphs in this paper are all connected).  More formally, $E$ consists of $C_{l_s}$-comets, $1\leq s\leq t$,  and for any vertex $v$ in $E$, there is at least a cycle, say,  $C_{l_k}$, such that $T(v) \cap C_{l_k}^0 \not = \emptyset$, and furthermore no cycle has an exit. 
\end{deff}

\begin{deff} \label{popyt}
A finite graph $E$ is called a {\it polycephaly} graph if $E$ consists of a multi-headed comets and/or an acyclic graph with its sinks attached to roses such that all the cycles and roses are mutually disjoint and any vertex is connected to at least a cycle or a rose. More formally, $E$ consists of $C_{l_s}$-comets, $1\leq s\leq h$, and a finite acyclic graph with sinks $\{v_{h+1},\dots,v_t\}$ together with  $n_s$-petal graphs $L_{n_s}$ attached to $v_s$, where  $n_s \in \mathbb N$ and $h+1\leq s \leq t$.  Furthermore any vertex $v$ in $E$ is connected to at least one of $v_s$, $h+1\leq s \leq t$, or  at least a cycle, i.e., there is  $C_{l_k}$, such that $T(v) \cap C_{l_k}^0 \not = \emptyset$, and  no cycle or a rose has an exit (see Definition~\ref{petaldef}). When $h=0$, $E$ does not have any cycle, and when $t=h$, $E$ does not have any roses.  
\end{deff}

\begin{remark}
Note that a cycle of length one can also be considered as a rose with one petal. This should not cause any confusion in the examples and theorems below. In some proofs (for example proof of Theorem~\ref{polyhead}), we collect all the cycles of length one in the graph as comet types and thus all the roses have either zero or more than one petals, i.e.,  $n_s=0$ or $n_s>1$, for any $h+1\leq s \leq t$. 
\end{remark}

\begin{example}\label{popyttr}
 Let $E$ be a polycephaly graph. 
\begin{enumerate}
\item If $E$ contains no cycles, and for any rose $L_{n_s}$ in $E$, $n_s=0$, then $E$ is a finite acyclic graph. 

\item If $E$ consists of only one cycle (of length $n$) and no roses, then $E$ is a $C_n$-comet. 

\item If $E$ contains no roses, then $E$ is a multi-headed comets. 

\item If $E$ contains no cycles, and for any rose $L_{n_s}$ in $E$, $n_s \geq 1$, then $E$ is a multi-headed rose graphs (see Figure~\ref{monster2}). 
\end{enumerate}
\end{example}

\begin{example}
The following graph is a (three-headed rose) polycephaly graph.

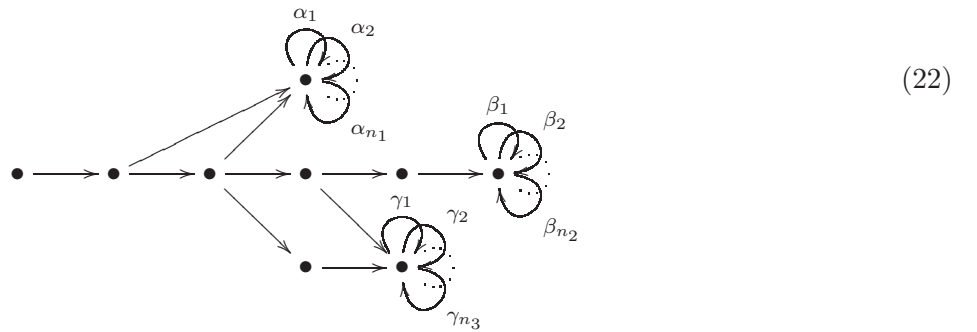
\begin{equation}\label{monster2}
\xymatrix{
 & & &&   \bullet  \ar@(ul,ur)^{\alpha_{1}}  \ar@(u,r)^{\alpha_{2}} \ar@{.}@(ur,dr) \ar@(r,d)^{\alpha_{n_1}}& \\
   &   \bullet \ar[r]  &  \bullet \ar[r] \ar[urr]   &  \bullet \ar[r] \ar[dr]  \ar[ur] & \bullet \ar[r]  \ar[dr] &    \bullet \ar[r]  & \bullet  \ar@(ul,ur)^{\beta_{1}}  \ar@(u,r)^{\beta_{2}} \ar@{.}@(ur,dr) \ar@(r,d)^{\beta_{n_2}}& \\
 &  &   &    & \bullet \ar[r]  &   \bullet  \ar@(ul,ur)^{\gamma_{1}} \ar@(u,r)^{\gamma_{2}} \ar@{.}@(ur,dr) \ar@(r,d)^{\gamma_{n_3}}&  
 }
 \end{equation}
By Theorem~\ref{sthfin}, this is a strongly graded ring. 
\end{example}

\begin{example}
The following graph is a (five-headed) polycephaly  graph, with two roses, two sinks and a comet.  
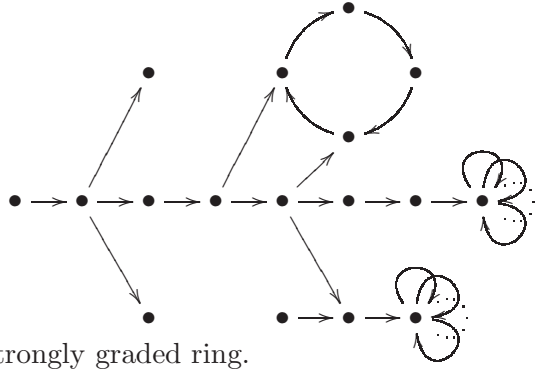
\begin{equation}\label{monster}
\xymatrix@=13pt{
& & & & & \bullet \ar@/^/[dr]      \\
& &  \bullet  & & \bullet \ar@/^/[ur] &&\bullet \ar@/^/[dl] &  \\
& &  &  & & \bullet \ar@/^/[ul]      \\
 \bullet \ar[r]  &   \bullet \ar[r]  \ar[uur] \ar[ddr] &   \bullet \ar[r]  &  \bullet \ar[r]   \ar[ruu] &  \bullet \ar[r] \ar[ddr] \ar[ru] & \bullet \ar[r]  &    \bullet \ar[r]  & \bullet  \ar@(ul,ur)  \ar@(u,r) \ar@{.}@(ur,dr) \ar@(r,d)& \\
  \\
& & \bullet  &   &  \bullet \ar[r]  & \bullet \ar[r]  &   \bullet  \ar@(ul,ur) \ar@(u,r) \ar@{.}@(ur,dr) \ar@(r,d)&  
 }
 \end{equation}
 By Theorem~\ref{sthfin}, this is not a strongly graded ring. 
\end{example}

We are in a position to classify (the graded structure of) polycephaly graphs. 

\begin{theorem}\label{polyhead}
Let $R$ be a  ring and    $E$ be a polycephaly graph consisting of cycles $C_{l_s}$, $1\leq s\leq h$, of length $l_s$
and an acyclic graph  with sinks $\{v_{h+1},\dots,v_t\}$ which are attached to $L_{n_{h+1}},\dots,L_{n_t}$, respectively.
For any $1\leq s\leq h$ choose $v_s$ (an arbitrary vertex) in $C_{l_s}$ and remove the edge $c_s$ with $s(c_s)=v_s$ from the cycle $C_{l_s}$.  Furthermore, for any $v_s$, $h+1\leq s \leq t$, remove the rose $L_{n_s}$ from the graph. 
In this new acyclic graph $E_1$, let $\{p^{v_s}_i \mid 1\leq i \leq n(v_s)\}$ be the set of all paths which end in $v_s$, $1\leq s \leq t$.
 Then there is a $\mathbb Z$-graded isomorphism
\begin{equation}\label{dampai}
\LL_R(E) \cong_{\gr}  \bigoplus_{s=1}^h  \M_{n(v_s)}\big(R[x^{l_s},x^{-l_s}]\big)(\overline{p_s}) \bigoplus_{s=h+1}^t \M_{n(v_s)} \big(\LL_R(1,n_s)\big) (\overline{p_s}),
\end{equation}
where $\overline{p_s}=\big(|p_1^{v_s}|,\dots, |p_{n(v_{s})}^{v_s}|\big)$.
\end{theorem}

\begin{remark}\label{jgfdfg}
Theorem~\ref{polyhead} shows that the Leavitt path algebra of a polycephaly graph decomposes into direct sum of three types of algebras. Namely, matrices over the field $K$ (which corresponds to acyclic parts of the graph, i.e, $n_s=0$ in~(\ref{dampai})), matrices over $K[x^l,x^{-l}]$, $l \in \mathbb N$, (which corresponds to comet parts of the graph) and matrices over Leavitt algebras $\LL(1,n_s)$, $n_s \in \mathbb N$ and $n_s \geq 1$ (which corresponds to rose parts of the graph). Also note that a cycle of length one can also be considered as a rose with one petal, which in either case, on the level of Leavitt path algebras, we obtain matrices over the algebra $K[x,x^{-1}]$. 
\end{remark}

\begin{proof}[Proof of Theorem~\ref{polyhead}] Let $C_{l_1},\dots, C_{l_h}$ be the cycles in the graph $E$. (By definition $h$ could be zero, i.e., $E$ contains no cycle, see Example~\ref{popyttr}.) And $\LL(1,n_{h+1}),\dots, \LL(1,n_t)$ be the roses in $E$ attached to $\{v_{h+1},\dots,v_t\}$, respectively. (Again $t$ could be $h$, i.e., $E$ contains no roses, or $n_s=0$, i.e., $v_s$ is a sink in $E$, see Definition~\ref{petaldef}.) 

Throughout the proof, we consider a cycle of length one in the graph as a comet and thus all $n_s>1$, $h+1\leq s \leq t$  (see Remark~\ref{jgfdfg}).
 Write $C_{l_s}=c_1^s\dots c_{l_s}^s$, $1\leq s\leq h$, where $c$'s are the edges constitute the cycle $C_{l_s}$, $C_{l_s}^0=\{v_1^s,\dots,v_{l_s}^s\}$, where $r(c_k^s)=v_k^s$ for all $k$, $s(c_1^s)=v_{l_s}^s$ and $s(c_k^s)=v_{k-1}^s$ for all $k\geq 1$.  For simplicity, throughout the proof, we denote the vertex $v_{l_s}^s$ by $v_s$, $1\leq s \leq h$. 

Remove all the edges $c_1^s$, $1\leq s\leq h$ from the cycles $C_{l_s}$ in the graph $E$. So $v_s=v_{l_s}^s$, $1\leq s\leq h$ becomes sinks. Furthermore, remove all the roses $\LL(1,n_{s})$, $h+1\leq s\leq t$ from this graph, and call the acyclic graph obtained in this way with sinks $\{v_1,\dots,v_h,v_{h+1},\dots,v_t  \}$ by $E_1$.   Define, for $1\leq s\leq h$, 
\[\I_{v_s}=\Big \{ \sum d \alpha C_{l_s}^k  \beta^* \mid d \in R, k\in \mathbb Z, \, \alpha, \beta \in E_1^*, r(\alpha)=v_{s}=r(\beta) \Big \} \subseteq \LL_R(E),\]
and for $h+1\leq s\leq t$,    
\[\I_{v_s}=\Big \{ \sum d \alpha x  \beta^* \mid d \in R, \, \alpha, \beta \in E_1^*, r(\alpha)=v_s=r(\beta), \, x\in \LL(1,{n_s}) \Big \} \subseteq \LL_R(E).\]

Note that since there are  complete graph homomorphisms from $E_1$, $C_{l_s}$, $1\leq s\leq h$ and $L_{n_s}$, $h+1\leq s\leq t$ to $E$ which induce  graded homomorphisms on the Leavitt path algebras level, in $\I_{v_s}$ we are actually working with the images of $\LL_R(E_1)$, $\LL_R(C_s)$ and  $\LL_R(1,{n_s})$ in $\LL_R(E)$. (Note that, by the graded uniqueness theorem (\cite[Theorem~4.8]{tomforde}), this homomorphisms are injective.) 

We first observe that  $\I_{v_s}$, $1\leq s\leq t$, are (two-sided) ideals of $\LL_R(E)$. It is enough to check that for monomials 
 $\alpha C_{l_s}^k \beta^* \in \I_{v_s}$,  $1\leq s\leq h$, $k\in \mathbb Z$ and $\alpha x \beta^* \in \I_{v_s}$, where $x\in \LL(1,n_s)$, $h+1\leq s\leq t$ and any monomial $\gamma \delta^*\in \LL_R(E)$,  $\gamma \delta^* \alpha C_{l_s}^k  \beta^*$ and  $\alpha C_{l_s}^k  \beta^* \gamma \delta^*$ are in $ \I_{v_s}$ and similarly  $\gamma \delta^* \alpha x  \beta^*$ and $\alpha x  \beta^* \gamma \delta^*$ are in $ \I_{v_s}$.

We first consider the case of roses, i.e.,  $x\in \LL(1,n_s)$ for $h+1\leq s\leq t$. Let $\gamma \delta^*$ be a nonzero monomial of $\LL_R(E)$ and $\alpha x \beta^* \in \I_{v_s}$. If $\gamma \delta^* \alpha x \beta^*\not = 0$, then either $\alpha=\delta p$ or $\delta=\alpha q$ for some paths $p$ and $q$ in $E$. 
Let us consider the first case $\alpha=\delta p$: Since $x\in \LL(1,{n_s})$ (recall that by our convention, $x$ is in fact in the image of  $\LL(1,{n_s})$ in $\LL(E)$), we know $r(\alpha)=r(p)=v_s$. If $s(p)=r(\delta)\not = v_s$, then $r(\gamma)=r(\delta)=s(p)\not = v_s$ and $\gamma p$ is a path in $E$. Furthermore, $\gamma p$ does not contain any vertices of cycles $C_{l_1},\dots,C_{l_h}$, neither any vertices of roses $L_{n_{h+1}},\dots, L_{n_{s-1}},  L_{n_{s+1}},\dots, L_{n_t}$, otherwise they would have an exit edge which can't be the case (see Definition~\ref{popyt} of polycephaly graphs). Therefore $\gamma p$ is a path in $E_1$. 
So $\gamma \delta^* \alpha x \beta^*=\gamma \delta^* \delta p x \beta^*=\gamma p x \beta^* \in \I_{v_s}$ as 
$\gamma p$ is a path in $E_1$ with $r(\gamma p)=v_s$. On the other hand, if $s(p)=v_s$, since $r(p)=v_s$ and $E_1$ is acyclic graph, $p=v_s$ and so $\alpha=\delta$. Therefore $r(\gamma)=r(\delta)=v_s$.  Repeating the same argument above,  $\gamma$ has to be a combination of a  path in $E_1$ which ends in $v_s$ and possibly continues in $L_{n_s}$, say $\gamma=q y$, where $q$ is a path in $E_1$ and $y $ a path in $L_{n_s}$. Therefore $\gamma \delta^* \alpha x \beta^*=\gamma \delta^* \delta x \beta^*=\gamma  x \beta^*=q (y x) \beta^* \in \I_{v_s}$

 In the second case that $\delta=\alpha q$, we have $\gamma \delta^* \alpha x \beta^*=\gamma q^* \alpha^* \alpha x \beta ^*=\gamma q^* x \beta ^*$. Since $\alpha$ is a path in $E_1$ which ends in $v_s$, $q$ needs to be a path in $L_{n_s}$. Therefore $r(\gamma)=r(q)=v_s$. So, as in the first case, $\gamma$ is a combination of a  path in $E_1$ which ends in $v_s$ and possibly continues in $L_{n_s}$. This shows that $\gamma q^* x\beta^* \in \I_{v_s}$. Therefore $\I_{v_s}$, $h+1\leq s\leq t$, is a left ideal of $\LL(E)$. In the same manner $\I_{v_s}$ is also a right ideal. 
 
We next consider the case of cycles, namely, let  $\gamma \delta^*$ be a nonzero monomial of $\LL_R(E)$ and 
$\alpha C_{l_s}^k  \beta^* \in \I_{v_s}$, for some $1\leq s \leq h$ and $k\in \mathbb Z$. Again, if $\gamma \delta^* \alpha C_{l_s}^k \beta^*\not = 0$, then either $\alpha=\delta p$ or $\delta=\alpha q$ for some paths $p$ and $q$ in $E$. 
We first consider the case $\alpha=\delta p$: We have $r(\alpha)=r(p)=s(C_{l_s})=v_{l_s}^s=v_s$. Three cases need to be considered: $s(p)=v_i^s$, $i\not=l_s$ (unless $l_s=1$); $s(p)=v_{l_s}^s$ and $s(p)  \in E_1^0 \backslash \cup_{s=1}^h C_{l_s}^0 \cup \{v_{h+1},\dots,v_t  \}$. For simplicity, throughout we denote the edges $c_i^s$ by $c_i$. 

Suppose $s(p)=v_i^s$, $i\not=l_s$, i.e., $p=c_{i+1}\dots c_{l_s}$.  We have $r(\gamma)=r(\delta)=s(p)=v_i^s$. Then $\gamma$ can be written as 
$\gamma=q C_{v_{i'}}^{k''}c_{i'+1}\dots c_i$, where $q$ is a path in $E_1$, $C_{v_{i'}}$ is a cycle with $s(C_{v_{i'}})=v_{i'}$, $k'' \in \mathbb N$, ($k''$ could be zero)  and $c_{i'+1}\dots c_i$ is a unique path in $C_{l_s}$ which does not contain $C_{v_i'}$, connecting $v_{i'}^s$ to $v_i^s$. 

If $v_{i'}^s=v_{l_s}^s=v_s$, then $\gamma=q C_{v_{i'}}^{k''}c_{1}\dots c_i$ and so 
\[\gamma \delta^* \alpha C_{l_s}^k  \beta^*=\gamma \delta^* \delta p C_{l_s}^k  \beta^*=\gamma p C_{l_s}^k  \beta^*=q C_{l_s}^{k''}c_{1}\dots c_i c_{i+1}\dots c_{l_s} C_{l_s}^k  \beta^*=q C_{l_s}^{k''+k+1}\beta^* \in \I_{v_s}.\] 

If the path $c_{i'+1}\dots c_i$ in $\gamma$ contains $c_1$ then 
\begin{equation}\label{daijoon}
\gamma \delta^* \alpha C_{l_s}^k  \beta^*=q C_{v_{i'}}^{k''}c_{i'+1}\dots c_{l_s} c_1 \dots  c_{i}c_{i+1}\dots c_{l_s} C_{l_s}^k  \beta^*=q C_{v_{i'}}^{k''}c_{i'+1}\dots c_{l_s}  C_{l_s}^{k+1}  \beta^*  .
\end{equation}
Note that $C_{v_{i'}}c_{i'+1}\dots c_{l_s}=c_{i'+1}\dots c_{l_s} C_{l_s}$. Plugging this in ~(\ref{daijoon}), we get 
\[\gamma \delta^* \alpha C_{l_s}^k  \beta^*=q c_{i'+1}\dots c_{l_s} C_{l_s}^{k''+k+1} \in \I_{v_s},\] as $q c_{i'+1}\dots c_{l_s}$ is a path in $E_1$. 

If the path $c_{i'+1}\dots c_i$ in $\gamma$ does not contains $c_1$ then 
\begin{equation}\label{daij}
\gamma \delta^* \alpha C_{l_s}^k  \beta^*=q C_{v_{i'}}^{k''}c_{i'+1}\dots c_{i} c_{i+1}\dots c_{l_s} C_{l_s}^k  \beta^*.
\end{equation}
Note that $C_{v_{i'}}c_{i'+1}\dots c_{i} c_{i+1}\dots c_{l_s}=c_{i'+1}\dots c_{i} c_{i+1}\dots c_{l_s} C_{l_s}$.  Plugging this into ~(\ref{daij}), we get
 \[
\gamma \delta^* \alpha C_{l_s}^k  \beta^*=q c_{i'+1}\dots c_{i} c_{i+1}\dots c_{l_s} C_{l_s}^{k''+k}  \beta^* \in \I_{v_s}.
\] as $q c_{i'+1}\dots c_{i} c_{i+1}\dots c_{l_s}$ is a path in $E_1$.

Suppose $s(p)=v_s$. Then  $p=v{_s}$, since $p$ is acyclic. Thus $r(\gamma)=s(p)=v_{s}$.  As in the previous case, $\gamma$ can be written as $\gamma=q C_{l_s}^{k'}$, where $q$ is a path in $E_1$.  Thus 
\[\gamma \delta^* \alpha C_{l_s}^k  \beta^*=\gamma \delta^* \delta p C_{l_s}^k  \beta^*=
\gamma  C_{l_s}^k  \beta^*= q C_s^{k'} C_{l_s}^k  \beta^*= q C_s^{k'+k} \beta^* \in \I_{v_s}.\] 

Finally, suppose $s(p)  \in E_1^0 \backslash \cup_{s=1}^h C_s^0 \cup \{v_{h+1},\dots,v_t  \}$. Then clearly $\gamma$ has to be a path in $E_1$, and so $\gamma p$ is a path in $E_1$ as well. So $\gamma \delta^* \alpha C_{l_s}^k  \beta^*=\gamma p  C_{l_s}^k  \beta^* \in  \I_{v_s}$. 

Next we consider the case $\delta=\alpha q$. Then $s(q)=r(\alpha)=v_{s}$. Since $q$ is a path with the source $v_{s}$ (i.e., on the cycle $C_s$), it has to be of the form $q=C_{l_s}^{k'}c_1\dots c_i$, where $k' \in \mathbb N$. Thus $r(\gamma)=r(\delta)=r(q)=v_i^s$. Therefore $\gamma$ (as in the previous case) can be written as $\gamma=p C_{v_{i'}}^{k''}c_{i'+1}\dots c_i$, where $p$ is a path in $E_1$, $C_{v_{i'}}$ is a cycle with $s(C_{v_{i'}})=v_{i'}$. So 
\begin{equation}\label{testshh}
\gamma \delta^* \alpha C_{l_s}^k  \beta^*=\gamma q^* \alpha ^* \alpha C_{l_s}^k  \beta^*=\gamma q^* C_{l_s}^k  \beta^*. 
\end{equation}
But $\gamma q^*=p C_{v_{i'}}^{k''}c_{i'+1}\dots c_i c_i^*\dots c_1^* C_{l_s}^{-k'}$. 
If the lengths of the paths $c_{i'+1}\dots c_i$ and $c_1\dots c_i$ are the same, then $v_{i'}=v_{s}$. Thus
$\gamma q^*=p C_{l_s}^{k''} C_{l_s}^{-k'}=pC_{l_s}^{k''-k'}$. Plugging this into~(\ref{testshh}), we have 
\[\gamma \delta^* \alpha C_{l_s}^k  \beta^* =pC_{l_s}^{k''-k'}C_{l_s}^k  \beta^* = pC_{l_s}^{k''-k'+k}\beta^*.\] Since $p$ is a path in $E_1$, it follows $pC_{l_s}^{k''-k'+k}\beta^* \in \I_{v_s}$. 

If the length of the path $c_{i'+1}\dots c_i$ is greater than   the length of $c_1\dots c_i$, then writing $c_{i'+1}\dots c_i=c_{i'+1}\dots c_{l_s} c_1\dots  c_i $ we have  
\begin{equation}\label{hgpoolad}
\gamma q^*=p C_{v_{i'}}^{k''}c_{i'+1}\dots c_{l_s} C_{l_s}^{-k'}.
\end{equation} 
  Otherwise 
\begin{equation}\label{hgpoolad2}
\gamma q^*=p C_{v_{i'}}^{k''}c_{i'+1}\dots c_i c_i^*\dots c_{i'+1}^* c_{i'}^*\dots c_1^* C_{l_s}^{-k'} =p C_{v_{i'}}^{k''}c_{i'}^*\dots c_{1}^* C_{l_s}^{-k'}.
\end{equation}
Note that $C_{v_{i'}}c_{i'+1}\dots c_{l_s}=c_{i'+1}\dots c_{l_s}C_{l_s}$ and $C_{v_{i'}}c_{i'}^*\dots c_1^*=c_{i'}^*\dots c_1^* C_{l_s}$.  Plugging these into Equations~\ref{hgpoolad} and~\ref{hgpoolad2}, respectively, we get 
$\gamma q^* =p c_{i'+1}\dots c_{l_s} C_{l_s}^{k''-k'}$ and \[\gamma q^* =p c_{i'}^*\dots c_1^* C_{l_s}^{k''-k'}=p  c_{i'+1}\dots c_{l_s} C_{l_s}^* C_s^{k''-k'}= p  c_{i'+1}\dots c_{l_s} C_{l_s}^{k''-k'-1}.\] Again plugging these two identities into~(\ref{testshh}), we have 
\[\gamma \delta^* \alpha C_{l_s}^k  \beta^*=p c_{i'+1}\dots c_{l_s} C_{l_s}^{k''-k'+k}\beta^*\] and 
\[\gamma \delta^* \alpha C_{l_s}^k  \beta^*=p  c_{i'+1}\dots c_{l_s} C_{l_s}^{k''-k'+k-1}\beta^*.\] Since $p c_{i'+1}\dots c_{l_s} $ is a path in $E_1$ with the range $v_s$, both of these equations are in $\I_{v_s}$.  
Therefore $\I_{v_s}$ $1\leq s\leq l$, is a left ideal of $\LL(E)$. In the same manner $\I_{v_s}$ is also a right ideal.

 Next we show that $\I_{v_s}\cong \M_{n(v_s)}(\LL_R(1,n_s))$, for $h+1\leq s\leq t$,  and  $\I_{v_s}\cong \M_{n(v_{s})}(R[x^{l_s},x^{-l_s}])$, where $l_s=|C_{l_s}|$, for $1\leq s\leq h$. Here $n(v_s)$  is the number of all the paths in $E_1$ which end in $v_s$.  
 
Rename the paths in $E_1$ which end in $v_s$, $h+1\leq s\leq t$, and consider the set  $P_{v_s}=\{p_i^{v_s} \mid 1\leq i \leq n(v_s)\}$  of all such paths. Throughout this paragraph, since we are working with a fixed $s$, we use $p_i$ instead of $p_i^{v_s}$. Thus $\I_{v_s}= \big \{\sum d p_i x p_j^* \mid d \in R, p_i,p_j \in P_{v_s}, x \in \LL(1,{n_s}) \big \}$. First observe that 
 \begin{equation}\label{mopo}
 p_j ^* p_t=\delta_{i t} v_s.
 \end{equation} This is clear if $t=j$. If $j \not = t$, then this product is zero unless, $p_t=p_j q$ or $p_j=p_t q'$, which is not possible, since $p_i$'s are coming from $E_1$ and $v_s$ is a sink in $E_1$. 
 
 Define a map $\phi: \I_{v_s} \rightarrow \M_{n(v_s)}(\LL_R(1,n_s))$, by sending 
 $p_ix p_j^*$ to $e_{ij}(x)$ and extend it linearly. (Note that if $n_s=0$, i.e, $L_{n_s}=v_s$, then $x=v_s$ and $
 \phi(p_iv_s p_j^*)=\phi(p_i p_j^*)=e_{ij}(1)$.) 
 This map is well-defined. Because, if 
 \[\sum_{i,j}  d_{ij} p_i x_{ij} p_j^*=\sum_{i,j}  d'_{ij} p_i x_{ij}' p_j^*,\] then fixing $1\leq m,n \leq n(v_s)$, and multiplying this equation from the left by $p_m^*$ and from the right by $p_n$ and using the identity~(\ref{mopo}) we get $d_{mn}x_{mn}=d'_{mn}x'_{mn}$. This shows that 
\[\phi(\sum_{i,j}  d_{ij} p_i x_{ij} p_j^*)=\sum_{i,j}e_{ij}(d_{ij}x_{ij})=\sum_{i,j}e_{ij}(d'_{ij}x'_{ij})=\phi(\sum_{i,j}  d'_{ij} p_i x_{ij}' p_j^*).\] So $\phi$ is well-defined.  Again using identity~(\ref{mopo}) we have $(p_i x p_j ^*)( p_f x' p_g^*)=\delta_{j f} p_i x x' p_g^*$. 
This shows that the map $\phi$ is a homomorphism. Now one can easily check that $\phi$ is bijective as well, so $\phi$ is in fact an isomorphism. Finally considering the grading $\M_{n(v_s)} (\LL_R(1,n_s))\big (|p^{v_s}_1|,\dots, |p^{v_s}_{n(v_s)}|)$, 
we check that $\phi$ is a graded map. (Recall that we denote $p_i^{v_s}$ by $p_i$.) Let $p_ixp_j^* \in \I_{v_s}$. Then $\deg(p_ixp_j^*)=\deg(x)+|p_i|-|p_j|$. Now by ~(\ref{hogr}), 
\[\deg(\phi(p_ixp_j^*))=\deg(e_{ij}(x))=\deg(x)+|p_i|-|p_j|\] which shows that $\phi$ respects the grading. 

Next we show that $\I_{v_s}\cong \M_{n(v_{s})}(R[x^{l_s},x^{-l_s}])$, where $l_s$ is the length of the cycle $C_{l_s}$. 
Rename the paths in $E_1$ which end in $v_{s}$, $1\leq s\leq h$, and consider the set  $P_{v_s}=\{p_i^s \mid 1\leq i \leq n(v_s)\}$  of all such paths. Throughout this paragraph, since we are working with a fixed $s$, we use $p_i$ instead of $p_i^{v_{s}}$. Thus $\I_{v_s}= \big \{\sum d p_i C_{l_s}^k p_j^* \mid d \in R, k\in \mathbb Z, p_i,p_j \in P_{v_s} \big \}$.  Similar to (\ref{mopo}), we have 
$p_j ^* p_t=\delta_{i t} v_{s}$.

Define a map $\phi: \I_{v_s} \rightarrow \M_{n(v_s)} (R[x^{l_s},x^{-l_s}])$, by sending 
 $p_iC_{l_s}^k p_j^*$ to $e_{ij}(x^{kl_s})$ and extend it linearly. (Note that if $l_s=|C_{l_s}|=0$, then $C_{l_s}=v_s$ and $
 \phi(p_iC_{l_s}^k p_j^*)=\phi(p_i p_j^*)=e_{ij}(1)$.) 
 This map is well-defined. Because, if 
 \[\sum_{i,j}  \sum_{k} d_{ijk} p_i C_{l_s}^{w_{ijk}} p_j^*=\sum_{i,j} \sum_{k} d'_{ijk} p_i C_{l_s}^{w'_{ijk}} p_j^*,\] where the inner sum on $k$ arranged in an increasing order on $w_{ijk}$ and $w'_{ijk}$, respectively, 
  then fixing $1\leq m,n \leq n(v_s)$, and multiplying this equation from the left by $p_m^*$ and from the right by $p_n$ and using the identity~(\ref{mopo}) we get 
 \[\sum_k d_{mnk} C_{l_s}^{w_{mnk}} = \sum_k d'_{mnk} C_{l_s}^{w'_{mnk}} .\] 
 Comparing the homogenous parts of each side of the equation, it immediately gives $w_{mnk}=w'_{mnk}$ and $d_{mnk} C_{l_s}^{w_{mnk}}= d'_{mnk} C_{l_s}^{w_{mnk}}$. Multiplying the last equation by ${C_{l_s}^{w_{mnk}}}^*$ we get $d_{mnk}v_{m_s}^s=d'_{mnk}v_{m_s}^s$. By~\cite[Proposition~3.4]{tomforde}, if $d v=0$, where $d \in R$ and $v \in E^0$, then $d=0$ (the proof of~\cite[Proposition~3.4]{tomforde} is written for a commutative ring $R$ with identity, however the same proof is valid for a non-commutative ring with identity).  This implies that $d_{mnk}=d'_{mnk}$. 
 This shows that 
\[\phi(\sum_{i,j}  \sum_{k} d_{ijk} p_i C_{l_s}^{w_{ijk}} p_j^*)=\sum_{i,j} \sum_k e_{ij}(d_{ijk} x^{l_sw_{ijk}})=\sum_{i,j} \sum_k e_{ij}(d'_{ijk} x^{l_sw'_{ijk}})=\phi(\sum_{i,j}  \sum_{k} d'_{ijk} p_i C_{l_s}^{w'_{ijk}} p_j^*).\] So $\phi$ is well-defined.  Again  we have $(p_i C_{l_s}^k p_j ^*)( p_f C_{l_s}^{k'} p_g^*)=\delta_{j f} p_i C_s^{k+k'} p_g^*$. 
This shows that the map $\phi$ is a homomorphism. Now one can easily check that $\phi$ is bijective as well, so $\phi$ is in fact an isomorphism. Finally considering the grading $\M_{n({v_s})} (R[x^{l_s},x^{-l_s}])\big (|p^{s}_1|,\dots, |p^{s}_{n({v_s})}|)$, 
we check that $\phi$ is a graded map. (Recall that we denote $p_i^{s}$ by $p_i$.) 
By definition, the degree of the homogeneous element $p_i C_{l_s}^k p_j^* \in \I_{v_s}$ is its length, i.e., $|p_i C_{l_s}^k p_j^*|=kl_s+|p_i|-|p_j|$ (note that $k \in \mathbb Z$). On the other hand, by ~(\ref{hogr}), \[\deg(\phi(p_i C_{l_s}^k p_j^*))=\deg(e_{ij}(x^{k l_s}))=k l_s+|p_i|-|p_j|.\] This shows that $\phi$  respects the grading and therefore it is a $\mathbb Z$-graded isomorphism.

We are left to show that \[\LL_R(E)=\bigoplus_{s=1}^t \I_{v_s}.\] 
Let $\alpha \beta ^* \in \LL_R(E)$ be a nonzero monomial.
We consider two possibilities: 

(i)  Let $r(\alpha)=v_{l_s}^s=v_s$ for $1\leq s \leq h$ (the source of the cycles), or  $r(\alpha)=v_s$ for $h+1\leq s \leq t$ (the source of the roses). It is easy to see that in the first case, $\alpha=p C_{l_s}^k$, where $p$ is a path in $E_1$ with $r(p)=v_s$, $k \in \mathbb N$ and in the second case $\alpha=p x$, where $p$ is a path in $E_1$ with $r(p)=v_s$ and $x \in \LL(1,{n_s})$. Since $r(\alpha)=r(\beta)$ it follows that either $\alpha=p_1 C_{l_s}^{k_1}$ and  $\beta=p_2 C_{l_s}^{k_2}$ and so $\alpha \beta^*=p_1 C_{l_s}^{k_1-k_2} p_2^* \in \I_{v_s}$ or $\alpha=p_1 x_1$ and  $\beta=p_2 x_2$ and so $\alpha \beta^*=p_1 x_1x_2^* p_2^* \in \I_{v_s}$

(ii) If $r(\alpha) \not \in \{v_{1},\dots,v_{h},v_{h+1},\dots,v_t  \}$, then $r(\alpha)$ is not a sink in $E_1$. Let $v:=r(\alpha) \in E_1^0$. Consider the set $\mathcal P= \{v_{1},\dots,v_{h},v_{h+1},\dots,v_t  \}$.  Note that $\mathcal P$ is dense in $E_1$, and since $E_1$ contains no cycle, the orbit of $v$ in $E_1$, $O_{\mathcal P}(v)$, consists of acyclic paths. Thus by Lemma~\ref{nicelem}, we have 
 \begin{equation*}
 v=\sum_{\gamma  \in O_{\mathcal P}(v)}\gamma r(\gamma) \gamma^*.
 \end{equation*}
  Thus, $\alpha\beta^*=\alpha v \beta^*=\sum \alpha \gamma r(\gamma) \gamma^* \beta ^*$. Note that 
  $r(\gamma) \in \{v_{1},\dots,v_{h},v_{h+1},\dots,v_t  \}$. Now by Part (i), each of $\alpha \gamma r(\gamma) \gamma^* \beta ^*$ is in some $\I_{v_s}$, $1\leq s \leq t$, so
  $\alpha\beta^* \in \sum_{s=1}^t \I_{v_s}$. 
Thus 
\begin{equation}\label{teh8Oct}
\LL_R(E)=\sum_{s=1}^t \I_{v_s}.
\end{equation}

To show that  Equation~\ref{teh8Oct} is a direct sum, let $\alpha x \beta^*  \in \I_{u}$ and $\gamma y \delta^* \in \I_{w}$,  where  $u\not =w$, 
$\{u,w\} \subseteq  \{v_{1},\dots,v_{h},v_{h+1},\dots,v_t  \}$,  and $\alpha, \beta$ and $\gamma, \delta$ be paths in $E_1$ with the range $u$ and $w$,  respectively. Since there is no path connecting $u$ to $w$ (and the other way around), there are no paths of the forms $\gamma=\beta q$ or $\beta=\gamma q'$, so $\alpha \beta^* \gamma \delta^*=0$, showing $\I_{u} \I_{w}=0$. Since $\LL_R(E)$ is unital, it now follows that~(\ref{teh8Oct}) is a direct sum.
%let $1=\sum_{s=1}^t \alpha_s$, where $\alpha_s \in \I_{v_s}$. If $\gamma \in \I_{v_i} \cap \sum_{s=1,  s\not = i}^t \I_{v_s}$, then $1 \gamma= \alpha_i \gamma =0$.  This establishes the isomorphism in ~(\ref{dampai}). 
\end{proof}

\begin{example} \label{nopain}
Consider the polycephaly graph $E$

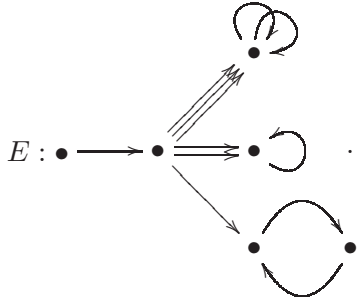
\[
 \xymatrix{
& & \bullet  \ar@(ul,ur)  \ar@(u,r) \\
E: \bullet \ar[r]  & \bullet   \ar@<1.5pt>[r]  \ar@<-1.5pt>[r] \ar@<0.5ex>[ur] \ar@<-0.5ex>[ur] \ar@<0ex>[ur] \ar[dr] & \bullet \ar@(rd,ru) &. \\ 
& &   \bullet \ar@/^1.5pc/[r] & \bullet \ar@/^1.5pc/[l]&  \\
 }
\]
\medskip 
Then by Theorem~\ref{polyhead}, 
\[\LL_K(E)\cong_{\gr} \M_5(K[x,x^{-1}])(0,1,1,2,2)  \oplus \M_4(K[x^2,x^{-2}])(0,1,1,2) \oplus \M_7(\LL(1,2))(0,1,1,1,2,2,2).\] By Theorem~\ref{sthfin}, this is a strongly graded ring. 
\end{example}
 
If there are no cycles and no petals in the graph (i.e., all $n_s=0$ in the definition), then the polycephaly graph reduces to a finite acyclic graph (Example~\ref{popyttr}) and thus 
 Theorem~\ref{polyhead} covers the Leavitt path algebra of acyclic graphs (see also Theorem~\ref{acyclicc}). For example, for an oriented $n$-line graph, we get  
\[\LL_K(\xymatrix@=13pt{ \bullet_1 \ar[r] &   \bullet_2 \ar@{.>}[r] &   \bullet_{n-1} \ar[r] & \bullet_n) } \cong_{\gr} \M_n(K)(0,1,\dots,n-1).\]

If a polycephay graph consists only of roses with nontrivial petals  (i.e., all $n_s \geq 1$), then  
by Theorem~\ref{sthfin}, the Leavitt path algebra of a  polycephaly graph is strongly graded. Furthermore we will see that if all the heads consist of only one loop, i.e., for any $s$, $L_{n_s}=L_1$, then its LPA is a group ring. 
For a simple case of polycephaly graph consisting of only one head, a special case of the following theorem shows that
\begin{align*}
&\xymatrix@=13pt{
\LL_K\Big (   \bullet_1 \ar[r]    &  \bullet_2 \ar@{.}[r]  &    \bullet_{m-1}  \ar[r] &  \bullet_m  \ar@(ul,ur)^{\alpha_{1}}  \ar@(u,r)^{\alpha_{2}} \ar@{.}@(ur,dr) \ar@(r,d)^{\alpha_{n}}
 }\Big ) \cong_{\gr} \\
& 
\xymatrix@=13pt{
  \LL_K\Big (   \bullet_1 \ar[r]    &  \bullet_2 \ar@{.}[r] & \bullet_{m-1} \ar[r]  &    \bullet_m   \Big ) \bigotimes_K   \LL_K\Big(\!\!\!\!\!&\! \! \! \! \bullet_1  \ar@(ul,ur)^{\alpha_{1}}  \ar@(u,r)^{\alpha_{2}} \ar@{.}@(ur,dr) \ar@(r,d)^{\alpha_{n}}
} \Big ) \cong_{\gr} \\
&
\M_m(K)(0,1,\dots,m-1)\otimes_K \LL_K(1,n) \cong_{\gr} \M_m(\LL_K(1,n))(0,1,\dots,m-1).
\end{align*}

This example (without taking into account the grading) is Proposition~13 in ~\cite{aap06} (see also~\cite[Lemma~5.1]{aalp}).

\begin{theorem}\label{polyheadnew}
Let $K$ be a  field and   $E$ be a polycephaly graph consisting of cycles $C_{l_s}$ for $1\leq s\leq h$ of length $l_s$
and an acyclic graph  with sinks $\{v_{h+1},\dots,v_t\}$ which are attached to $L_{n_{h+1}},\dots,L_{n_t}$, respectively. Let $E_1$ be the acyclic graph after removing the appropriate edges from $E$ as in Theorem~\ref{polyhead}. 
\begin{enumerate}[\upshape(1)]

\item If $E$ contains no cycle and  for $1\leq s \leq t$, $L_{n_s}=L_n$, then
\begin{equation}\label{xxdampai}
\LL_K(E)\cong_{\gr} \LL_K(E_1)\otimes_K \LL_K(1,n).
\end{equation}

\item  The Leavitt path algebra $\LL_K(E)$ is a group ring
\[\LL_K(E)\cong_{\gr} \bigoplus_{s=1}^t \M_{n(v_s)}(K)[\mathbb Z],\]
if and only if for any $1\leq s \leq h$, $|C_{l_s}|=1$ and for any $h+1\leq s \leq t$, $L_{n_s}=L_1$. 

\end{enumerate}
\end{theorem}

\begin{proof}

(1)  In order to prove (\ref{xxdampai}), first note that there are obvious complete graph homomorphisms from $E_1$ and $L_n$ to $E$ which in effect induce ring homomorphisms from $\LL_K(E_1)$ and $\LL_K(1,n)$ to $\LL_K(E)$, respectively.  For any sink $v_s$ in $E_1$, let $P_s=\{p^{v_s}_i \mid 1\leq i \leq n(v_s)\}$ be the set of all paths which end in $v_s$ and set $P=\bigcup_{1\leq s\leq t}P_s$. Using the  the proof of Theorem~\ref{polyhead} for the graph $E_1$, i.e., a polycephaly graph with all heads consisting of no loops, it follows that all nonzero $p_ip_j^*$, where $p_i, p_j \in P$ is a $K$-basis for $\LL_K(E_1)$ (clearly $p_ip_j^*\not = 0$ if and only if  $p_i, p_j \in P_s$ for some $1\leq s \leq t$). 
 
  Define $\phi: \LL_K(E_1) \otimes_K \LL_K(1,n) \rightarrow \LL_K(E)$ by $k_{ij}p_i p_j^* \otimes x  \mapsto p_i k_{ij} x p_j^*$ for $p_i,p_j\in P$, $k_{ij}\in K$, $x\in \LL_K(1,n)$
 and extend it linearly (here the images of $p_i\in \LL_K(E_1)$ and $x\in \LL_K(1,n)$ in $\LL_K(E)$ are defined by using the homomorphisms above). This is a well-defined map and one can easily check that it is also a $K$-algebra homomorphism which respects the grading (see~\ref{tengr}).  Injectivity of $\phi$ follows from the following general fact: If $A_s$, $s\in I$ and $B$  are rings with identities with (non-unital) ring homomorphisms $\phi_s:A_s \rightarrow B$, such that the images of $A_s$ in $B$ are mutually orthogonal, then the map $\textstyle{\bigoplus_s}\phi_s :\textstyle{\bigoplus_s}A_s \rightarrow B$ is a ring homomorphism and furthermore if $\phi_s$'s are all injective, then so is  $\textstyle{\bigoplus_s}\phi_s$. Now in our setting $A_s=M_{n(v_s)}(K)\otimes_K \LL_K(1,n)$, $B=\LL_K(E)$ with $\phi_s(e_{ij}(k_{ij})\otimes x)=p_i(k_{ij} x)p_j^*$. 
 Using part (1), observe that $\textstyle{\bigoplus_{1\leq s \leq t}}A_s \cong  \LL_K(E_1) \otimes_K \LL_K(1,n)$ and  $\phi=\textstyle{\bigoplus_i}\phi_i$. Since $\LL_K(1,n)$ is simple (\cite{vitt65,aap05}), $A_s$ is simple (see~\cite[IX, Theorem~6.2]{hungerford}) and so all $\phi_s$'s are injective. So $\phi$ is injective. 
 It is easy to see that $\phi$ is surjective as well and so we are done. (In fact one does not need to use the simplicity of the algebras to prove that $\phi_s$'s are injective. One can easily establish the injectivity, using the fact that $\M_n(K)$ and $\LL_K(1,n)$ are $K$-algebras. This is important, as in \S\ref{secitwo}  we carry over this proof to the case of weighted Leavitt path algebras where these algebras are not necessarily simple.)
 
 (2) Suppose that  for any $1\leq s \leq h$, $|C_{l_s}|=1$ and $h+1\leq s \leq t$, $L_{n_s}=L_1$.  Since a cycle of length one is a loop, we thus have $C_{l_s}=L_1$ as well.  Since $\LL_K(1,1)\cong K[x,x^{-1}]$ by Part~(1) we get
\[ \LL_K(E)\cong_{\gr} \LL_K(E_1)\otimes_K \LL_K(1,1) \cong_{\gr} 
 \bigoplus_{s=1}^t \M_{n(v_s)}(K)\otimes_K K[x,x^{-1}] \cong_{\gr}
 \bigoplus_{s=1}^t \M_{n(v_s)}(K)[\mathbb Z],\] with the grading as described in~(\ref{hgogt}).

Conversely, suppose that $\LL_K(E)$ is a group ring over $\mathbb Z$. Then we shall show that for any $1\leq s \leq h$, $|C_{l_s}|=1$ and 
$h+1\leq s \leq t$, $L_{n_s}=L_1$.   Recall the grading on a group ring from~(\ref{hgogt}). In our setting the grade group is $\mathbb Z$. Thus  there is an invertible element  of degree $1$ which commutes with all elements of the ring.  
 So let $x\in \LL_K(E)_1$ be an invertible element which is in the center of $\LL_K(E)$.   Using the graded isomorphism of~(\ref{dampai}), one can write $x=(x_1,\dots,x_h,x_{h+1},\dots,x_r)$,  where $x_s \in \M_{n(v_s)}(K[x^{l_s},x^{-l_s}])$ for $1\leq s \leq h$ and 
 $x_s \in \M_{n(v_s)}(\LL_K(1,n_s))$ for  $h+1\leq s \leq t$ and $\deg(x_s)=1$, $1\leq s \leq t$. Since $x$ is in the center of $\LL_K(E)$, all the  $x_s$'s are in the center of 
$\M_{n(v_s)}(K[x^{l_s|},x^{-l_s}])$ and  $\M_{n(v_s)}(\LL_K(1,n_s))$, respectively. But if there is an $s$, $h+1\leq s \leq t$, such that $n(v_s)>1$, then 
 the center of $\M_{n(v_s)}(\LL_K(1,n_s))$ is $KI_{n(v_s)}$ (see~\cite[Theorem~4.2]{acrow}).   Clearly $KI_{n(v_s)} \subseteq \M_{n(v_s)}(\LL_K(1,n_s))_0$. This prevents $x_s$ to be a homogeneous element of degree $1$ in $\LL_K(E)$. Thus $n(v_s)=1$, $h+1\leq s \leq t$. 
 
On the other hand suppose $|C_{l_s}|>1$ for some $1\leq s \leq h$.  
Denote $v_s$ by $u$ and as usual consider the set $P_{u}=\{p \mid p \text{ is a path in } E_1 \text{ with }  r(p)=u\}$.   Then by the proof of Theorem~\ref{polyhead},  
\[\I_{v_s}= \big \{\sum d p_i C_{l_s}^k p_j^* \mid d \in K, k\in \mathbb Z, p_i,p_j \in P_{u} \big \}=\M_{n(u)}(K[x^{l_s},x^{-l_s}]).\]

Thus $x_s \in \I_{v_s}$ which is invertible in $\I_{v_s}$ and commutes with all the elements of zero degree.  In particular, $x_su=ux_s$. Writing $x_s$ as a sum of monomials of the form $p_i C_{l_s}^k p_j^*$ (which are linearly independent, see also the proof of Theorem~3.3 in~\cite{ajis}) this implies that any monomial in $x$ which starts from $u$ also ends in $u$. Since $x_sx_s^{-1}=1_{\I_{v_s}}$, there is a monomial in $x$, say $p_iC_{l_s}^kp_j^*$, which does start from $u$ (If not then $0=ux_sx_s^{-1}=u$ which is a contradiction). Thus this monomial should also ends in $u$.  This forces that $p_i=p_j=u$ and therefore $C_{l_s}^k$ is a monomial in $x_s$. But then the monomial $C_{l_s}^k$ could not be in $\LL_K(E)_1$ as the length of $C_{l_s}$  is greater than $1$. This is a contradiction.  So $|C_{l_s}|=l_s=1$ for any $1\leq s \leq h$. This finishes the proof.
\end{proof}

\begin{example}
Let $K$ be a field. By Theorem~\ref{polyheadnew}  the Leavitt path algebra of the following graph is the group ring   
$\M_2(K)^7[\mathbb Z]$, 
\[
\LL_K\Big ( \; \xy 
/r3pc/: {\xypolygon7{~<<{}~>{}{\circlearrowright}}}
\endxy \; \Big ) \cong_{\gr} \M_2(K)^7[\mathbb Z],
\]
where $\M_2(K)^7$ stands for $7$ copies  of $\M_2(K)$, with the shifting $(0,1)$, with the grading as described in~(\ref{hgogt}).

\subsection{Complete classification of acyclic and multi-headed comet graphs}
In~\cite{aap07}  finite dimensional Leavitt path algebras (which are LPA's associated to finite acyclic graphs) and  in \cite{ajis} locally finite just infinite Leavitt path algebras (which are LPA's associated to $C_n$-comets) are characterized. However in both cases, the natural $\mathbb Z$-grading of Leavitt path algebras are not taken into account. Here,  building on their method, we characterize the Leavitt path algebras of these graphs, taking into account their natural grading, thus uncovering the natural structure that the Leavitt path algebras carry (see Introduction \S\ref{introf}).  

The following theorem characterizes the $\mathbb Z$-graded structure of acyclic Leavitt path algebras. 
\begin{theorem}\label{acyclicc}
Let $R$ be a ring and $E$ be a finite acyclic graph with sinks $\{v_1,\dots,v_t\}$. For any sink $v_s$, let $\{p^{v_s}_i \mid 1\leq i \leq n(v_s)\}$ be the set of all paths which end in $v_s$.  Then there is a $\mathbb Z$-graded isomorphism
\begin{equation}\label{skyone23}
\LL_R(E) \cong_{\gr} \bigoplus_{s=1}^t \M_{n(v_s)} (R)\big (|p^{v_s}_1|,\dots, |p^{v_s}_{n(v_s)}|\big).
\end{equation}

Furthermore, let $R$ be a division ring,  $F$ be another acyclic graph with sinks  $\{u_1,\dots,u_k\}$ and  $\{p^{u_s}_i \mid 1\leq i \leq n(u_s)\}$ be the set of all paths which end in $u_s$. Then $\LL_R(E) \cong_{\gr} \LL _R(F)$ if and only if $k=t$, and after a permutation of indices, $n(v_s)=n(u_s)$ and  $\{|p^{v_s}_i| \mid 1\leq i \leq n(v_s)\}$ and $\{|p^{u_s}_i| \mid 1\leq i \leq n(u_s)\}$ present the same list. 
\end{theorem}
\begin{proof}
The graded isomorphism ~(\ref{skyone23}) follows from Theorem~\ref{polyhead}, as the graph $E$ is a polycephaly graph with no cycle, and   $t$ petals $L_{n_1},\dots,L_{n_t}$ with  $n_s=0, 1\leq s\leq t$, i.e,  the vertices $\{v_1,\dots,v_t\}$ (see Definition~\ref{petaldef}). 
 
 For the second part, if $E$ and $F$ have the same number of sinks and the lengths of paths which end to sinks are the same as described in the statement of the Theorem, then by ~(\ref{skyone23}) (and~(\ref{kjhgf})), one obtains the same graded rings for $\LL_R(E)$ and $\LL_R(F)$, thus they are graded isomorphic. 
For the converse, assume that $\LL_R(E)$ is graded isomorphic to $\LL_R(F)$. We need to use~\cite[Lemma~3.8, p.37]{lamfc}, which states that if $A$ is a ring where $A=B_1\oplus \dots \oplus B_r= C_1 \oplus \dots \oplus C_s$, and $B_i$ and $C_j$ are indecomposable ideals then $r=s$ and after a permutation $B_i=C_i$.  
Writing $\LL_R(E)$ and $\LL_R(F)$ in the form of (\ref{skyone23}) and applying the lemma to the isomorphism $\LL_R(E) \cong_{\gr} \LL_R(F)$, we get $k=t$ and after a suitable permutation, for all $1\leq s \leq t$, \[\M_{n(v_s)} (R )\big (|p^{v_s}_1|,\dots, |p^{v_s}_{n(v_s)}|\big)\cong_{\gr} 
\M_{n(u_s)} (R)\big (|p^{u_s}_1|,\dots, |p^{u_s}_{n(u_s)}|\big).\] 
By Theorem~\ref{clasigr}, $n(v_s)=n(u_s)$ (with $\Ga=0$ as $R$ has a trivial grading), and  $|p^{v_s}_j|=|p^{u_s}_j|+m$, where $1\leq j \leq n(v_s)$ and $m\in \mathbb N$. However, since there are paths of lengths zero in both lists, this forces $m=0$ and thus $\{|p^{v_s}_i| \mid 1\leq i \leq n(v_s)\}$ and $\{|p^{u_s}_i| \mid 1\leq i \leq n(u_s)\}$ present the same list. This completes the proof. 
\end{proof}

\begin{example} \label{niroi}

The above theorem shows that the Leavitt path algebras of the graphs $E_1$ and $E_2$ below with coefficients from the field $K$ are graded isomorphic to 
$\M_{5} (K )\big (0,1,1,2,2)$ and thus $\LL(E_1) \cong_{\gr} \LL(E_2)$. However $\LL(E_3)\cong_{\gr} \M_5(K)(0,1,2,2,3)$, so the only if part of the Theorem implies that $\LL(E_3)$ is not graded isomorphic to the others.

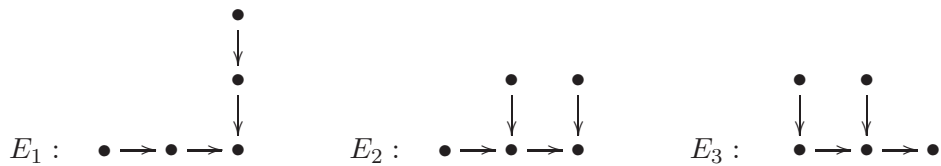
\[\xymatrix@=13pt{
& & \bullet \ar[d] \\
& & \bullet \ar[d] & &  & \bullet \ar[d] & \bullet \ar[d] &  & &   \bullet \ar[d] &  \bullet \ar[d]\\
E_1: \quad \bullet \ar[r] & \bullet \ar[r]  & \bullet  & & E_2: \quad  \bullet \ar[r] & \bullet \ar[r] &\bullet && 
E_3: & \bullet \ar[r] & \bullet \ar[r] &\bullet
}\]
\end{example}

\begin{deff}\label{mmjjsspp}
Let $E$ be an acyclic graph with a unique sink and $F$ be a row finite graph. Denote by $E\otimes F$ a graph obtained by attaching $E$ from the sink  to each vertex of $F$.
\end{deff}

\begin{example}
Consider the graphs

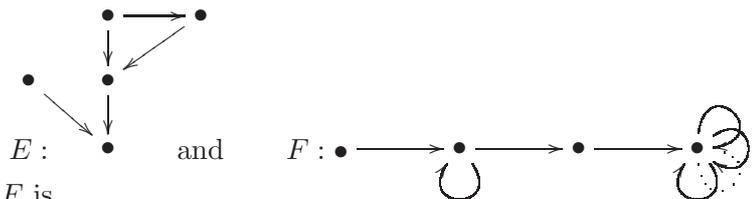
\begin{equation*}
\xymatrix@=13pt{
& \bullet \ar[d] \ar[r] & \bullet \ar[dl] & \\
\bullet \ar[dr] & \bullet\ar[d] & \\
E: & \bullet &
\text{ and } & F: \bullet \ar[rr] &&  \bullet \ar[rr] \ar@(dr,dl) && \bullet\ar[rr]  && \bullet \ar@(u,r) \ar@(ur,dr)  \ar@{.}@(r,d) \ar@(dr,dl)&.
}
\end{equation*}

Then the graph $E\otimes F$ is

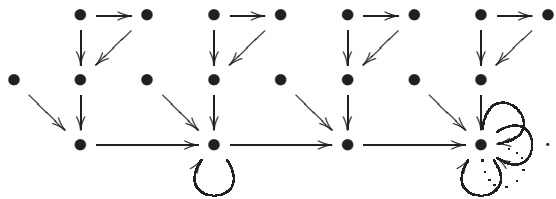
\begin{equation*}
\xymatrix@=13pt{
&\bullet \ar[d] \ar[r] & \bullet \ar[dl] &  \bullet \ar[d] \ar[r] & \bullet \ar[dl] &   \bullet \ar[d] \ar[r] & \bullet \ar[dl] &  \bullet \ar[d] \ar[r] & \bullet \ar[dl] &\\
\bullet \ar[dr] & \bullet\ar[d] & \bullet \ar[dr] & \bullet\ar[d] &   \bullet \ar[dr] & \bullet\ar[d] &   \bullet \ar[dr] & \bullet\ar[d] &\\
&   \bullet \ar[rr] & &   \bullet \ar[rr] \ar@(dr,dl) && \bullet\ar[rr]  && \bullet \ar@(u,r) \ar@(ur,dr)  \ar@{.}@(r,d) \ar@(dr,dl)&. 
}
\end{equation*}
\end{example}

\medskip 

\begin{lemma}\label{tensoir}
Let $K$ be a field, $E$ be an acyclic graph with a unique sink and $F$ be a row finite graph. Then \[\LL_K(E\otimes F) \cong_{\gr} \LL_K(E) \otimes_K \LL_K(F).\]
\end{lemma}
\begin{proof}
Let $u$ be the sink of the acyclic graph $E$ and let $\{p_i \mid 1\leq i \leq n \}$ be the set of all paths in $E$ which end in $u$. Then in $E\otimes F$ for each $v\in F^0$, by definition, we have one copy of the above set with the range $v$, which we denote them by
$\{p_i^v \mid 1\leq i \leq n \}$.
Recall from Theorem~\ref{acyclicc}, that the set of $p_ip_j^*, 1\leq i \not = j \leq n(v)$, is the basis of $\LL(E)$ and
\[\LL(E) \cong_{\gr} \M_{n(v)}(K)(|p_1|,\dots, |p_{n(v)}|).\] Define $\phi: \LL(E) \otimes_K \LL(F) \rightarrow \LL(E\otimes F)$ by $\phi(p_ip_j^* \otimes x  )=p_i^{s(x)} x {p_j^{r(x)}}^*$, where $x$ is a word in $\LL(F)$ and extend it linearly. Since the set of $p_ip_j^*$  is a basis of $\LL(E)$ and ${p_j^v}^*p_t^w=\delta_{jt} \delta_{vw} v$, it is an easy exercise to see that this map is well-defined, and is a $K$-algebra isomorphism. Now considering the grading defined on the tensor product of graded algebras (see \ref{tengr}), the map $\phi$ is clearly graded as well.
\end{proof}

\begin{example}
Let $K$ be a field. If $E$ and $F$ are acyclic graphs with unique sinks, then since tensor products on the level of algebras commute, we have
\[\LL_K(E\otimes F) \cong_{\gr} \LL_K(E)\otimes_K \LL_K(F) \cong_{\gr} \LL_K(F) \otimes_K \LL_K(E) \cong_{\gr} \LL_K(F\otimes E).\] For example, if $E$ is $\xymatrix@=10pt{\bullet \ar[r] & \bullet }$ and $F$ is $\xymatrix@=10pt{\bullet \ar[r] & \bullet \ar[r]& \bullet}$
then the Leavitt path algebras of the following two graphs are graded isomorphic.
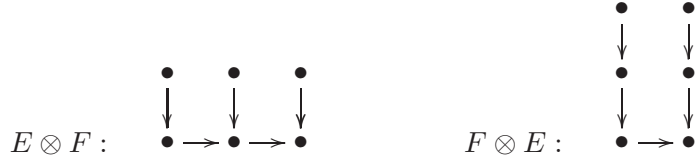
\begin{equation*}
\xymatrix@=13pt{
& &&&&&& \bullet \ar[d]&\bullet \ar[d]\\
& \bullet \ar[d] & \bullet\ar[d] &\ar[d]  \bullet \ar[d]&&& & \bullet \ar[d] & \bullet \ar[d]& \\
E \otimes F : & \bullet \ar[r] & \bullet \ar[r] & \bullet &&& F\otimes E: &  \bullet \ar[r] & \bullet&
}
\end{equation*}
In fact by  Theorem~\ref{acyclicc}, the Leavitt path algebras of both graphs are $\M_6(K)(0,1,1,2,2,3)$.
 \end{example}

\begin{theorem}\label{grcomet}
Let $R$ be a ring and  $E$ a multi-headed comets consisting of cycles $C_{l_s}$ for $1\leq s\leq t$ of length $l_s$. 
For any $1\leq s\leq t$ choose $v_s$ (an arbitrary vertex) in $C_{l_s}$ and remove the edge $c_s$ with $s(c_s)=v_s$ from the cycle $C_{l_s}$. In this new graph, let $\{p^{v_s}_i \mid 1\leq i \leq n(v_s)\}$ be the set of all paths which end in $v_s$.  Then 
there is a $\mathbb Z$-graded isomorphism 
\begin{equation}\label{skyone4}
\LL_R(E) \cong_{\gr} \bigoplus_{s=1}^t \M_{n(v_s)} (R[x^{l_s},x^{-l_s}] )\big (|p^{v_s}_1|,\dots, |p^{v_s}_{n(v_s)}|\big).
\end{equation}

Furthermore, let $R$ be a division ring and $F$ be another  multi-headed comets consisting of cycles $C_{l'_s}$ for $1\leq s\leq k$ of length $l'_s$. 
For any $1\leq s\leq k$ choose $u_s$ (an arbitrary vertex) in $C_{l'_s}$ and remove the edge $\alpha'_s$ with $s(\alpha'_s)=u_s$ from the cycle $C_{l'_s}$. In this new graph, let $\{p^{u_s}_i \mid 1\leq i \leq n(u_s)\}$ be the set of all paths which end in $u_s$. 
Then $\LL_R(E) \cong_{\gr}\LL_R(F)$ if and only if $k=t$, and after a permutation of indices, $l_s=l'_s$ and $n(v_s)=n(u_s)$ for any $1\leq s\leq t$ and 
$\{|p^{v_s}_i| \mid 1\leq i \leq n(v_s)\}$ can be obtained from  $\{|p^{u_s}_i| \mid 1\leq i \leq n(u_s)\}$ by applying {\upshape (\ref{pqow1})} and {\upshape (\ref{pqow2})}. 
\end{theorem}

\begin{proof}
The graded isomorphism ~(\ref{skyone4}) follows from Theorem~\ref{polyhead}, as the graph $E$ is a polycephaly graph with $t$ cycles, and   no roses  (see Definition~\ref{petaldef}). 

For the second part, if $E$ and $F$ have the same number of cycles of the same lengths  as described in the statement of the Theorem  and $\{|p^{v_s}_i| \mid 1\leq i \leq n(v_s)\}$ can be obtained from  $\{|p^{u_s}_i| \mid 1\leq i \leq n(u_s)\}$ by applying {\upshape (\ref{pqow1})} and {\upshape (\ref{pqow2})} then by ~(\ref{skyone4}), one obtains the same graded rings for $\LL_R(E)$ and $\LL_R(F)$, thus they are graded isomorphic. 

For the converse, assume that $\LL_R(E)$ is graded isomorphic to $\LL_R(F)$. Writing $\LL_R(E)$ and $\LL_R(F)$ in the form of (\ref{skyone4}), since each component is an indecomposable ideal, similar to the proof of Theorem~\ref{acyclicc}, using~\cite[Lemma~3.8, p.37]{lamfc}, it follows $k=t$ and after a suitable permutation of indices,  for all $1\leq s \leq t$, \[\M_{n(v_s)} (R[x^{l_s},x^{-l_s}] )\big (|p^{v_s}_1|,\dots, |p^{v_s}_{n(v_s)}|\big)\cong_{\gr} 
\M_{n(u_s)} (R[x^{l'_s},x^{-l'_s}])\big (|p^{u_s}_1|,\dots, |p^{u_s}_{n(u_s)}|\big).\]

Now Proposition~\ref{prohungi} forces $n(v_s)=n(u_s)$ and $ R[x^{l_s},x^{-l_s}] \cong_{\gr} R[x^{l'_s},x^{-l'_s}]$. Thus $l_s=l'_s$. Since $R[x^{l_s},x^{-l_s}]$ is a graded division ring, the rest  follows from Theorem~\ref{clasigr}. 
\end{proof}

\begin{example}\label{cycskew}
By Theorem~\ref{grcomet}, the Leavitt path algebra of a cycle $C_n$ of length $n\geq 1$ with coefficients in a division ring $R$ is 
$\mathbb Z$-graded isomorphic to $\M_n(R[x^n,x^{-n}])(0,1,\dots,n-1)$. We can describe this algebra in a more familiar way. 
Let $\{\alpha_i \mid 1\leq i \leq n \}$ be the set of all the edges of this graph, $s(\alpha_i)=v_i$ for all $i$, and $r(\alpha_i)=v_{i+1}$ for $i < n$, and $r(\alpha_n)=v_1$.  One can easily check that $x=\sum_{i=1}^n\alpha_i \in \LL_R(C_n)_1$ is invertible, namely $x \overline x =\overline x x=1$.  Since the grade group is the cyclic group $\mathbb Z$, this implies that all homogeneous components contain invertible elements. Therefore $\LL_R(C_n)$ is a crossed-product (see~\S\ref{gradingsec}). Furthermore, from the graph of $C_n$, it is easy to see that the only monomials of degree zero are the vertices, which we know are $R$-linearly independent. Therefore the ring of homogeneous elements of degree zero, $\LL_R(C_n)_0$ is $\bigoplus_{i=1}^n Rv_i$, i.e, it is $n$-copies of $R$. Furthermore $x v_i\overline x=(\alpha_1+\cdots+\alpha_n)(v_i)(\alpha_1^*+\cdots+\alpha_n^*)=v_{i-1}$, unless $i=1$ which in this case $xv_1\overline x=v_n.$ This together with the fact that $\LL_R(C_n)$ is cross-product, shows that $\LL_R(C_n)$ is graded isomorphic to the skew-group ring $R^n \star_\sigma  \mathbb Z$, where $\sigma(k_1,\dots,k_n)=(k_2,\dots,k_n,k_1)$. Observe that if $n=1$, then $\sigma=\id$ and the LPA reduces to the group ring $R[\mathbb Z]=R[x,x^{-1}]$. (See also Theorem~\ref{polyheadnew}.)
\end{example}

\begin{example}\label{noncori}
By Theorem~\ref{sthfin}, the Leavitt path algebra of the graph $E$ with coefficients in a field $K$,
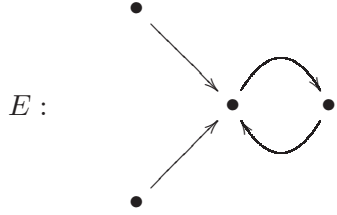
\begin{equation*}
\xymatrix{
& \bullet \ar[dr] &\\
E:  &  & \bullet \ar@/^1.5pc/[r] & \bullet \ar@/^1.5pc/[l] &\\
& \bullet \ar[ur] &&\\}
\end{equation*}
is strongly graded. By Theorem~\ref{grcomet}, 
\begin{equation}\label{wokj}
\LL_K(E) \cong_{\gr}\M_4(K[x^2,x^{-2}])(0,1,1,1).
\end{equation} However this algebra is not  a crossed-product. Set $B=K[x,x^{-1}]$ with the grading $B=\textstyle{\bigoplus_{n\in \mathbb Z}} Kx^n$ and consider $A=K[x^2,x^{-2}]$ as a graded subring of $B$ with 
$A_n=Kx^n$ if $n\equiv 0 \pmod{2}$, and $A_n=0$ otherwise.  
Using the graded isomorphism of ~(\ref{wokj}), by (\ref{mmkkhh})  a homogeneous element of degree $1$ in $\LL_K(E)$ has the form 
\begin{equation*}
\begin{pmatrix}
A_{1} & A_{2} & A_2 & A_2\\
A_{0} & A_{1} & A_1 & A_1\\
A_{0} & A_{1} & A_1 & A_1\\
A_{0} & A_{1} & A_1 & A_1
\end{pmatrix}.
\end{equation*}
Since $A_1=0$, the determinants of these matrices are zero, and thus no homogeneous   element of degree $1$ is  invertible. Thus $\LL_K(E)$ is not crossed-product (see~\S\ref{gradingsec}). 

Now consider the following graph from Introduction \S\ref{introf},

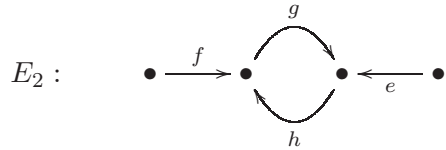
\begin{equation*}
\xymatrix{
E_2:  & \bullet \ar[r]^f & \bullet \ar@/^1.5pc/[r]^g & \bullet \ar@/^1.5pc/[l]^h & \bullet \ar[l]^e& \\
}
\end{equation*}
By Theorem~\ref{grcomet}, 
\begin{equation}\label{wokj5}
\LL_K(E) \cong_{\gr}\M_4(K[x^2,x^{-2}])(0,1,1,2).
\end{equation}
Using the graded isomorphism of ~(\ref{wokj5}), by (\ref{mmkkhh})   homogeneous elements of degree $0$ in $\LL_K(E)$ have the form 
\begin{equation*}\LL_K(E)_0=
\begin{pmatrix}
A_{0\phantom{-}} & A_{1\phantom{-}} & A_{1\phantom{-}} & A_{2\phantom{-}}\\
A_{-1} & A_{0\phantom{-}} & A_{0\phantom{-}} & A_{1\phantom{-}}\\
A_{-1} & A_{0\phantom{-}} & A_{0\phantom{-}} & A_{1\phantom{-}}\\
A_{-2} & A_{-1} & A_{-1} & A_{0\phantom{-}}
\end{pmatrix}=
\begin{pmatrix}
K\phantom{x^{-2}} & 0 & 0 & Kx^2\\
0\phantom{x^{-2}} & K & K & 0\phantom{x^{2}}\\
0\phantom{x^{-2}} & K & K & 0\phantom{x^{2}}\\
Kx^{-2} & 0 & 0 & K\phantom{x^{2}}
\end{pmatrix}.
\end{equation*}
In the same manner, homogeneous elements of degree one have the form,
\begin{equation*}\LL_K(E)_1=
\begin{pmatrix}
0 & Kx^2 & Kx^2 & 0\phantom{x^{2}}\\
K & 0\phantom{x^{2}} & 0\phantom{x^{2}} & Kx^2\\
K & 0\phantom{x^{2}} & 0\phantom{x^{2}} & Kx^2\\
0 & K\phantom{x^{2}} & K\phantom{x^{2}} & 0\phantom{x^{2}}
 \end{pmatrix}.
\end{equation*}
Choose 
\begin{equation*}u=
\begin{pmatrix}
0 & 0 & x^2 & 0\phantom{^{2}}\\
0 & 0 & 0\phantom{^{2}} & x^2\\
1 & 0 & 0\phantom{^{2}} & x^2\\
0 & 1 & 0\phantom{^{2}} & 0\phantom{^{2}}
 \end{pmatrix} \in \LL(E)_1
\end{equation*}
and observe that $u$ is invertible (this matrix corresponds to the element $g+h+fge^*+ehf^* \in \LL_K(E)_1$). 

Thus $\LL_K(E)$ is crossed-product (and therefore a skew group ring as the grading is cyclic), i.e., \[\LL_K(E)\cong_{\gr} \bigoplus_{i\in \mathbb Z}\LL_K(E)_0 u^i\] and a simple calculation shows that one can describe this algebra as follows:
\[\LL_K(E_2)_0 \cong \M_2(K)\times \M_2(K),\] and 
\begin{equation}\label{coiu}
\LL_K(E_2) \cong_{\gr} \big ( \M_2(K)\times \M_2(K) \big ) \star_\tau \mathbb Z,
\end{equation}
where, 
 \begin{equation*}
\tau(\left( 
\begin{matrix} 
a_{11} & a_{12}\\ 
a_{21} & a_{22}  
\end{matrix} 
\right) ,
\left( 
\begin{matrix} 
b_{11} & b_{12}\\ 
b_{21} & b_{22}  
\end{matrix} 
\right))=(
\left( 
\begin{matrix} 
b_{22} & b_{21}\\ 
b_{12} & b_{11}  
\end{matrix} 
\right) , 
\left( 
\begin{matrix} 
a_{22} & a_{21}\\ 
a_{12} & a_{11}  
\end{matrix} 
\right) ).
\end{equation*}
\end{example}

\begin{corollary}
Let $R$ be a division ring and  $E$  be a $C_n$-comet with the cycle $C$ of length $n \geq 1$. Let $\{p_i \mid 1\leq i \leq m\}$ be the paths as in Theorem~{\upshape \ref{grcomet}}. Let $\varepsilon_1+n\mathbb Z,\dots,\varepsilon_k+n\mathbb Z$ be the distinct elements of all elements $|p_i|+n\mathbb Z$, $1\leq i \leq m$.  For each $\varepsilon_l$ let $r_l$ be the number of $i$ such that $|p_i|  \equiv \varepsilon_l \pmod{n\mathbb Z}$.
Then  
\begin{equation*}
\LL_R(E)_0 \cong \M_{r_1}(R)\times \dots \times \M_{r_k}(R).
\end{equation*}
Furthermore, $\LL_R(E)_0$ is simple if and only if the head of the comet is a loop, i.e., $|C|=1$. 
\end{corollary}
\begin{proof}
This follows from Theorem~\ref{grcomet} and (\ref{urnha}). 
\end{proof}

\begin{example}
This example shows that even by taking the $\mathbb Z$-grading into consideration, the Leavitt path algebra functor does not distinguish between the graphs below and would give a $\mathbb Z$-graded isomorphism between their associated Leavitt path algebras. Let
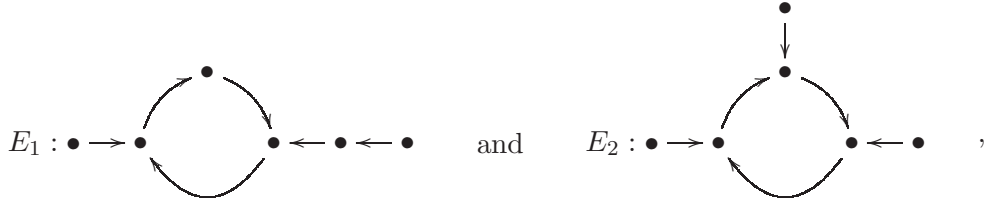
\begin{equation*}
\xymatrix@=13pt{
& & & & & & &&&    \bullet \ar[d] &&&  \\
& & \bullet \ar@/^/[dr] &&&&&&&  \bullet \ar@/^/[dr] \\
E_1 : \bullet \ar[r] & \bullet \ar@/^/[ur] &&\bullet \ar@/^1.7pc/[ll] & \bullet \ar[l] & \bullet \ar[l] & \text{ and }& 
E_2 : \bullet \ar[r] & \bullet \ar@/^/[ur] &&\bullet \ar@/^1.7pc/[ll] & \bullet \ar[l] &, 
}
\end{equation*}

%\bigskip
\bigskip 
 
then by Theorem~\ref{grcomet} we have 
\[\LL_K(E_1) \cong_{\gr} \LL_K(E_2) \cong_{\gr} \M_6(K[x^3,x^{-3}])(0,1,1,2,2,3).\]
\end{example}

\section{Weighted Leavitt Path Algebras}\label{secitwo}

In this section we introduce the {\it weighted} Leavitt path algebras (wLPA for short) starting from a weighted graph (i.e., a graph which each edge comes with some copies of itself). This is a graded algebra, graded by a free abelian group (i.e., some copies of $\mathbb Z$) which in the special case of a graph with  weights $1$ (or unweighted), it gives the Leavitt path algebras and in its simplest form, when the graph $E$ has only one vertex and $n+k$ loops of weights $n$,  $\LL(E)$ recovers the algebra constructed by Leavitt~\cite[p.30]{vitt62} and ~\cite[p.322]{vitt57} which is of type $(n,k)$ (see Example~\ref{lat}).  (In fact we can define weighted Leavitt path algebras with $G$-grading where $G$ is any group, see Remark~\ref{ggrad}.)  In this note, after defining wLPA, we establish their basic properties.  One could then ask for characterization of an wLPA in terms of its underlying graph and its module and algebraic structure. These shall be dealt with in future papers.

We begin by the definition of a weighted graph.

\begin{deff}\label{weighted}
A {\it weighted graph} $E=(E^0,E^{\st}, E^1,r,s,w)$ consists of three countable sets, $E^0$ called {\it vertices}, $E^{\st}$ {\it structured edges} and $E^1$ {\it edges} and maps $s,r:E^{\st}\rightarrow E^0$,  and  a {\it weight map} $w:E^{\st} \rightarrow \mathbb N$ such that $E^1= \coprod_{\alpha \in E^{\st}} \{\alpha_i\mid 1\leq i \leq w(\alpha)\}$, i.e., for any $\alpha\in E^{\st}$, with $w(\alpha)=k$, there are $k$ distinct elements $\{\alpha_1,\dots,\alpha_k\}$, and $E^1$ is the disjoint union of  all such sets for all $\alpha \in E^{\st}$. 
\end{deff}

%Therefore if  $\alpha \in E^s$ and $w(\alpha)=k$, then there are $k$ edges $\alpha_1,\dots \alpha_k\in E^1$ with the same source and range as $\alpha$, i.e., $s(\alpha_i)=s(\alpha)$ and $r(\alpha_i)=r(\alpha)$ for $1\leq i \leq k$.

%We refer to the maps $r_{E^s}$ and $s_{E^s}$ if we need to distinguish them from those on $E^1$.  

We sometimes write $(E,w)$ to emphasis the graph is  weighted. If $s^{-1}(v)$ is a finite set for every $v \in E^0$, then the graph is called {\it row-finite}.  In this note we will consider only the row-finite graphs. In this setting, if the number of vertices, i.e.,  $|E^0|$,  is finite, then the number of edges, i.e.,  $|E^1|$, is finite as well and we call $E$ a {\it finite} graph. 

\begin{deff}{\sc Weighted Leavitt path algebras.} \label{iLPA} \\
For a weighted graph $E$ and a ring $R$ with identity, we define the {\it weighted Leavitt path algebra of $E$}, denoted by $\LL_R(E,w)$, to be the algebra generated by the sets $\{v \mid v \in E^0\}$, $\{ \alpha_1,\dots \alpha_{w(\alpha)} \mid \alpha \in E^{\st} \}$ and $\{ \alpha_1^*,\dots \alpha_{w(\alpha)}^* \mid \alpha \in E^{\st} \}$ with the coefficients in $R$, subject to the relations 

\begin{enumerate}
\item $v_iv_j=\delta_{ij}v_i \textrm{ for every } v_i,v_j \in E^0$.

\item $s(\alpha)\alpha_i=\alpha_i r(\alpha)=\alpha_i \textrm{ and }
r(\alpha)\alpha_i^*=\alpha_i^*s(\alpha)=\alpha_i^*  \textrm{ for all } \alpha \in E^{\st}$ and $1\leq i \leq w(\alpha)$.

\item $\sum_{\{\alpha \in E^{\st} \mid s( \alpha)=v\}} \alpha_i \alpha_j^*=\delta_{ij}s(\alpha)$ for fixed 
$1 \leq i,j \leq \max\{w(\alpha)\mid \alpha\in E^{\st}, s( \alpha)=v\}$. 

\item $\sum_{1\leq i \leq \max\{ w(\alpha),w(\alpha')\}} \alpha_i^* \alpha'_i=\delta_{\alpha \alpha'}r(\alpha)$, for all $\alpha, \alpha' \in E^{\st}$.
\end{enumerate}
\end{deff}
Here the ring $R$ commutes with the generators $\{v,\alpha, \alpha^* \mid v \in E^0,\alpha \in E^1\}$. Also in relations (3) and (4), we set $\alpha_i$ and $\alpha_i^*$ zero whenever $i>w(\alpha)$. When the coefficient ring $R$ is clear from the context, we simply write $\LL(E,w)$ instead of $\LL_R(E,w)$. When $R$ is not commutative, then we consider $\LL_R(E,w)$ as a left $R$-module.

\begin{example} \label{calcu} We compare the relations of the weighted Leavitt path algebra $\LL(E,w)$ and the usual Leavitt path algebras $\LL(E')$ and $\LL(E'')$ in the following:

\begin{equation*}
\xymatrix{
(E,w):&u \ar@/^/[r]^{\alpha_1,\alpha_2}  \ar@/_/[r]_{\beta_1,\beta_2} & v & &E':&u \ar@/^/[r]|{\alpha_1} \ar@/^1pc/[r]|{\alpha_2} \ar@/_/[r]|{\beta_1} \ar@/_1pc/[r]|{\beta_2} & v &&
E'':&u \ar@/^/[r]^{\alpha_1}  \ar@/_/[r]_{\beta_1}  & v
}
\end{equation*}

\begin{align}
\alpha_1\alpha_1^*+\beta_1\beta_1^*=u& & \alpha_1\alpha_1^*+\alpha_2\alpha_2^*+\beta_1\beta_1^*+\beta_2\beta_2^*=u&& \alpha_1\alpha_1^*+\beta_1\beta_1^*=u \notag\\
\alpha_2\alpha_2^*+\beta_2\beta_2^*=u&&\alpha_1^*\alpha_1=\alpha_2^*\alpha_2=\beta_1^*\beta_1=\beta_2^*\beta_2=v&& \alpha_1^*\alpha_1=\beta_1^*\beta_1=v\notag\\
\alpha_1\alpha_2^*+\beta_1\beta_2^*=0&& \alpha_i^*\alpha_j=\beta_i^*\beta_j=0 \textrm{ if } i \not = j&& \alpha_1^*\beta_1=\beta_1^*\alpha_1=0
\notag\\
\alpha_2\alpha_1^*+\beta_2\beta_1^*=0&& \alpha_i^*\beta_j=\beta_j^*\alpha_i=0 \textrm{ for all } i,j &&\notag\\
\alpha_1^*\alpha_1+\alpha_2^*\alpha_2=v&\notag\\
\beta_1^*\beta_1+\beta_2^*\beta_2=v&\notag\\
\alpha_1^*\beta_1+\alpha_2^*\beta_2=0&\notag\\
\beta_1^*\alpha_1+\beta_2^*\alpha_2=0&
\notag
\end{align}
The graphs $E'$ and $E''$ are acyclic graphs whose Leavitt path algebras were characterized in Theorem~\ref{acyclicc}. Note that in $\LL(E,w)$, relations (3) and (4) in Definition~\ref{iLPA} amounts to 
\[
\left(
\begin{matrix}
\alpha_1& \beta_1\\
\alpha_2& \beta_2\\
\end{matrix}\right) 
\left(
\begin{matrix}
\alpha_1^*& \alpha_2^*\\
\beta_1^*& \beta_2^*\\
\end{matrix}\right) =
\left(
\begin{matrix}
u& 0\\
0& u\\
\end{matrix}\right)
\text{ \quad and \quad       } 
\left(\begin{matrix}
\alpha_1^*& \alpha_2^*\\
\beta_1^*& \beta_2^*\\
\end{matrix}\right)
\left(
\begin{matrix}
\alpha_1& \beta_1\\
\alpha_2& \beta_2\\
\end{matrix}\right) 
=
\left(
\begin{matrix}
v& 0\\
0& v\\
\end{matrix}\right).
\]
\end{example}

\begin{example} \label{lpa} {\sc Leavitt Path Algebras.} \\
Let the weight map $w:E^{\st}\rightarrow \mathbb N$ be the constant function $w(\alpha)=1$ for any $\alpha \in E^{\st}$. Then $\LL_R(E,w)$ is the usual Leavitt path algebra (with the coefficients in the ring $R$) as defined in~\cite{aap05} and~\cite{amp}.
\end{example}

\begin{example} \label{lat} {\sc The Leavitt algebra of type $(n,k)$.}\\
Let $R$ be a division ring. For positive integers $n$ and $k$, let the structured edges $E^{\st}$ of a graph $E$ consist of $n+k$ loops, i.e., $s(y)=r(y)$ for  $y \in E^{\st}$ and let the weight function be the constant map $w(y)=n$ for all $y\in E^{\st}$. We visualize this data as follows:
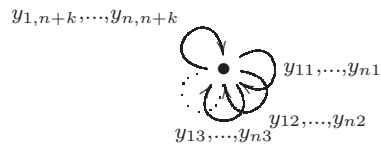
\begin{equation*}
\xymatrix{
& \bullet \ar@{.}@(l,d) \ar@(ur,dr)^{y_{11},\dots,y_{n1}} \ar@(r,d)^{y_{12},\dots,y_{n2}} \ar@(dr,dl)^{y_{13},\dots,y_{n3}} \ar@(l,u)^{y_{1,n+k},\dots,y_{n,n+k}}& 
}
\end{equation*}
Then the weighted Leavitt path algebra associated to $E$, $\LL_R(E,w)$, is the algebra constructed by Leavitt in~\cite[p.190]{vitt56}, for $n=2$ and $k=1$, where he showed that this algebra has no zero divisors, in
~\cite[p.322]{vitt57}, for arbitrary $n$ and $k=1$ and
in~\cite[p.130]{vitt62}  for arbitrary $n$ and $k$ and established that  these algebras are domains and of type $(n,k)$. Recall that a ring $A$ is of type $(n,k)$ if $n$ and $k$ are the least positive integers such that $A^n\cong A^{n+k}$ as $A$-modules. To recover Leavitt's algebra from Definition~\ref{iLPA} (and to arrive to his notations), let $E^{\st}=\{y_1,\dots,y_{n+k}\}$ be the structured edges and denote $(y_s)_r=y_{rs}\in E^1$, for $1\leq r\leq n$ and $1\leq s\leq n+k$. Denote $y_{rs}^*=x_{sr}$ and arrange the $y$'s and $x$'s in the matrices
\begin{equation} \label{breaktr}
Y=\left( 
\begin{matrix} 
y_{11} & y_{12} & \dots & y_{1,n+k}\\ 
y_{21} & y_{22} & \dots & y_{2,n+k}\\ 
\vdots & \vdots & \ddots & \vdots\\ 
y_{n1} & y_{n2} & \dots & y_{n,n+k} 
\end{matrix} 
\right), \qquad 
X=\left( 
\begin{matrix} 
x_{11\phantom{n+{},}} & x_{12\phantom{n+{},}} & \dots & x_{1n\phantom{n+{},}}\\ 
x_{21\phantom{n+{},}} & x_{22\phantom{n+{},}} & \dots & x_{2n\phantom{n+{},}}\\ 
\vdots & \vdots & \ddots & \vdots\\ 
x_{n+k,1} & x_{n+k,2} & \dots & x_{n+k,n} 
\end{matrix} 
\right) 
\end{equation} 
Then Condition~(3) of Definition~\ref{iLPA} precisely says that $Y\cdot X=I_{n,n}$ and Condition~(4) is equivalent to 
$X\cdot Y=I_{n+k,n+k}$ which is how Leavitt defines his algebra. We denote this algebra by $\LL_R(n,k+1)$. (See Definition~\ref{wpetaldef}; Cohn's notation in~\cite{cohn11} for this algebra is $V_{n,n+k}$.)
\end{example}

%\begin{remark}
%In the similar manner, one can define weighted graph algebras and weighted graph $C^*$-algebras. 
%\end{remark} 

\begin{remark} {\sc $G$-weighted Leavitt Path Algebras.} \label{ggrad}\\
The idea of structured edges in Definition~\ref{weighted} is to bundle together certain edges in order to be able to define relations~(3) and~(4) in Definition~\ref{iLPA}. Instead of introducing the structured edges, one can partition the (in-coming and out-going) edges and define the $G$-graded Leavitt path algebras, for any arbitrary group $G$, such that when $G=\mathbb Z$ or $G=\textstyle{\bigoplus_n}\mathbb Z$, we obtain the LPAs or wLPAs, respectively. The construction is as follows:

Let $E$ be a directed graph and $G$ be an arbitrary group with the identity element $e$. Let $w:E^1\rightarrow G$ be a {\it weight} map and further define $w(\alpha^*)=w(\alpha)^{-1}$, for any edge $\alpha \in E^1$ and $w(v)=e$ for $v\in E^0$.  Write $\alpha_g$ for an edge $\alpha$ with the weight $g$ (i.e., $w(\alpha)=g$).  For any $v\in E^0$, consider partitions of the sets $s^{-1}(v)$ and $r^{-1}(v)$ (if they are not empty) into $\bigcup_{i\in I_v} \mathcal P_i(v)$ and $\bigcup_{j\in J_v} \mathcal Q_j(v)$  for some index sets $I_v$ and $J_v$, respectively.  Since $E$ is a row-finite graph, each set $\mathcal P_i(v)$, $i\in I_v$, is finite. We further assume that each set $\mathcal Q_j(v)$, $j \in J_v$, is also finite.

Let $\LL_R(E,G)$ be the algebra generated by the sets $\{v \mid v \in E^0\}$ and $\{ \alpha, \alpha^* \mid \alpha  \in E^1 \}$ with the coefficients in $R$, subject to the relations

\begin{enumerate}
\item $v_iv_j=\delta_{ij}v_i \textrm{ for every } v_i,v_j \in E^0$.

\item $s(\alpha)\alpha=\alpha r(\alpha)=\alpha \textrm{ and }
r(\alpha)\alpha^*=\alpha^*s(\alpha)=\alpha^*  \textrm{ for all } \alpha \in E^1$.

\item $\sum_{\{ i \in I_v, \alpha_g, \alpha_h' \in \mathcal P_i(v) \} } \alpha_g {\alpha'}_{h^{-1}}^*=\delta_{g h} \delta_{\alpha \alpha'} v$ for $v \in E^0$ and any fixed  $g, h \in G$. 

\item $\sum_{g\in G} \alpha_{g^{-1}}^* \alpha'_g=\delta_{\alpha \alpha'}r(\alpha_g)$, where $\alpha_g \in \mathcal Q_i(v), \alpha'_g \in \mathcal Q_j(w)$ for all $v,w \in E^0$ and $i\in I_v$, $j \in J_w$.

\end{enumerate}

In (3), $\delta_{\alpha \alpha'}=1$ if $\alpha_g=\alpha'_g$, where $\alpha_g,\alpha'_g \in \mathcal P_i(v)$, for all $i \in I_v$ and $0$ otherwise. In (4) $\delta_{\alpha \alpha'}=1$ if $\alpha_g=\alpha'_g$ for all $g\in G$ and $0$ otherwise. Also the ring $R$ commutes with the generators $\{v,\alpha, \alpha^* \mid v \in E^0,\alpha \in E^1\}$. Also in relations (3) and (4), % we set $\alpha_g=0$ if $w(\alpha)\not = g$. Also 
we set $\alpha_g=0$ if $\alpha_g \not \in \mathcal P_i(v)$ or $\mathcal Q_i(v)$. 

For a path $p=\alpha_1\alpha_2\dots\alpha_k$, set $w(p)=w(\alpha_1)w(\alpha_2)\dots w(\alpha_k)$ 
 and extend it to the monomials of the algebra in an obvious manner. This defines a $G$-graded algebra which we call {\it $G$-weighted Leavitt path algebra}. 
 
For example, consider the following $G$-weighted graph: 
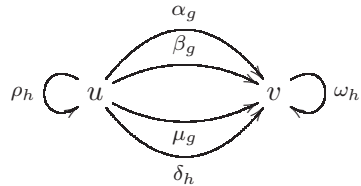
\begin{equation*}
\xymatrix{
  u  \ar@(lu,ld)_{\rho_h} \ar@/^0.9pc/[rr]^{\beta_g} \ar@/^2pc/[rr]^{\alpha_g} 
\ar@/_0.9pc/[rr]_{\mu_g} \ar@/_2pc/[rr]_{\delta_h}    && v \ar@(ru,rd)^{\omega_h} 
 }
\end{equation*}
We write some of the relations coming from Relations (3) and (4) of the definition above. Let $\mathcal P_1(u)=\{\alpha_g,\beta_g,\rho_h \}$ and $\mathcal P_2(u)=\{\mu_g,\delta_h\}$. Then for $g$ fixed, Relation (3) gives,  
\begin{align*}
\alpha_g \alpha^*_{g^{-1}}+\mu_g \mu^*_{g^{-1}}=u\\
\beta_g \beta^*_{g^{-1}}+\mu_g \mu^*_{g^{-1}}=u\\
\rho_h \rho^*_{h^{-1}}+\delta_h \delta^*_{h^{-1}}=u\\
\alpha_g \beta^*_{g^{-1}}+\mu_g \mu^*_{g^{-1}}=0\\
\beta_g \alpha^*_{g^{-1}}+\mu_g \mu^*_{g^{-1}}=0.\\
\end{align*}
Furthermore for $g$ and $h$ fixed, we obtain:
\begin{align*}
\alpha_g \rho^*_{h^{-1}}+\mu_g \delta^*_{h^{-1}}=0\\
\rho_h \alpha^*_{g^{-1}}+\delta_h \mu^*_{g^{-1}}=0\\
\beta_g \rho^*_{h^{-1}}+\mu_g \delta^*_{h^{-1}}=0\\
\rho_h \beta^*_{g^{-1}}+\delta_h \mu^*_{g^{-1}}=0.\\
\end{align*}
Let $\mathcal Q_1(v)=\{\alpha_g,\beta_g,\delta_h \}$ and  $\mathcal Q_2(v)=\{\mu_g,\gamma_h\}$.  For $\mathcal Q_2(v)$ fixed, Relation (4) gives, 
\begin{align*}
\mu^*_{g^{-1}} \mu_g+\gamma^*_{h^{-1}} \gamma_h=v
\end{align*}
For $\mathcal Q_1(v)$ fixed we obtain:
\begin{align*}
\alpha^*_{g^{-1}} \alpha_g+\delta^*_{h^{-1}} \delta_h=v\\
\beta^*_{g^{-1}} \beta_g+\delta^*_{h^{-1}} \delta_h=v\\
\alpha^*_{g^{-1}} \beta_g+\delta^*_{h^{-1}} \delta_h=0\\
\beta^*_{g^{-1}} \alpha_g+\delta^*_{h^{-1}} \delta_h=0.
\end{align*}
Also for $\mathcal Q_1(v)$ and $\mathcal Q_2(v)$ fixed we have:
\begin{align*}
\alpha^*_{g^{-1}} \mu_g+\delta^*_{h^{-1}} \gamma_h=0\\
\mu^*_{g^{-1}} \alpha_g+\gamma^*_{h^{-1}} \delta_h=0\\
\beta^*_{g^{-1}} \mu_g+\delta^*_{h^{-1}} \gamma_h=0\\
\mu^*_{g^{-1}} \beta_g+\gamma^*_{h^{-1}} \delta_h=0.
\end{align*}

One can write further relations by considering $\mathcal Q_1(u)$, $\mathcal Q_1(v)$ and $\mathcal Q_2(v)$ as well.

 It is clear that if $G=\mathbb Z$ and the weight map is the constant map assigning $1$ to each edge and moreover each partition set has only one element, then the above relations  give back the usual Leavitt path algebra.  Also, for a weighted graph with structured edges, 
 considering  \[\mathcal P_i(v)=\big\{ \alpha_{i1},\dots,\alpha_{ik} \mid \alpha_i \in E^{\st}, s(\alpha_i)=v, w(\alpha_i)=k\big\},\] and similarly for $\mathcal Q$'s, and moreover 
 setting
 $G=\textstyle{\bigoplus_n} \mathbb Z$, where $n=\max\{w(\alpha)\mid \alpha \in E^{\st}\}$, and assigning $(0,\dots,0,1,0,\dots)$, where $1$ is in the $k$-th component, for the weight of $\alpha_{ik}$,  $1\leq k \leq w(\alpha_i)$, we retrieve the weighted Leavitt path algebra in Definition~\ref{iLPA}.

Finally note that (as it is demonstrated in the above example) if in some partition set $\mathcal P_i(u)$ (or $\mathcal Q_i(u)$), there is more than one element with the same weight, i.e., there are $\alpha_g \not = \beta_g\in \mathcal P_i(u)$, then each of the relations (3) and (4) will give more than one relation. This clearly does not happen in the setting of wLPA of Definition~\ref{iLPA}. 

We refer the reader to \cite{green} for a study of  $G$-graded algebras arising from graphs with relations, in particular \cite[Theorem~3.4]{green}.

 \end{remark}

We now turn to weighted Leavitt path algebras and establish some of their basic properties. 

\begin{proposition} Let $E$ be a weighted graph and $\LL_R(E)$ be a weighted Leavitt path algebra with coefficients in a ring $R$. Then we have  

\begin{enumerate}[\upshape(1)]

\item $\LL_R(E)$ is a $\textstyle{\bigoplus_n} \mathbb Z$-graded ring with an involution where $n=\max\{w(\alpha)\mid \alpha \in E^{\st}\}$. 

\item $\LL_R(E)$ is a ring with local identities. If $E$ is finite, then $\LL_R(E)$ is a ring with identity. 
\end{enumerate}
\end{proposition}

\begin{proof}
(1) For the free ring generated by  $\{v \mid v \in E^0\}$, $\{ \alpha_1,\dots \alpha_{w(\alpha)} \mid \alpha \in E^{\st} \}$ and $\{ \alpha_1^*,\dots \alpha_{w(\alpha)}^* \mid \alpha \in E^{\st} \}$, with the coefficients in $R$, set for $v \in E^0$, $\deg(v)=0$, for $\alpha \in E^{\st}$, $1\leq i\leq w(\alpha)$, $\deg(\alpha_i)=(0,\dots,0,1,0,\dots)$ and $\deg(\alpha_i^*)=(0,\dots,0,-1,0,\dots) \in  \textstyle{\bigoplus_n} \mathbb Z$, where $n=\max\{w(\alpha)\mid \alpha \in E^{\st}\}$ and 
$1$ and $-1$ are in the $i$-th component, respectively. This defines a $\textstyle{\bigoplus_n} \mathbb Z$-grading on this free ring where $n=\max\{w(\alpha)\mid \alpha \in E^{\st}\}$ ($n$ could be infinite).  Note that all the relations in Definition~\ref{iLPA} involve  homogeneous elements, so the quotient of this algebra by the homogeneous ideal generated by these relations, i.e., $\LL_R(E,w)$ is also a graded ring. 

 To show that $\LL(E,w)$ is equipped with an involution, define a homomorphism from the free ring generated by $\{v \mid v \in E^0\}$, $\{ \alpha_1,\dots \alpha_{w(\alpha)} \mid \alpha \in E^{\st} \}$ and $\{ \alpha_1^*,\dots \alpha_{w(\alpha)}^* \mid \alpha \in E^{\st} \}$ with the coefficients in $R$, to $\LL(E,w)^{\op}$, the opposite ring of $\LL(E,w)$,  by sending $r\mapsto r$,  $v\mapsto v$, $\alpha_i\mapsto \alpha_i^*$ and $\alpha_i^*\mapsto \alpha_i$, where $r\in R$, $v\in E^0$ and $\alpha_i \in E^1$. One can see that all the relations in Definition~\ref{iLPA} are in the kernel of this ring homomorphism, thus inducing a homomorphism of order two from $\LL(E,w)$ to $\LL(E,w)^{\op}$. 

(2) Recall that a ring $A$ has local identities if for any finite subset $S \subseteq A$, there is an idempotent $e\in A$ such that $S \subseteq eAe$. The set of all such idempotents is called a set of local identities for $A$. 

Note that if $a_i$'s are mutually orthogonal idempotents in a ring $A$ such that $A=\sum a_iA=\sum Aa_i$, then the set of $\sum_{\textrm{finite}} a_i$ is a set of local identities of this ring. If the number of $a_i$'s is finite then $\sum a_i$ is an identity for this ring. Now it is easy to see that the set of vertices in $E$ is such a system of idempotents for $\LL_R(E,w)$. 
\end{proof}

 When $R$ is a division ring, by constructing a representation of $\LL_R(E)$, one can show that the vertices of a graph $E$ are linearly independent in $\LL_R(E)$ and the edges and ghost edges are not zero (see Lemma~1.5 in \cite{goodearl}). In the next theorem we will carry over this in the generalized setting of weighted Leavitt path algebras  and therefore covering the special case of LPA as a corollary (by setting the weight map the constant map $1$). 

\begin{theorem}\label{liniin}
Let $R$ be a division ring and $E$ be a weighted graph. Then the vertices of $E$ are $R$-linearly independent in the weighted Leavitt path algebra $\LL_R(E,w)$. 
\end{theorem}
\begin{proof}
Let $X$ be a (left) vector space over $R$ with an infinite countable basis. In the ring $\End_R(X)$, we will find nonzero elements $\{\Pp_v \mid v\in E^0\}$, $\{ \A_1,\dots \A_{w(\alpha)} \mid \alpha \in E^{\st} \}$ and $\{ \A_1^*,\dots \A_{w(\alpha)}^* \mid \alpha \in E^{\st} \}$ which satisfy the relations in Definition~\ref{iLPA} of a weighted Leavitt path algebra. Since the basis of $X$ is an infinite countable set, we can decompose $X=\bigoplus_{v\in E^0}X_v$, where each $X_v$ is a subspace of $X$ with an infinite countable basis.  Moreover, for each $v\in E^0$, we further  decompose $X_v$,
 \begin{equation}\label{decio}
 X_v=\bigoplus_{\{\alpha\in E^{\st},s(\alpha)=v \}} Y_{\alpha_i}, \text{ for fixed } 1\leq i \leq \max \big \{w(\alpha) \mid \alpha \in E^{\st}, s(\alpha)=v \big \},
 \end{equation}
 \begin{equation}\label{decio2}
 X_v=\bigoplus_{ 1\leq i \leq w(\alpha) } Z_{\alpha_i} \text{ for all } \alpha \in E^{\st} \text { with } r(\alpha)=v,
 \end{equation}
where each $Y_{\alpha_i}$ and $Z_{\alpha_i}$ is a subspace of $X_v$ with an infinite countable basis, 
unless in the first equation $v$ is a sink and in the second $v$ is a source. Note that in Equation~\ref{decio} we only consider a copy $Y_{\alpha_i}$  when $\alpha_i\not = 0$, i.e., $i\leq w(\alpha)$. 
Now for each $v\in E^0$, define $\Pp_v \in \End_R(X)$ as projection of $X$ onto $X_v$. Clearly $\Pp_v\Pp_v=\Pp_v$ and $\Pp_v\Pp_u=0$ for $v\not = u$. This immediately implies the set $\{\Pp_v, v \in E^0\}$ is $R$-linearly independent in $\End_R(X)$.

For $\alpha \in E^{\st}$ with $1 \leq i \leq w(\alpha)$, define $\A_i$ as follows 
\begin{equation}\label{mnft}
\A_i:X \longrightarrow Z_{\alpha_i} \stackrel{\theta_{\A_i}}{\longrightarrow} Y_{\alpha_i} \longrightarrow X,
\end{equation} where the first map is the projection, $\theta_{\A_i}$ is an isomorphism (which exists, as $Z_{\alpha_i}$ and $Y_{\alpha_i}$ are vector spaces with bases of the same cardinality), and the last map is the injection. Similarly define  $\A_i^*$ as follows \[\A_i^*:X \longrightarrow Y_{\alpha_i} \stackrel{\theta_{\A_i}^{-1}}{\longrightarrow} Z_{\alpha_i} \longrightarrow X.\]

The maps $\{ \A_1,\dots \A_{w(\alpha)} \mid \alpha \in E^{\st} \}$ and $\{ \A_1^*,\dots \A_{w(\alpha)}^* \mid \alpha \in E^{\st} \}$ satisfy the relations in Definition~\ref{iLPA} of a weighted Leavitt path algebra. We check these relations:

Let $\alpha \in E^{\st}$ with $s(\alpha)=v$ and let $1\leq i \leq w(\alpha)$. We will show that $\Pp_v\A_i=\A_i$. For $x\in X$, the definition of $\A_i$ in (\ref{mnft}) shows that $\A_i(x) \in Y_{\alpha_i} \subseteq X_v$. Since $\Pp_v$ is the projection on $X_v$, i.e., its restriction to $X_v$ is the identity, we get   $\Pp_v\A_i=\A_i$. The other relations of type (2) in Definition~\ref{iLPA} follow similarly. 

We now verify type (3) relations in the definition.  Fix $v\in E^0$ and consider $\{ \alpha \in E^{\st} \mid s(\alpha)=v\}$. For $x \in X_v$, by Equation~\ref{decio}, one can write $x$ as a direct sum $x=\sum_{\{s(\alpha)=v\}} y_{\alpha_i}$,  where $y_{\alpha_i} \in Y_{\alpha_i}$ and $i$ is fixed. Now $\A_i^*(x)=\A_i^*(\sum y_{\alpha_i})=\theta_{\alpha_i}^{-1}(y_{\alpha_i})$ and so 
$\A_i\A_i^*(\sum y_{\alpha_i})=\theta_{\alpha_i}(\theta_{\alpha_i}^{-1}(y_{\alpha_i}))=y_{\alpha_i}$. This implies that 
$\sum_{\{s(\alpha)=v\}}\A_i\A_i^*(x)=\sum y_{\alpha_i}=x$. Clearly if $x \not \in X_v$ then $\sum\A_i\A_i^*(x)=0$. Putting these together we get $\sum\A_i\A_i^*(x)=\Pp_v$. On the other hand choose fixed $i$ and $j$ such that $i \not = j$. Write  $x=\sum_{\{s(\alpha)=v\}} y_{\alpha_j}$, where $y_{\alpha_j} \in Y_{\alpha_j}$. Then as in above, 
$\A_j^*(x)=\A_j^*(\sum y_{\alpha_j})=\theta_{\alpha_j}^{-1}(y_{\alpha_j})\in Z_{\alpha_j}$. But by Equation~\ref{decio2}, $X_{r(\alpha)}=\bigoplus Z_{\alpha_i}$, where $1 \leq i \leq w(\alpha)$.  So $\A_i(\A_j^*(x))=0$ and therefore $\sum_{\{s(\alpha)=v\}}\A_i\A_j^*(x)=0$ for any $x\in X_v$ and thus for any $x\in X$. This shows that the maps $\{ \A_1,\dots \A_{w(\alpha)} \mid \alpha \in E^{\st} \}$ and $\{ \A_1^*,\dots \A_{w(\alpha)}^* \mid \alpha \in E^{\st} \}$ satisfy relation (3) in Definition~\ref{iLPA}. The rest is similar.

Now because of the universality of $\LL_R(E,w)$, there is a $R$-homomorphism $\Phi: \LL_R(E,w) \rightarrow \End_R(X)$ such that
$\Phi(v)=\Pp_v$, $\Phi(\alpha_i)=\A_i$ and $\Phi(\alpha_i^*)=\A_i^*$. Since $\{\Pp_v, v \in E^0\}$ are nonzero elements in $\End_R(X)$ and are $R$-linearly independent, then the set of vertices in $\LL_R(E,w)$ should be so as well.  
  \end{proof}

\begin{deff}\label{wewe} For a directed graph $E$, define the {\it opposite graph}, $E^{\op}$ as a graph with the same set of vertices and edges as $E$ (for  an edge $\alpha$ in $E$, denote the corresponding edge in $E^{\op}$ with $\alpha^{\op}$), such that for an edge $\alpha^{\op}$ in $E^{\op}$, $s(\alpha^{\op})=r(\alpha)$ and $r(\alpha^{\op})=s(\alpha)$. This means that $E^{\op}$ is obtained from $E$ by simply reversing the arrows. 
\end{deff}

It is not clear in general how the algebras  $\LL(E)$ and $\LL(E^{\op})$ are related. For example, for 
\[
\xymatrix{
E:\bullet \ar[r] & \bullet  \ar@/^/[r]  \ar@/_/[r] & \bullet 
}
\]
one obtains $\LL_K(E)\cong_{\gr}\M_5(K)(0,1,1,2,2)$, whereas, for 
\[
\xymatrix{
E^{\op}:\bullet  & \bullet  \ar[l] & \bullet \ar@/^/[l]  \ar@/_/[l]
}
\] we have $\LL_K(E^{\op})\cong_{\gr}\M_4(K)(0,1,2,2)$ (see Theorem~\ref{acyclicc}, see also Example~\ref{opex}).

\begin{deff}

For a directed graph $E$, the {\it weighted graph associated to $E$}, denoted by $E_w$, is obtained by considering all the edges with the same source and the same range in $E$ as one structured edge with appropriate weight (i.e., the number of these edges) in $E_w$. More formally, $E_w$ has $E^0$ as the set of vertices and if in $E$, for any $u,v\in E^0$,  
$s^{-1}(u)\cap r^{-1}(v) \subseteq E^1$ is non-empty, then there is $\alpha\in E_w^s$ with $s(\alpha)=u$, $r(\alpha)=v$, and 
$w(\alpha)=|s^{-1}(u)\cap r^{-1}(v)|$. 
\end{deff}

\begin{example}
Consider the graph $E$ with two vertices and with no loops
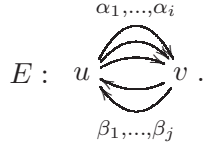
\[
\xymatrix{
E:  \!\!\!\!\!\!\!\!&\!\!\!\! u  \ar@/^/[r] \ar@/^1pc/[r] \ar@/^1.5pc/[r]^{\alpha_1,\dots ,\alpha_i}  & v \ar@/^/[l]  \ar@/^1.2pc/[l]^{\beta_1,\dots,\beta_j}
}.
\] Then one can see that the map $E\longrightarrow E^{\op}, (u\mapsto v, v\mapsto u, \alpha_i \mapsto \alpha_i^{\op}, \beta_i \mapsto \beta_i^{\op})$ induces an isomorphism on the level of LPAs, i.e., $\LL(E) \cong \LL(E^{\op})$. Now consider the associated weighted graph $E_w$ of $E$ (see Definition~\ref{wewe}).  The map $E_w \longrightarrow E^{\op},  (u\mapsto u, v\mapsto v, \alpha_i \mapsto \alpha_i^*, \beta_i \mapsto \beta_i^*)$ gives the isomorphism, $\LL(E_w,w)\cong \LL(E^{\op})$.  In the same manner, one can see that for any graph $E$ with one vertex, $\LL(E)\cong \LL(E^{\op})\cong \LL(E_w,w)$, i.e., 
\[\xymatrix{
\LL\big(\!\!\!\!\!\!\!\!\!\!\!\!\!& \bullet\ar@(ul,ur)^{\alpha_{1},\dots,\alpha_{n}}
},w\big)\cong
\xymatrix{
\LL\big(\!\!\!&   \bullet \ar@{.}@(l,d) \ar@(ur,dr)^{\alpha_{1}} \ar@(r,d)^{\alpha_{2}} \ar@(dr,dl)^{\alpha_{3}} 
\ar@(l,u)^{\alpha_{n}}& 
}\big).\]

Also note that these isomorphisms are not graded as in the latter example $\LL(E_w,w)$ is $\textstyle{\bigoplus_n} \mathbb Z$-graded, whereas $\LL(E)$ is just $\mathbb Z$-graded. 

\end{example}

Consider the category $\mathcal G^w$ with objects all row-finite weighted graphs and morphisms, the {\it complete weighted graph homomorphisms}, i.e., a morphism $f:E\rightarrow F$ consists of a map $f^0:E^0 \rightarrow F^0$ and $f^1:E^{\st} \rightarrow F^{\st}$ such that $r(f^1(\alpha))=f^0(r(\alpha))$, $s(f^1(\alpha))=f^0(s(\alpha))$ and $w(\alpha)=w(f^1(\alpha))$
for any $\alpha \in E^{\st}$, additionally, $f^0$ is injective and $f^1$ restricts to a bijection from $s^{-1}(v)$ to $s^{-1}(f^0(v))$ for every $v\in E^0$ which emits edges.   One can check that a morphism $f$ preserves the relations in Definition~\ref{iLPA}, and thus induces a graded ring homomorphism $\LL_R(E,w)\rightarrow \LL_R(F,w)$. Thus, when $R$ is commutative,  we have a functor $\LL:\mathcal G^w \rightarrow \mathcal A$, where $\mathcal A$ is the category of (non-unital) $R$-algebras. A {\it weighted subgraph} $X$ of the weighted graph $E$, is a weighted graph $X=(X^0,X^{\st},X^1,r_X,s_X,w_X)$ such that $X^0 \subseteq E^0$, $X^{\st} \subseteq E^{\st}$, and $s_X,r_X$ are the restrictions of   $s_E,r_E$ on $X^{\st}$ and $w_X(\alpha)=w_E(\alpha)$ for any $\alpha \in X^{\st}$. 
A weighted subgraph $X$ is called {\it complete}, if $x\in X^0$ and $s^{-1}_X(v)\not = \emptyset$, then $s^{-1}_X(v)=s^{-1}_E(v)$. In this case the inclusion map $X \hookrightarrow E$ is a complete graph homomorphism.

In the case of multi-headed weighted rose graphs, a statement similar to Theorem~\ref{polyhead} can be obtained. For this, we need a weighted version of multi-headed rose graphs. 

\begin{deff}\label{wpetaldef}
A {\it weighted rose with $k$-petals}, is a weighted graph which consists of one vertex, $k$ structured loops and the weight map $w$. We denote this weighted graph by $L_{(k,w)}$. If $w$ is a constant map $1$, then the weighted rose reduces to the usual $L_k$ (see Definition~\ref{petaldef} and Example~\ref{lpa}) and if there are $n+k-1$ petals and $w$ is the constant map assigning $n$ to each petal, then the wLPA of $L_{(n+k-1,w)}$ is the Leavitt algebra of type $(n,k-1)$, denoted by $\LL(n,k)$ (see Example~\ref{lat} and compare this with Definition~\ref{petaldef} for consistency). The weighted Leavitt path algebra of $L_{(n,w)}$ with the coefficients in $R$ is denoted by $\LL_R{(n,w)}$
\end{deff}

\begin{deff} \label{wpopyt}
A {\it  multi-headed weighted rose} graph $E$ consists of an unweighted  finite acyclic graph $E_1$ with sinks $\{v_1,\dots,v_t\}$ together with weighted  $n_s$-petal graphs $L_{(n_s,w_s)}$, $1\leq s \leq t$, attached to $v_s$, where  $n_s \in \mathbb N$. If $n_s=0$ for all $1\leq s \leq t$, then the graph is  finite acyclic. Note that by definition the weighted map $w$ of $E$  is $w(\alpha)=1$ if $\alpha \in E_1^{\st}$ and $w(\alpha)=w_s(\alpha)$ if $\alpha \in L_{(n_s,w_s)}^{\st}$. 
\end{deff}

\begin{theorem}\label{wpolyhead}
Let   $E$ be a  multi-headed weighted rose  graph consisting of an acyclic graph $E_1$ with sinks $\{v_1,\dots,v_t\}$ which are attached to $L_{(n_1,w_1)},\dots,L_{(n_t,w_t)}$, respectively.
For any  $v_s$,  let $\{p^{v_s}_i \mid 1\leq i \leq n(v_s)\}$ be the set of all paths in $E_1$ which end in $v_s$.  Let $l=\max\{w_s(L_{(n_s,w_s)}^{\st}) \mid 1\leq s \leq t\}$. 

\begin{enumerate}[\upshape(1)]

\item If $R$ is a ring then there is a $\bigoplus_l \mathbb Z$-graded isomorphism
\begin{equation}\label{wdampai}
\LL_R(E,w) \cong_{\gr} \bigoplus_{s=1}^t \M_{n(v_s)} \big(\LL_R(n_s,w_s)\big)\big (|p^{v_s}_1|,\dots, |p^{v_s}_{n(v_s)}|\big ),
\end{equation}
where $|p^{v_s}_i|$ denotes the element $(|p^{v_s}_i|,0,\dots,0) \in \bigoplus_l \mathbb Z$.

\medskip
\item If $R=K$ is a field and  for any $1\leq s \leq t$, $L_{(n_s,w_s)}=L_{(n,w_n)}$, then
\begin{equation}\label{wxxdampai}
\LL_K(E,w)\cong_{\gr} \LL_K(E_1)\otimes_K \LL_K(n,w_n).
\end{equation}

%\item Let $R=K$ be a field. Then  $\LL_K(E)$ is a group ring
%$$\LL_K(E)\cong_{\gr} \bigoplus_{s=1}^t \M_{n(v_s)}(K)[\mathbb Z],$$
%if and only if for any $1\leq s \leq t$, $L_{n_s}=L_1$. 

\end{enumerate}
\end{theorem}
\begin{proof}
The proofs are quite similar to the ones in Theorem~\ref{polyhead} and~\ref{polyheadnew} and we give a sketch. 

(1) Consider the obvious complete weighted graph homomorphisms from  $E_1$ and $L_{(n_s,w_s)}$ to $(E,w)$ and identify $\LL_R(E_1,w)$ and $\LL_R(n_s,w_s)$ with their images in $\LL_R(E,w)$, respectively. 
Define 
\[\I_{v_s}=\Big \{ \sum k \alpha x  \beta^* \mid k \in R, \, \alpha, \beta \in E_1^*, r(\alpha)=v_s=r(\beta), \, x\in \LL_R(n_s,w_s) \Big \} \subseteq \LL_R(E,w).\]   
We observe that $\I_{v_s}$ is an ideal of $\LL_R(E,w)$. It is enough to check that for monomials $\alpha x \beta^* \in \I_{v_s}$ and $\theta \in \LL_R(E,w)$,   $\theta \alpha x \beta^*$ and $\alpha x \beta^* \theta$ are in $ \I_{v_s}$. Let $\theta \alpha x \beta^*\not = 0$. Then one can decompose $\theta$ as $\gamma y \delta ^*$, where $\gamma$ and $\delta$ are paths in $E_1$ with $r(\gamma)=r(\delta)=v_s$ and $y\in \LL(n_s,w_s)$. Now a similar argument as in the proof of Theorem~\ref{polyhead} shows that $\I_{v_s}$'s are  ideals of $\LL_R(E,w)$, and furthermore $\LL_R(E,w)$ is a direct sum of such ideals, and that $\I_{v_s}$ are graded isomorphic to $\M_{n(v_s)}(\LL_R(n_s,w_s))$. 

(2) The proof is similar to part (1) of  Theorem~\ref{polyheadnew}. 
\end{proof}

\begin{example}

By Theorem~\ref{wpolyhead} and Example~\ref{lat} we have   

%& \bullet \ar@{.}@(l,d) \ar@(ur,dr)^{y_{11},\dots,y_{n1}} \ar@(r,d)^{y_{12},\dots,y_{n2}} \ar@(dr,dl)^{y_{13},\dots,y_{n3}} \ar@(l,u)^{y_{1,n+k},\dots,y_{n,n+k}}& 

\[
 \xymatrix{
& & 
\bullet \ar@{.}@(l,d) \ar@(ur,dr)^{y_{11},\dots,y_{n1}} \ar@(r,d)^{y_{12},\dots,y_{n2}} \ar@(dr,dl)^{\, \, y_{13},\dots,y_{n3}} \ar@(l,u)^{y_{1,n+k},\dots,y_{n,n+k}}& \\
\LL_K\Big( \; \bullet \ar[r]  & \bullet   \ar@<1.5pt>[r]  \ar@<-1.5pt>[r] \ar@<0.5ex>[ur] \ar@<-0.5ex>[ur] \ar@<0ex>[ur] \ar[dr] &  \bullet  \ar@(dl,dr)  \ar@(d,r) \ar@(dr,ur)
 &  ,w\Big)  \cong_{\gr} \M_3(K) \times \M_5(\LL_K(1,3)) \times \M_7(\LL_K(n,k+1)).
 \\ & & \bullet  }
\]
%with the shifted grading $(0,1,2), (0,1,1,2,2)$ and $(0,1,1,1,2,2,2)$ respectively. 
\end{example}

One can extend Definition ~\ref{mmjjsspp} and  Theorem~\ref{tensoir} in the setting of weighted Leavitt path algebras.  

\begin{deff}
Let $E$ be an acyclic graph with a unique sink and $F$ be a weighted graph. Denote by $E\otimes F$ a graph obtained by attaching $E$ from the sink  to each vertex of $F$.
\end{deff}

\begin{lemma}\label{tensoir5}
Let $K$ be a field, $E$ be an acyclic graph with a unique sink and $F$ be a weighted graph. Then \[\LL_K(E\otimes F) \cong_{\gr} \LL_K(E) \otimes_K \LL_K(F).\]
\end{lemma}
\begin{proof}
The proof is quite similar to Theorem~\ref{tensoir} so we omit the argument. 
\end{proof}

\begin{example}
Consider the oriental line graphs $E_1$ and $E_2$ with $n$ and $k+1$ vertices respectively and the graph $E$ consisting of $n(k+1)$ loops, all of weight $n$. Then by Lemma~\ref{tensoir5}
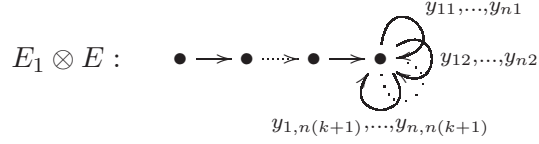
\begin{equation*}
\xymatrix@=13pt{
E_1\otimes E: &  \bullet \ar[r] &  \bullet \ar@{.>}[r] &  \bullet \ar[r]& \bullet \ar@(u,r)^{y_{11},\dots,y_{n1}} \ar@(ur,dr)^{y_{12},\dots,y_{n2}}  \ar@{.}@(r,d) \ar@(dr,dl)^{y_{1,n(k+1)},\dots,y_{n,n(k+1)}}&
}
\end{equation*}
and $\LL(E_1\otimes E)\cong \LL(E_1) \otimes \LL(E)=\M_n(A)$, where $A=V_{n,n(k+1)}$.
Furthermore $E_2\otimes E_1 \otimes E$ is

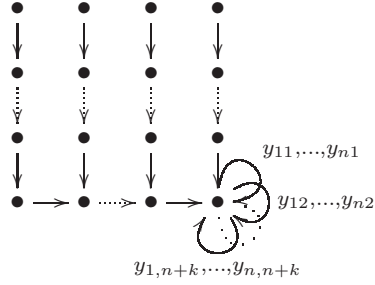
\begin{equation*}
\xymatrix@=13pt{
&\bullet \ar[d] &  \bullet \ar[d] &   \bullet \ar[d] &   \bullet \ar[d]    &\\
&\bullet \ar@{.>}[d] &   \bullet \ar@{.>}[d] &    \bullet \ar@{.>}[d] &   \bullet \ar@{.>}[d]  &\\
& \bullet\ar[d]  & \bullet\ar[d] &  \bullet\ar[d] &   \bullet\ar[d] &\\
&   \bullet \ar[r] &   \bullet \ar@{.>}[r] & \bullet\ar[r]  & \bullet \ar@(u,r)^{y_{11},\dots,y_{n1}} \ar@(ur,dr)^{y_{12},\dots,y_{n2}}  \ar@{.}@(r,d) \ar@(dr,dl)^{y_{1,n+k},\dots,y_{n,n+k}}&
}
\end{equation*}
and $\LL(E_2\otimes E_1\otimes E)\cong \LL(E_2) \otimes \LL(E_1\otimes E)=\M_{k+1}(A)\otimes \M_n(A)\cong\M_{n(k+1)}(A)$. But because $A=V_{n,n(k+1)}$, we have $A^n \cong A^{n(k+1)}$.
Considering the endomorphism rings, we get \[\M_n(A)\cong \M_{n(k+1)}(A)\cong \M_n(\M_{k+1}(A)).\]  But $A$ is not isomorphic to $\M_{k+1}(A)$ as the first ring is a domain. This gives an example of two rings $R$ and $S$ such that $\M_n(R)\cong \M_n(S)$ for some $n$, but $R\not \cong S$ (see \cite[p.225]{cohn11}).
\end{example}

\begin{lemma} \label{catepf} \hfill

\begin{enumerate}[\upshape(1)]

\item Every row-finite weighted graph is a direct limit of a directed system of finite weighted graphs. 

\item Every weighted Leavitt path algebra is a direct limit of weighted Leavitt path algebras  corresponding to finite weighted graphs. 
\end{enumerate}
\end{lemma}
\begin{proof}
(1) Let $E$ be a row-finite weighted graph. First note that a union of complete subgraphs of $E$ is again a complete subgraph of $E$. Moreover, $E$ is the union of finite subgraphs (not necessarily complete).  Finally, any finite subgraph is contained in a finite complete subgraph. For, if $X$ is a finite subgraph of $E$, then consider the subgraph $Y$ of $E$ as follows:
$Y^0=X^0\cup \{r(\alpha)\mid \alpha \in E^{\st}, s(\alpha)\in X^0 \}$ and $Y^{\st}=\{\alpha \in E^{\st} \mid s(\alpha)\in X^0\}$. One can easily see that $Y$ is a complete subgraph of $E$. Putting these three facts together, it follows that $E$ is the direct limit of the directed system of its finite complete subgraphs. 

(2) Let $E=\varinjlim X_i$, where $\{X_i\}$ is the directed system of complete subgraphs of $E$, by (1). Then it is easy to observe that $\LL(E,w)\cong \varinjlim\LL(X_i,w)$. 
\end{proof}

\begin{remark}
To any weighted graph, one can associate a directed graph by simply considering the weight as the number of edges connecting the adjacent vertices. (In Example~\ref{calcu}, $E'$ is the directed graph obtained from the weighted graph $E$). One can check that this defines a (forgetful) functor  $n : \mathcal G^w \rightarrow \mathcal G$. It is not known whether there is a functor which relates the corresponding wLPA to LPA, i.e., whether there is a functor such that the  following diagram is commutative
\[\xymatrix{
\mathcal G^w \ar[r]^{\LL(-,w)} \ar[d]_{n} & \mathcal A \ar@{.>}[d]^?\\
\mathcal G \ar[r]_{\LL(-)} & \mathcal A 
} \]
In the same manner, recall that one can associate to a graph a weighted graph (see Definition~\ref{wewe}), so a similar question can be raised here too:
\[
\xymatrix{
\mathcal G \ar[r]^{\LL(-)} \ar[d]_{w} & \mathcal A \ar@{.>}[d]^?\\
\mathcal G^w \ar[r]_{\LL(-,w)} & \mathcal A 
}  
\]

\end{remark}

For a ring $A$ with identity, the monoid $\V(A)$ is defined as the set of isomorphism classes of finitely generated projective left $A$-modules equipped with the direct sum as the binary operation. When $A$ is not unital, i.e., does not have identity, one defines $\V(A)$ as the set of equivalent classes of idempotents in $M_{\infty}(A)$ with $[e]+[f]=
\left[ \begin{array}{ccc} 
e & 0  \\ 
0 & f 
\end{array} \right]$, where $e \sim e'$ if there are $x,y \in M_{\infty}(A)$ such that $e=xy$ and $e'=yx$. Here $M_{\infty}(A)$ are matrices over $A$ with finitely many nonzero entries. There is a corresponding construction based on finitely generated projective modules as well, see~\cite[p.163]{amp}. In~\cite[Theorem~3.5]{amp}, Ara, Moreno and Pardo show that for a directed graph $E$, $\V(\LL_K(E))$ coincides with a monoid naturally constructed from the graph $E$ and further this monoid is a refinement monoid and thus separative \cite[Theorem~6.3]{amp}. A similar construction is valid in the setting of wLPA, however we will see this monoid is neither refinement nor separative in general. Recall that in an abelian monoid $M$, for $x,y \in M$, we denote $y\leq x$  if there is $z\in M$ such that $x=y+z$. Then $M$ is called {\it separative} if for elements $x,y,z \in M$,  $x+z=y+z$ and $z\leq nx$ and $z\leq ny$, for some positive integer $n$, implies that $x=y$. $M$ is called a {\it refinement monoid} if  $x_1+x_2=y_1+y_2$, $x_1,x_2,y_1,y_2 \in M$, then there are $z_{i,j}, 1 \leq i,j \leq 2$ such that $x_i=z_{i1}+z_{i2}$ and $y_j=z_{1j}+z_{2j}$, for $1 \leq i,j \leq 2$.

In the following theorem we will use Bergman's machinery \cite[p. 38 and Theorem~3.3]{berg1}: Let $A$ be a $K$-algebra and $P$ and $Q$ be finitely generated projective $A$-modules. Then there is a $K$-algebra $B:=A\langle i,i^{-1}:\overline P \cong \overline Q \rangle$, with an algebra homomorphism $A\rightarrow B$ such that there is a universal isomorphism $i:\overline P \rightarrow \overline Q$, where $\overline M=B\otimes_A M$ for a left $A$-module $M$. Then Bergman's Theorem~5.2 in \cite{berg1} states that $\V(B)$ is the quotient of $\V(A)$ modulo the  relation $[P]=[Q]$.

\begin{theorem}\label{kspcts} Let $K$ be a field and $E$ be a  weighted graph. Let $M_E$ be the abelian monoid generated by $\{v \mid v \in E^0 \}$ subject to the relations 
\begin{equation}\label{phgqcu}
n_vv=\sum_{\{\alpha\in E^{\st} \mid s(\alpha)=v \}} r(\alpha),
\end{equation}
for every $v\in E^0$ that emits edges, where $n_v= \max\{w(\alpha)\mid \alpha\in E^{\st}, \, s( \alpha)=v\}$. Then there is a natural monoid isomorphism $\V(\LL_K(E,w))\cong M_E$. Furthermore, if $E$ is finite, then $\LL_K(E,w)$ is hereditary. 
\end{theorem}
\begin{proof}
Define a map $\psi:E^0\rightarrow \V(\LL_K(E,w))$ by $\psi(v)=[v]$ and extend this to the map from the free monoid on $E^0$ to $\V(\LL_K(E,w))$. This induces a map $\psi_E:M_E\rightarrow \V(\LL_K(E,w))$. To see this, we need to show that if $v$ emits edges, then $n_vv$ and $\sum_{\{\alpha\in E^{\st} \mid s(\alpha)=v \}} r(\alpha)$ map to the same element in $\V(\LL_K(E,w))$, where $n_v= \max\{w(\alpha)\mid \alpha\in E^{\st}, \,s( \alpha)=v\}$. Let $\{\alpha_1,\dots,\alpha_s\}$ be all the structured edges which are emitted from $v$. Consider the matrices $Y=(\alpha_{ij})_{1\leq j \leq s, 1\leq i\leq w(\alpha_j)}$, where $\alpha_{ij}={(\alpha_j)}_i \in E^1$ and $X=(Y^*)^t$, where ${}^t$ is the transpose operation (this is a similar arrangement as in ~(\ref{breaktr})). Then the conditions of Definition~\ref{iLPA} guarantee  that $Y.X=n_v[v]$ and $X.Y=\sum_{\{\alpha\in E^{\st} \mid s(\alpha)=v \}} [r(\alpha)]$. So $\psi_E$ is well-defined. But $\LL_K(E,w)$ is a direct limit of graph algebras corresponding to finite graphs (see Lemma~\ref{catepf}). Thus it is enough to prove $\psi_E$ is an isomorphism for a finite graph $E$. So let $E$ be a finite graph and $\{v_1,\dots,v_m\}$ be the set of vertices which emit edges. Let $A_0=\prod_{v\in E^0}K$. Consider the following two finitely generated projective $A_0$-modules, $P=n_{v_1}(A_0v_1)$  and $Q=\textstyle{\bigoplus_{\{\alpha\in E^{\st} \mid s(\alpha)=v_1 \}}}A_0r(\alpha)$. Using the Bergman's machinery, there exists an algebra $A_1=A_0 \langle i,i^{-1}: \overline P \cong \overline Q\rangle$ with a universal isomorphisms \[i: \overline P:=A_1\otimes_{A_0} P \rightarrow\overline Q:=A_1\otimes_{A_0} Q.\] In fact  this algebra is $\LL(X_1,w)$, where $X_1$ is a graph with the same vertices as $E$ and where $v_1$ emits the same structured edges (thus the same edges) as in $E$ and other vertices do not emit any edges. Namely, if $\{\alpha_1,\dots,\alpha_s\}$ is all the structured edges which emit from $v_1$ then the right multiplication by  the matrix $Y=(\alpha_{ij})_{1\leq j \leq s, 1\leq i\leq w(\alpha_j)}$, where $\alpha_{ij}={(\alpha_j)}_i \in E^1$, gives the map 
\[i:\overline P=n_{v_1}(A_1 v_1) \rightarrow \overline Q =\textstyle{\bigoplus_{\{\alpha\in E^{\st} \mid s(\alpha)=v_1 \}}} A_1r(\alpha),\] and $X=(Y^*)^t$, where ${}^t$ is the transpose operation gives $i^{-1}$. Now \cite[Theorem~5.2]{berg1} asserts that $\V(A_1)$ is obtained from $\V(A_0)$ by adding the relation $[P]=[Q]$.  Translating this to our setting, we get that $\V(A_1)$ is the monoid generated by the set $\{[v] \mid v\in E^0\}$ subject to the relation $n_{v_1} [v_1]=\sum_{\{\alpha\in E^{\st} \mid s(\alpha)=v_1 \}} [r(\alpha)].$ 

We repeat this process to cover the whole graph. To be precise, let $A_k=\LL(X_k,w)$, $k\geq 1$, where $X_k$ is the graph with the same vertices as $E$, but only the first $k$ vertices $\{v_1,\dots,v_k\}$ emit structured edges. By induction, $\V(A_k)$ is an abelian group generated by $\{[v] \mid v\in E^0\}$ subject to the relation 
$n_{v_i} [v_i]=\sum_{\{\alpha\in E^{\st} \mid s(\alpha)=v_i \}} [r(\alpha)]$, where $1\leq i \leq k$. Then  $A_{k+1}=A_k \langle i,i^{-1}: \overline P \cong \overline Q\rangle$ with $P=n_{v_{k+1}}(A_{k} v_{k+1})$ and $Q =\textstyle{\bigoplus_{\{\alpha\in E^{\st} \mid s(\alpha)=v_{k+1}\}}}A_k r(\alpha)$.  So by \cite[Theorem~5.2]{berg1}, $\V(A_{k+1})$ is the monoid generated by all the vertices of $E$ subject to relations corresponding to $\{v_1,\dots,v_{k+1}\}$. Now Bergman's theorem also implies that the global dimension of $A_{k+1}$ is the same as the global dimension of $A_{k}$, therefore they are all equal to the global dimension of $A_0$ which, being a semisimple algebra, has dimension  zero. This implies that $\LL_K(E,w)$ is hereditary. 
\end{proof}

As mentioned above, the monoids associated to Leavitt path algebras are refinement and separative monoids. This is not the case for weighted Leavitt path algebras. Consider the following weighted graphs: 
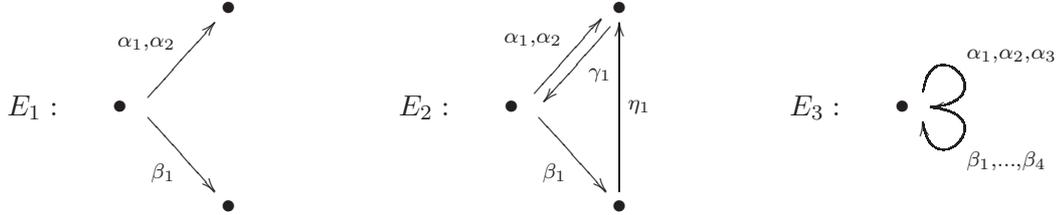
\[ 
\xymatrix{
&& \bullet &&&& \bullet \ar@<0.5ex>[dl]^{\gamma_1} && \\
E_1: & \! \! \! \! \! \! \! \! \bullet \ar[ur]^{\alpha_1,\alpha_2} \ar[dr]_{\beta_1} &&& E_2: & \! \! \! \! \! \! \! \!  \bullet \ar@<0.5ex>[ur]^{\alpha_1,\alpha_2}  \ar[dr]_{\beta_1} &&&E_3: & \! \! \! \! \! \! \! \! 
 \bullet  \ar@(u,r)^{\alpha_{1},\alpha_2,\alpha_3} \ar@(r,d)^{\beta_1,\dots,\beta_4} & \\
&& \bullet &&&& \bullet \ar[uu]_{\eta_1} 
}
\]
One can easily show that $M_{E_1}=\big \langle (1,0),(0,1),(1/2,1/2) \big \rangle \subseteq \mathbb Q\times \mathbb Q$ which is not a refinement monoid but it is separative, $M_{E_2}\cong \mathbb N$ which it is a refinement and separative monoid, and $M_{E_3}$ which is not even separative.

 We mention here (without proof) that although $M_E$ is not, in general, a refinement monoid, but it is  {\it weighted refinement}. Recall from Theorem~\ref{kspcts} that $M_E$ is a monoid obtained from the free monoid $F$ generated by the set of vertices of $E$ subject to the equivalence relation $\thicksim$ generated by the relation~(\ref{phgqcu}). 
Then one can prove that, if  $\alpha_1+\alpha_2 \thicksim \beta_1+\beta_2$, where $\alpha_1,\alpha_2,\beta_1,\beta_2 \in F$ and $\alpha_1$, $\beta_1$ are {\it weighted elements}, i.e., if $v$ appears in the expression of $\alpha_1$ then $n_v v$ also appears in the expression, where  $n_v=\max\{w(\alpha)\mid \alpha\in E^{\st}, \, s( \alpha)=v\}$,
then there are $\mu_{i,j}, 1 \leq i,j \leq 2$ such that $\alpha_i=\mu_{i1}+\mu_{i2}$ and $\beta_j=\mu_{1j}+\mu_{2j}$, for $1 \leq i,j \leq 2$.

 %The support of an element $\alpha$ in $F$, denoted by $\Supp(\alpha)$, is the set of basis elements appearing in the expression of $\alpha$. We call $\alpha$ a {\it weighted element} if $x\in\Supp(\alpha)$ then $n_x x \in \Supp(\alpha)$. A weighted monoid is {\it weighted refinement}, 

In recent years there have been several algebraic constructions motivated by algebras introduced by Leavitt~\cite{vitt56,vitt57,vitt62}. Very  recently Ara and Goodearl~\cite{aragood} have introduced Leavitt path algebras associated to separated graphs. The idea is to partition the edges emitted from each vertex and put the relations (3) and (4) in Definition~\ref{llkas} on each partition separately. Their construction covers the usual Leavitt path algebras when the graph is ``not separated'', and this is the only instance that their construction coincides with weighted Leavitt path algebras when the weight map is the constant 1 (see Example~\ref{lpa}). 
\end{example}

We finish the paper by determining the Grothendieck group, $K_0$, of weighted Leavitt path algebras. When the graph is unweighted, i.e., all edges have weight one, this recovers $K_0$ group of Leavitt path algebras. 

For an abelian monoid $M$, we denote by $M^{+}$ the group completion of $M$. This gives a left adjoint functor to the forgetful functor from the category of abelian groups to abelian monoids.
Let $F$ be a free abelian monoid generated by a countable set $X$. The nonzero elements of $F$ can be written as $\sum_{t=1}^n x_t$ where $x_t \in X$. Let $r_i, s_i$, $i\in I \subseteq \mathbb N$, be elements of $F$. We define an equivalence relation on $F$ denoted by $\langle r_i=s_i\mid i\in I \rangle$ as follows: Define a binary relation $\rightarrow$ on $F\backslash \{0\}$,  $r_i+\sum_{t=1}^n x_t \rightarrow s_i+\sum_{t=1}^n x_t$, $i\in I$ and generate the equivalence relation on $F$ using this binary relation. Namely, $a \sim a$ for any $a\in F$ and for $a,b \in F \backslash \{0\}$, $a \sim b$ if there is a sequence $a=a_0,a_1,\dots,a_n=b$ such that for each $t=0,\dots,n-1$ either $a_t \rightarrow a_{t+1}$ or $a_{t+1}\rightarrow a_t$. We denote the quotient monoid by
$F/\langle r_i=s_i\mid i\in I\rangle$.
Then one can see that there is a canonical group isomorphism
\begin{equation} \label{monio}
\Big (\frac{F}{\langle r_i=s_i\mid i\in I \rangle}\Big)^{+} \cong \frac{F^{+}}{\langle r_i-s_i\mid i\in I \rangle}.
\end{equation}

\begin{deff}\label{adji}
Let $(E,w)$ be a weighted graph and  $N'$ be the adjacency  matrix
$(n_{ij}) \in \mathbb Z^{E^0\oplus E^0}$ where $n_{ij}$ is the
number of structured edges from $v_i$ to $v_j$. Furthermore let
$I'_w$ be the weighted identity matrix defined as $(a_{ij})$, where
$a_{ij}=0$ for $i\not = j$ and $a_{ii}=n_{v_i}$.
Here $n_v= \max\{w(\alpha)\mid \alpha\in
E^{\st}, \, s( \alpha)=v\}$ for every
$v\in E^0$ that emits edges and zero otherwise.
Let $N^t$ and $I_w$ be
the matrices obtained from $N'$ and $I'_w$ by first taking the transpose and then removing the columns
corresponding to sinks, respectively.
\end{deff}
Clearly the adjacency matrix depends on the ordering we put on
$E^0$. We usually fix
 an ordering on $E^0$ such that the elements of
$E^0 \backslash \sink$ appear first in the list follow with elements
of the $\sink$.

 Multiplying the matrix $N^t-I_w$ from the left defines a
homomorphism $\mathbb Z^{E^0\backslash \sink} \longrightarrow
\mathbb Z^{E_0}$, where  $\mathbb Z^{E^0\backslash \sink} $ and
$\mathbb Z^{E^0}$ are the direct sum of copies of $\mathbb Z$
indexed by $E^0\backslash \sink$ and $E^0$, respectively. The next
theorem shows that the cokernel of this map gives the Grothendieck
group of weighted Leavitt path algebras. Clearly, for the graph with
weights one (unweighted), this recovers $K_0$ group of Leavitt path
algebras.

\begin{theorem}\label{wke}
Let $K$ be a field, $(E,w)$ be a  weighted graph  and $\LL(E,w)$ be a weighted Leavitt path algebra. Then
\begin{equation}
K_0(\LL(E,w))\cong \coker\big(N^t-I_w:\mathbb Z^{E^0\backslash \sink} \longrightarrow \mathbb Z^{E^0}\big).
\end{equation}
\end{theorem}
\begin{proof}
Let $M_E$ be the abelian monoid generated by $\{v \mid v \in E^0 \}$ subject to the relations
\begin{equation}\label{phgqcu2}
n_vv=\sum_{\{\alpha\in E^{\st} \mid s(\alpha)=v \}} r(\alpha),
\end{equation}
for every $v\in E^0\backslash \sink$, where $n_v=
\max\{w(\alpha)\mid \alpha\in E^{\st}, \, s( \alpha)=v\}$ as in Theorem~\ref{kspcts}. Arrange
$E^0$ in a fixed order such that the elements of $E^0\backslash
\sink$ appear first in the list follow with elements of $\sink$. The
relations ~(\ref{phgqcu2}) can be then written as $N^t \overline v_i= I_w
\overline v_i$, where  $v_i \in E^0\backslash \sink $ and
$\overline v_i$ is the $(0,\dots,1,0,\dots)$ with $1$ in the $i$-th
component.   Therefore,
\[M_E\cong \frac{F}{ \langle N^t \overline v_i = I_w \overline v_i , v_i \in E^0 \backslash \sink \rangle},\] where $F$ is the free abelian monoid generated by the vertices of $E$.
By Theorem~\ref{kspcts} there is a natural monoid isomorphism $\V(\LL_K(E,w))\cong M_E$. So using~(\ref{monio}) we have,
\begin{equation}\label{pajd}
K_0(\LL(E,w))\cong \V(\LL_K(E,w))^{+}\cong M_E^{+}\cong \frac{F^+}{ \langle (N^t-I_w) \overline v_i, v_i \in E^0 \backslash \sink \rangle}.
\end{equation}
Now $F^{+}$ is $\mathbb Z^{E^0}$ and it is easy to see that the
denominator in ~(\ref{pajd}) is the image of $N^t-I_w:\mathbb
Z^{E^0\backslash \sink} \longrightarrow \mathbb Z^{E^0}$.
\end{proof}

\begin{example}
Consider the following weighted graph where all the edges have weight one except the top right edge which has weight two.
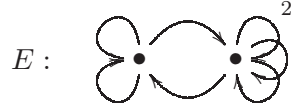
\begin{equation*}
\xymatrix{
 E: &  \bullet \ar@(u,l) \ar@(d,l) \ar@/^1.2pc/[r] & \ar@/^1.2pc/[l] \bullet \ar@(u,r)^{2} \ar@(ur,dr) \ar@(r,d) &
}
\end{equation*}

\bigskip

Then we have
$N=\left(
\begin{matrix}
2 & 1\\
1 & 3
\end{matrix}
\right)
\text{  and  }
I_w=\left(
\begin{matrix}
1 & 0\\
0 & 2
\end{matrix}
\right).$
Therefore the matrix $N-I_w$ is equivalent to
$\left(
\begin{matrix}
0 & 0\\
0 & 1
\end{matrix}
\right)
$. Thus by Theorem~\ref{wke}, $K_0(\LL(E,w))\cong \mathbb Z$.
\end{example}

\begin{example}
Consider the following weighted graph consisted of an acyclic graph (of weight one) and the Leavitt algebra $V_{n,n+k}$ attached to the only sink of the acyclic graph.
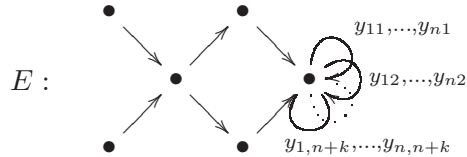
\begin{equation*}
\xymatrix@=13pt{
& \bullet \ar[dr] & & \bullet \ar[dr] & \\
E: & & \bullet \ar[ur] \ar[dr] && \bullet  \ar@(u,r)^{y_{11},\dots,y_{n1}} \ar@(ur,dr)^{y_{12},\dots,y_{n2}}  \ar@{.}@(r,d) \ar@(dr,dl)^{\qquad \qquad y_{1,n+k},\dots,y_{n,n+k}}&\\
& \bullet \ar[ur] & & \bullet \ar[ur]&
 }
\end{equation*}

By Lemma~\ref{tensoir5}, $\LL(E,w)\cong \M_9(V_{n,n+k})$. Therefore by Morita theory and  Theorem~\ref{wke},
\[K_0(\LL(E,w))\cong K_0(V_{n,n+k})\cong \mathbb Z/k\mathbb Z.\]
\end{example}

\begin{example}
Here is an example of a ring $A$ such that $K_0(A)\not=0$ but $[A]=0$. Let $R=V_{n,n+k}$, where $k\geq 2$. By Theorem~\ref{wke}, $K_0(R)=\mathbb Z/k\mathbb Z$ and $k[R]=0$. The ring $A=\M_k(R)$ is Morita equivalent to $R$ (using the assignment 
$P_{\M_k(R)} \mapsto P\otimes_{\M_k(R)} R^k$). Thus $K_0(A)\cong K_0(R)\cong   \mathbb Z/k\mathbb Z$. Under this assignment, $[A]=[\M_k(R)]$ is sent to $k[R]$ which is zero, thus $[A]=0$.
\end{example}

\end{document}